\documentstyle [12pt] {article}

\newenvironment{namelist}[1]{%
\begin{list}{}
    {
     
     \settowidth{\labelwidth}{#1}
     \setlength{\leftmargin}{1.1\labelwidth}
    }
  }{%
\end{list}}

\newcommand{\mapright}[1] {\stackrel{#1} {\hbox to 15pt{\rightarrowfill}}}
\newcommand{\maprightu}[1] {\stackrel{#1} {\hbox to 40pt{\rightarrowfill}}}
\newcommand{\mapdownr}[1] {\Big\downarrow
  \rlap{$\vcenter{\hbox{$\scriptstyle#1$}}$}}
\newcommand{\mapvert}[1] {\Big\vert
  \rlap{$\vcenter{\hbox{$\scriptstyle#1$}}$}}
\newcommand{\mapupr}[1] {\Big\uparrow
  \rlap{$\vcenter{\hbox{$\scriptstyle#1$}}$}}

\newcommand{\mapsear}[1] {\Large\searrow
  \rlap{$\vcenter{\hbox{$\scriptstyle#1$}}$}}
\newcommand{\mapswr}[1] {\Large\swarrow
  \rlap{$\vcenter{\hbox{$\scriptstyle#1$}}$}}
\newcommand{\mapnear}[1] {\Large\nearrow
  \rlap{$\vcenter{\hbox{$\scriptstyle#1$}}$}}

\newcommand{\CC}{\Sym{\cal C}}
\newcommand{\PP}{\Sym{\cal P}}
\newcommand{\DD}{\Sym{\cal D}}

\newcommand{\OO}{\Sym{\cal O}}

\newcommand{\WW}{\Sym{\cal W}}
\newcommand{\SS}{\Sym{\cal S}}
\newcommand{\HH}{\Sym{\cal H}}

\newcommand{\Sym}[1]{\ifmmode{#1}\else \mbox{${#1}$}\fi}
\newcommand{\Vjp}{\bigbreak}
\newcommand{\ov}{\overline}
\newcommand{\va}{\varepsilon}
\newcommand{\wt}{\widetilde}
\newcommand{\wh}{\widehat}
\newcommand{\Fix}{\mbox{Fix}\,}

\newcommand\setC{\mbox{C\hskip-4.5pt\vrule height 6.7pt depth0pt \hskip 4.1pt}}

\newcommand\setR{\mbox{\vrule height 7.2pt depth 0pt \hskip-0.8pt R}}
\newcommand\setP{\mbox{\vrule height 7.2pt depth 0pt \hskip-0.8pt P}}
\newcommand\setZ{\bf Z}

\begin{document}

\baselineskip 18pt

\centerline{\LARGE \bf Feynman Integral, Knot Invariant}

\centerline{\LARGE \bf and Degree Theory of Maps}

\vskip 1cm
\centerline{\large \bf Su-Win Yang}

\vskip 1cm
\centerline{\bf
Department of Mathematics}

\centerline{\bf
National Taiwan University}

\centerline{\bf
Taipei, Taiwan}

\centerline{\bf
E-mail: swyang@math.ntu.edu.tw}

\vskip 12pt
\centerline{September 10, 1997}

\vskip 1cm
\Vjp
\begin{abstract}
The universal Vassiliev invariant from the perturbative Chern-Simons theory is
actually a knot invariant without any correction term.
The anomaly considered by Bott and Taubes is proved
to be zero.
\end{abstract}
\newpage

\noindent
{\Large \bf Contents}

\begin{description}
  \item[0] {\bf Introduction}
  \item[1] {\bf Degree theory and knot invariant}
\begin{description}
  \item[1.1] Knot graph and Configuration space
\begin{description}
    \item[1.1.1] Knot graph
    \item[1.1.2] Equivalence of knot graphs
    \item[1.1.3] Configuration space of knot graphs
\end{description}
  \item[1.2] Boundary of Configuration Spaces
\begin{description}
    \item[1.2.1] Descriptions of Boundaries
    \item[1.2.2] The canonical map $\Phi:C(\Gamma ; A) \longrightarrow B_k$
\end{description}
  \item[1.3] Identification Rules
\begin{description}
    \item[1.3.1] Free vertice
    \item[1.3.2] Univalent inner vertice
    \item[1.3.3] Bi-valent inner vertice
    \item[1.3.4] Edge Identification (Type III)
    \item[1.3.5] Non-edge Identification (Type IV)
    \item[1.3.6] Summary of the identification maps
\end{description}
  \item[1.4] Stable configuration space
\begin{description}
    \item[1.4.1] Identifications
    \item[1.4.2] Combined configuration spaces
\end{description}
\end{description}
\item[2] {\bf $P(k)$-structure}
\begin{description}
  \item[2.1] The canonical $P(k)$-structures
  \item[2.2] Stable $P(r)$-extension of configuration space
\end{description}
\item[3] {\bf Infinitesimal knot graph}
\begin{description}
  \item[3.1] The translation and dilation relation
  \item[3.2] The combined spaces $\WW_0(n)$
  \item[3.3] The $S^1$-action on $\WW_0(n)$
  \item[3.4] The space of infinitesimal knot graphs
\end{description}
\item[4] {\bf $S^1$-equivariant $P(k)$-structure}
\begin{description}
  \item[4.1] The $S^1$-action on $P(k)$
  \item[4.2] $P(k)$-space with a semifree $S^1$-action
  \item[4.3] Preliminary to the proof of Theorem (4.2)
  \item[4.4] Proof of Theorem (4.2)
  \item[4.5] The space $S(k)-F$
\end{description}
\item[5] {\bf Normal bundle of the fixed point set in $\WW_0$}
\begin{description}
  \item[5.1] Inhomogeneous coordinate of $P(k)$
  \item[5.2] Boundary of $W_0(\Gamma)$
  \item[5.3] Identification maps of vector bundles
\begin{description}
    \item[5.3.1] Identification map of type II
    \item[5.3.2] Identification map of type III
    \item[5.3.3] Identification map of type IV
    \item[5.3.4] Identification of type I
\end{description}
  \item[5.4] Conclusion, Over all conclusion
\end{description}
\item[6] {\bf Finite structure group}
\begin{description}
  \item[6.1] Finite structure group
  \item[6.2] A trivialization of $D(\Gamma)$
  \item[6.3] Boundary behavior of the trivialization
\begin{description}
    \item[6.3.1] Having base points
    \item[6.3.2] Inner vertices only
    \item[6.3.3] (Continued)
\end{description}
  \item[6.4] Isotopy of trivializations
  \item[6.5] Transition maps of vector bundle $\DD(n)$
\begin{description}
    \item[6.5.1] Identification of Type 0
    \item[6.5.2] Identification of Type I
    \item[6.5.3] Identification of Type II
    \item[6.5.4] Identification of Type III
    \item[6.5.5] Identification of Type IV
    \item[6.5.6] Conclusion
\end{description}
\end{description}
\end{description}
\newpage

\noindent{\Large \bf \S 0  \ Introduction}

\vskip 6pt
Following the works of Bar-Natan [3], Bott-Taubes [5]
and Altschuler-Freidel [1],
we know how to use Feynman diagrams and their associated
integral to get the beautiful and natural Universal Vassiliev Invariant
in the infinite dimensional algebra of chord diagram
(or Feynman diagram).
But, there still have some defect.
A correction term coming from the
integrals over the spaces of totally concentrated
Feynman diagrams should be considered. The main purpose of this article
is to show that the correction term is equal to zero.

\vskip 20pt
\Vjp
\noindent{\Large \bf  0.1}

\vskip 6pt
Roughly speaking, a Feynman diagram is a type of graph in $\setR^3$ with partial
vertices staying on a knot $K$. If $\Gamma$ denotes a Feynman diagram,
the Feynman integral $I(\Gamma, K)$ is a measure for the
configuration space $C(\Gamma, K)$ of graphs, equiv. to $\Gamma$, on $K$.
$|\Gamma|$ denotes the number of automorphisms of $\Gamma$.
($|\Gamma|$ is also equal to the multiplicity of $\Gamma$ appearing in
$C(\Gamma, K)$.)
The Universal Vassiliev invariants, proposed  by Bar-Natan [3], also
by Kontsevich [8], can be written as

\begin{eqnarray*}
&&Z(K)=\sum\limits_{\Gamma}
\frac{I(\Gamma, K)}{|\Gamma|} [\Gamma], \\
&&\mbox{(the summation is over equivalence classes of Feynman diagrams)}
\end{eqnarray*}
where $K$ is a knot and $[\Gamma]$
denotes the corresponding element
of $\Gamma$ in the algebra of (chord) diagram.
Atschuler and Freidel [1] consider the following correction term:
\[
\alpha =\frac{1}{2} \sum\limits_{\Gamma}
\frac{f_\Gamma}{|\Gamma|} [\Gamma],
\]
which is quite similar to the above $Z(K)$ except
the integral part $f_\Gamma$.
$f_\Gamma$, invented by Bott and Taubes [5],
is a measure for the ``universal space''
$W(\Gamma)$ of totally collapsed Feynman diagrams which
is equivalent to $\Gamma$.
The ``universal space'' $W(\Gamma)$ has the same dimension as
$C(\Gamma, K)$, and $f_\Gamma$ is also a Feynman integral similar to
$I(\Gamma, K)$.
The interesting thing is that the correction term $\alpha$ is
independent of $K$ and what are the values of $\alpha$
and $f_\Gamma$, for each $\Gamma$.
Altschuler and Freidel showed that, when $\Gamma$
is not connected as a graph or the order of $\Gamma$ is even,
$f_\Gamma$ is always zero.

\vskip 6pt
In this article, we propose a method to show that the
sequence
$\alpha=\frac{1}{2} \sum\limits_{\Gamma}
 \frac{f_\Gamma}{|\Gamma|}[\Gamma]$
is equal to zero in the algebra of diagrams.
A weight system $\omega$ is a integral value function for the
diagrams such that
$\sum\limits_{\Gamma} \omega (\Gamma) C(\Gamma, K)$
forms a ``cycle in the homology theory''.
Thus, to show that $\alpha=0$ in the algebra of diagram,
it is enough to show that $\omega(\alpha)=
\frac{1}{2} \sum\limits_{\Gamma} \frac{f_\Gamma}{|\Gamma|}
\omega(\Gamma)$ is zero, for any weight system $\omega$.

\vskip 6pt
As above, we use $W(\Gamma)$ denote the space of totally
collapsed Feynman diagrams equivalent to $\Gamma$.

\vskip 6pt
By the result of Bott and Taubes,
$\sum\limits_{\Gamma} \omega(\Gamma)C(\Gamma, K)$
can not form a cycle
($W(\Gamma)$ is related to the anomalous boundary of
$C(\Gamma, K)$). But, when we consider the analogue for
$W(\Gamma)$, $\sum\limits_{\Gamma}\omega(\Gamma)W(\Gamma)$
does form a cycle under suitable interpretation.
And the value $\omega(\alpha)=\frac{1}{2} \Sigma
\frac{f_\Gamma}{|\Gamma|} \omega(\Gamma)$ is equal to the degree
of canonical map from $W(\Gamma)$ to a special space
$B_k$, which looks like the classifying space of diagrams with
$k$ edges.

\vskip 6pt
Furthermore, $W(\Gamma)$ has a natural $SO(3)$-action and
it is a fibre bundle over $S^2$ with fibre $W_0(\Gamma)$,
a $S^1$-space. Thus, $W(\Gamma)=W_0(\Gamma) $
$\mathrel {\mathop\times\limits_{S^1}}$ $SO(3)$,
and the canonical map
$\Phi_0:W_0(\Gamma)\longrightarrow B_k$
determines
$\Phi:W(\Gamma)\longrightarrow B_k$ completely.

\vskip 6pt
Now, we shall move the problem to $W_0(\Gamma)$ and
$\Sigma \omega(\Gamma)W_0(\Gamma)$.
The $S^1$-action on
$W_0(\Gamma)$ is semifree, that is, the isotropic groups are trivial
or the whole $S^1$.
Let $H(\Gamma)$ denote the fixed point set of the
$S^1$-action on $W_0(\Gamma)$ and $D(\Gamma)$ be the
normal vector bundle of $H(\Gamma)$ in $W_0(\Gamma)$.
Then $\Sigma \omega (\Gamma)D(\Gamma)$ also forms a
vector bundle over $\Sigma \omega(\Gamma)H(\Gamma)$ and
is the normal bundle of $\Sigma\omega(\Gamma)H(\Gamma)$
in $\Sigma\omega(\Gamma)W_0(\Gamma)$.
And, the interesting facts are the following:
(i) the ``degree'' of $\Sigma\omega(\Gamma)W(\Gamma)$
is equal to a Chern number of the normal bundle
$\Sigma\omega(\Gamma)D(\Gamma)$,
(ii) all the possible normal vector bundles has structure group
isomorphic to the semi-direct product of
$\bigoplus\limits_k \setZ_2$ and the symmetry group $\Sigma_k$.
This concludes that the possible degrees are all zero.

\medskip
{\bf Remark :}\ \
In the above consideration, we may restrict the diagrams $\Gamma$
with the same order $n$.
But, under the order restriction,
$W(\Gamma)$ may also have different dimensions.
Thus, we should consider the stable one
$W(\Gamma)^{(r)}=W(\Gamma) \times (\prod\limits_{r}\setR\setP^2)$,
and arrange them into spaces with the same dimension.

\medskip
{\bf Acknowledgement:}\ \
The author wishes to thank Ping-Zen Ong for discussing with him most of the
topics in this article.

\vskip 6pt
In the remaining sections of the introduction,
we talk about the tools and strategies using in the article.

\vskip 20pt
\Vjp
\noindent{\Large \bf \ 0.2 \
Knot graph and its configuration space:}

\vskip 6pt
Suppose $K$ is a oriented knot or a oriented simple curve.
A knot graph $\Gamma$ on $K$ is a finite $1$-dimensional
simplicial complex  with two kinds of vertices in
$V_0(\Gamma)$ and $V_1(\Gamma)$, elements in $V_0(\Gamma)$
are called base points which are always stay in $K$,
and elements in $V_1(\Gamma)$ are called inner vertices which are only
assumed to be points of $\setR^3$ (could be in $K$).
$V(\Gamma)=V_0(\Gamma) \cup V_1(\Gamma)$
is total vertice set.
The edge set $E(\Gamma)$ is a finite set of unordered pairs of
distinct points in $V(\Gamma)$, that is, for
$e \in E(\Gamma)$, $e=\{v, w\}$
with $v \ne w$ in $V(\Gamma)$.
For convenience, choose a linear on $K$ such that the
increasing order is the same as the given orientation.

\vskip 6pt
Suppose $\Gamma_1$ is a knot graph on $K_1$ and $\Gamma_2$
is a knot graph on $K_2$.
A simplicial map $g:\Gamma_1 \longrightarrow \Gamma_2$
is said to be an equivalence of knot graphs, if
$g$ is an isomorphism of simplicial complexes,
$g(V_0(\Gamma_1))=V_0(\Gamma_2)$,
$g(V_1(\Gamma_1))=V_1(\Gamma_2)$,
$g(E(\Gamma_1))=E(\Gamma_2)$, and $g$ preserves
the linear orders of $V_0(\Gamma_1)$ and $V_0(\Gamma_2)$,
which are inherited from the linear orders of $K_1$ and $K_2$,
respectively.

\vskip 6pt
Let $\ov C(\Gamma, K)=\{\Gamma':\Gamma'$
is a knot graph on $K$ and $\Gamma'$ is equivalent to
$\Gamma$ by an equivalence $g\}$,
it is called the configuration space of $\Gamma$ (over $K$).

\vskip 6pt
The interesting thing for knot graph and its configuration space
is to find a measure for the configuration spaces, i.e. integral
on $\ov C(\Gamma, K)$ and a method to collect the knot
graphs $\Gamma_1, \Gamma_2, \cdots$ such that their measure
$\sum\limits_i I(\ov C(\Gamma_i, K))$
constitutes a  knot invariant.

\vskip 6pt
For convenience, consider $C(\Gamma, K)=\{ (g, \Gamma'):\Gamma'$
is a knot graph on $K$ and $g:\Gamma \longrightarrow \Gamma'$
is an equivalence $\}$, it is a covering space of $\ov C(\Gamma, K)$.
Sometimes, for simplicity, just denote it by $C(\Gamma)$.

\vskip 6pt
Whenever $\Gamma$ is a knot graph, there is a natural map
$\phi_\Gamma : E(\Gamma) \longrightarrow \setR \setP^2$
$=\{$ lines in $\setR^3$ through the origin $\}$, defined by:
$e=\{v,w\} \in E(\Gamma)$, $\phi_\Gamma(e)$ is the line passing through
$(v-w)$ and the origin.
If $\Gamma$ has $k$ edges, then the image of
$\phi_\Gamma$ determines an element in
$\prod\limits_k \setR\setP^2 \big / \Sigma_k$, the $k$-fold
symmetric product of $\setR\setP^2$.
Let $B_k=\prod\limits_k \setR\setP^2 \big / \Sigma_k$,
it looks like a classifying space of knot graphs.
Then, there is a canonical smooth map $\Phi:C(\Gamma, K) \longrightarrow
B_k$, $\Phi(\Gamma')=Im(\phi_{\Gamma'})$.
(of course, also for $\ov C(\Gamma, K)$).
$C(\Gamma, K)$ is not compact. But, it is very fortunate that there is
a compactification (Fulton-MacPherson [6]) such that
$\Phi$ can be smoothly extended to this compactification
and the codimension $1$ boundary is the union of
$C(\Gamma; A)$, where $A$ run over all the subsets of $V(\Gamma)$,
containing at least two points.
(details see section 1.2)

\medskip
{\bf Definition :} \ \
(i) $B_k=\prod\limits_k \setR\setP^2 \big / \Sigma_k$,
$\Sigma_k$, the symmetry group acting on
$\prod\limits_k \setR\setP^2$ by permuting the coordinates.

(ii) $(X, f)$ is said to be a $B_k$-space, if $X$ is space and
$f:X \longrightarrow B_k$ is continuous.

(iii) $P(k)=\prod\limits_k \setR\setP^2$.

(iv) $(X,f)$ is said to be a $P(k)$-space, if
$X$ is a space and  $f:X \longrightarrow P(k)$ is continuous.

\medskip
{\bf Examples:} \ \ $C(\Gamma)$ has a canonical $B_k$-structure
$\Phi:C(\Gamma) \longrightarrow B_k$.
When we choose an order for the edges in $\Gamma$,
we have a $P(k)$-structure $\Psi:C(\Gamma) \longrightarrow P(k)$
in the obvious way, but there are $k!$ different $P(k)$-structures
dependent on the $k!$ different choices of orders.
$\ov C (\Gamma)$ also has a canonical $B_k$-structure,
but may not have a $P(k)$-structure. The boundaries
$C(\Gamma; A)$ are similar to the situation for $C(\Gamma)$.

\vskip 6pt
To get an a priori measure for $(C(\Gamma), \Phi)$, we need
to introduce the line bundles for the $B_k$-apces and the
associated local coefficient for these $B_k$-spaces.

\vskip 6pt
Let $L'_k$ be the orientation line bundle over $P(k)$
(thus $L'_k$ is non-orientable) and $L_k$ be the line bundle
over $B_k$ ``push down'' from $L'_k$, that is,
$L_k=L'_k \big / \Sigma_k$.

\vskip 6pt
For any $B_k$-space $(X,f)$, let $L(X)$ denote the
pull-back $f^{\#}(L_k)$ of $L_k$ by $f$.
$(X, f)$ and $(X', f')$ are two $B_k$-space,
$\tau:X \longrightarrow X'$ is said to be a
$B_k$-map, if
$f' \circ \tau = f$. A $B_k$-map $\tau:X \longrightarrow X'$
determines uniquely a bundle map $\ov \tau:L(X)
\longrightarrow L(X')$ such that $\ov f' \circ \ov \tau=\ov f$,
where $\ov f: L(X) \longrightarrow L_k$ and
$\ov f': L(X') \longrightarrow L_k$ are the bundle
maps covering $f$ and $f'$ respectively.

\medskip
{\bf Definition :} \ \
Suppose $X$ is a manifold.
A $B_k$-orientation for a $B_k$-manifold $(X, f)$ is an
orientation for the total space of the line bundle $L(X)$.

\vskip 6pt
Let $\OO_k$ denote the orientation sheaf of $L_k$,
then $f^{\#} \OO_k$ is the orientation sheaf of
$L(X)$ over $X$.
$H_{2k}(B_k, \OO_k) \simeq \setZ$ and
$f_* : H_{2k}(X, f^{\#} \OO_k)$
$\longrightarrow$ $H_{2k}(B_k, \OO_k)$
provides a most natural measure for the $2k$-dimensional cycles
in $X$.

\vskip 6pt
In Chapter 1, we shall use the configuration spaces constructing
a big ``knot invariant'' $B_k$-space, and the
weight systems of Vassiliev invariant will provide
``tautological'' cycles in this $B_k$--space and the
associated degrees from the homology theory are the
knot invariants.

\vskip 6pt
The ``troubles'' of ``knot invariant'' $B_k$-space
will be solved by introducing another
big $B_k$-space $\WW$ whose degrees are always trivial.
($B_k$-space $\WW$ is defined in Chapter 3.)

\Vjp
\vskip 10pt
\noindent{\Large \bf  0.3 \
Construction of ``knot invariant'' $B_k$-space:}

\vskip 6pt
The main tool of this construction is the identification
map for the boundary $C(\Gamma; A)$ of configuration space
$C(\Gamma)$.

\vskip 6pt
For any graph $\Gamma$ and a subset $A$ of vertice set $V(\Gamma)$,
$A(\Gamma)$ denotes the subgraph of $\Gamma$ inside $A$.
Precisely, the vertice set of $A(\Gamma)$ is $A$ and the edge
set $E(A(\Gamma))$ is equal to $\{e=\{v,w\} \in E(\Gamma):v, w \in A\}$.
The codimension $1$ boundary $C(\Gamma; A)$ is a fibre bundle over
$C(\Gamma \big / A)$ with fibre $D(A(\Gamma), x)$.(details see section 1.2
and 1.3)

\vskip 6pt
The identification maps for $C(\Gamma; A)$ have two kinds:
\begin{namelist}{(ii)}
\item[{\rm (i)}]
First kind: \ \ When $|A|\ge 3$ and $A(\Gamma)$ has one
of the following three type of vertice:
free vertice, univalent inner vertice or bivalent inner vertice,
there is a smooth map:
\[
\tau:C(\Gamma; A)\longrightarrow C(\Gamma; A).
\]
\item[{\rm (ii)}]
Second kind: \ \ When $|A|= 2$ (that is, $A$ has exactly two vertices),
there is a smooth map:
\[
\tau:C(\Gamma; A)\longrightarrow C(\Gamma\big / A)\times P(A),
\]
where
\[
P(A)=\left\{\begin{array}{ll}
               \setR \setP^2, & \mbox{if $A(\Gamma)$ has one edge,} \\
               \{\mbox{one point}\}, & \mbox{if $A(\Gamma)$ has no edge.}
            \end{array}
     \right.
\]
\end{namelist}

\vskip 6pt
It is easy to see that, when $X$ is a $B_k$-space,
$X^{(r)} = X \times B_r$ is a $B_{k+r}$-space.
And the identification map can be extended directly to
$\tau:C(\Gamma; A)^{(r)} \longrightarrow C(\Gamma;A)^{(r)}$
(first kind) and $\tau:C(\Gamma; A)^{(r)} \longrightarrow
C(\Gamma \big / A)^{(r+e(A))}$ (second kind),
where $e(A)$ denote the number of edge in $A(\Gamma)$.

\medskip
{\bf Definition :}\ \  Suppose $\Gamma$ is a knot graph.
\[
\mbox{order } (\Gamma) = |E(\Gamma)| - |V_1(\Gamma)|.
\]

\vskip 6pt
Now, we restrict to the knot graphs of order $n$.
Let $\ov\SS(n)$ be the disjoint union of $C(\Gamma)^{(3n-k)}$
($k$ is the edge number of $\Gamma$), for all
$\Gamma$ with order $n$ and having no free vertice,
univalent inner vertice and bivalent inner vertice.
Furthermore, let $\SS(n)$ denote the
quotient space of $\ov \SS(n)$ by identifying $\alpha$ with
$\tau (\alpha)$, for all possible identification maps $\tau$
of first kind and second kind, and for all $\alpha$ in the
corresponding boundary $C(\Gamma; A)$.
Because the identification rule to identify some points
in $\ov \SS(n)$ which have the same image in $B_{3n}$
(note: $C(\Gamma)^{(3n-k)}$ is a $B_{3n}$-space),
$\SS(n)$ is also a natural $B_{3n}$-space.
Let $\Phi:\SS(n) \longrightarrow B_{3n}$ also denote the
structure map. Then,
$\Phi_* : H_{6n}(\SS(n), \Phi^{\#} (\OO_{3n}))\longrightarrow
H_{6n}(B_{3n}, \OO_{3n})$ is supposed to give a lot of
knot invariant.
Unfortunately, there is none.
Originally, for any weight system $\omega$ with integral coefficient,
one expects that
$\Sigma \omega (\Gamma)C(\Gamma)$
is a tautological cycle in
$(\SS(n), \Phi^{\#}(\OO_{3n}))$
and $\Phi_* (\Sigma \omega(\Gamma) C(\Gamma))$ will
provide a knot invariant for $K$.
But, Bott and Taubes [5] showed that
$\Sigma \omega (\Gamma)C(\Gamma)$ has extra boundary
$C(\Gamma; A), A=V(\Gamma)$,
the total vertice set, which is not considered
in the above identifications.
We find that
$\Phi_* (\Sigma \omega (\Gamma)C(\Gamma))$
actually could also be a cycle in $(B_{3n}, \OO_{3n})$.
But, it is hard to find further ``tautological'' identifications in
the space $\SS(n)$.
Thus, we propose to consider additional configuration spaces of the
totally collapsed knot graphs $W(\Gamma, K)$ such that
$\Sigma \omega (C(\Gamma, K)+W(\Gamma, K))$
could be a cycle in some $B_{3n}$-space which contains $\SS(n)$.

\vskip 6pt
But, the trouble of $W(\Gamma, K)$ is that
$W(\Gamma, K)$ depends on a domain $D(K)$ in $S^2$.
The boundary of $D(K)$ is the
closed curve of unit tangent of $K$.
Such a domain has two choice $D_1(K)$ and $D_2(K)$,
$D_1(K) \cup D_2(K)=S^2$.
And, $W(\Gamma, K)$ has also two choice:
$W(\Gamma, D_1(K))$ and
$W(\Gamma, D_2(K))$.
If one show that
$\Sigma\omega(\Gamma)(W(\Gamma, D_1(K))+W(\Gamma, D_2(K))=
\Sigma \omega (\Gamma)W(\Gamma, S^2)$
always gets zero value in
$H_{6n}(B_{3n}, \OO_{3n})$  (under $\Phi_*$).
Then, the $\Phi_*$-value of
$\Sigma\omega(\Gamma)(C(\Gamma, K)+W(\Gamma, K))$ for $K$ is a
knot invariant.

\vskip 6pt
(Sorry, the above argument is quite rough, the correct
form of tautological cycle might be much complicated.)

\vskip 6pt
In the following, we write $W(\Gamma)$ for
$W(\Gamma, S^2)$ $=$ $W(\Gamma, D_1(K))$ $\cup$ $W(\Gamma$, $D_2$ $(K))$,
which is independent of $K$.

\vskip 6pt
Now, we apply the same trick for $W(\Gamma)$,
instead of $C(\Gamma)$,
$(\dim C(\Gamma)=\dim W(\Gamma))$.
And, for some technical reason, we consider the
$P(k)$-structure of $W(\Gamma)$, instead of $B_k$-structure.

\Vjp
\vskip 10pt
\noindent{\Large \bf 0.4 \
The $P(3n)$-space $\WW_0(n)$ and $\WW(n)$:}

\vskip 6pt
Instead of a knot, we consider the line $l_0$
of $z$-axis and the configuration space
$C(\Gamma, l_0)$ of knot graphs on $l_0$, which
is equivalent to $\Gamma$.
When we translate a knot graph on $l_0$ in the line
direction or dilate a knot graph,
the knot graph does not change its image in
$P(k)$ or $B_k$. Thus, $W_0(\Gamma)=C(\Gamma, l_0) \big/
\mbox{(translation and dilation relation)}$ is also
a $P(k)$-space. But, now, we should consider
all $P(k)$-structures by interchanging the order of edges,
that is, $\Sigma_k \cdot \Psi$,
where $\Psi:W_0(\Gamma) \longrightarrow P(k)$ is
some canonical $P(k)$-structure.

\vskip 6pt
Let $\ov \WW_0(n)$ be the disjoint union of all
$W_0(\Gamma) \times P(3n-k) \times \Sigma_{3n} \cdot \Psi$
and $\WW_0(n)$ be the quotient space of $\ov \WW_0(n)$ by the
identification maps with a reasonable choice of
$P(3n)$-structure and the extended translation and dilation relations.
(see section 3.1.)

\medskip
{\bf Explanation:}\ \ The extended translation and dilation relation
is defined on $W_0(\Gamma)$, for the splittable knot graph $\Gamma$.
The purpose is to reduce the dimension of $W_0(\Gamma)$ by one or two
such that we may forget the splittable knot graphs in $\WW_0(n)$.
A knot graph is said to be splittable, if it is a union of two
subgraphs whose intersection is empty or exactly a base point.

\Vjp
\noindent{\large \bf \ 0.4.1 \ $S^1$ action on the spaces
$W_0(\Gamma), \WW_0(n)$}

\vskip 6pt
Identify $S^1$ as the $SO(2)$-subgroup of $SO(3)$,
which contains the rotations around the $z$-axis.
Then $S^1$ naturally acts on
$C(\Gamma , l_0)$, and hence on $W_0(\Gamma )$.
The identification maps in the boundaries
of $W_0(\Gamma )$ are all $S^1$-equivariant.
Thus, the $S^1$-action is also well-defined on
$\WW_0(n)$.

\vskip 6pt
$l_0$ is oriented by the direction $(0,0,1)$.
For any $x \in S^2$, we have the line $l_x=\{tx:t \in \setR\}$
with orientation $x$.
We may define $W_x(\Gamma )$ and $\WW_x(n)$ in the
same way.
Then, $W(\Gamma)=\bigcup\limits_{x \in S^2} W_x(\Gamma)$
with suitable topology, it is a fibre bundle over $S^2$
with fibre $W_0(\Gamma)$, and
$\WW(n)=\bigcup\limits_{x \in S^2} \WW_x(n)$, also.
Because $SO(3)$ acts on $S^2$ transitively,
we may write $W(\Gamma)=SO(3)$
$\mathrel {\mathop\times\limits_{S^1}}$
$W_0(\Gamma)$ and
$\WW(n)=SO(3)$
$\mathrel {\mathop\times\limits_{S^1}}$
$\WW_0(n)$, where the topology of
$W(\Gamma)$ and $\WW(n)$ is obvious, instead of the loose union.
Actually, $SO(3)$ has a natural action on
$\bigcup\limits_{x \in S^2} W_x(\Gamma)$,
which is the same as the $SO(3)$-action
from
$SO(3)$
$\mathrel {\mathop\times\limits_{S^1}}$
$W_0(\Gamma)$ (also for $\WW(n)$).
$SO(3)$ also acts on $\setR\setP^2$ and $P(k)$.
The $P(k)$-structures of $W(\Gamma)$ and $\WW(n)$ are also $SO(3)$-equivariant.
Thus, one might expect that the $S^1$-equivariant
$P(3n)$-structure of $\WW_0(n)$ could control the
$SO(3)$-equivariant $P(3n)$-structure of $\WW(n)$ nicely.

\Vjp
\noindent{\large \bf  0.4.2 \  Orbit space of $S^1$-action:}

\vskip 6pt
Consider the orbit spaces $\WW_0(n) \big / S^1$ and
$P(3n) \big / S^1$ and the map $\ov \Psi :\WW_0(n)\big / S^1$
$\longrightarrow$ $P(3n) \big / S^1$ induced by
$P(3n)$-structure $\Psi:\WW_0(n)\longrightarrow P(3n)$.
We find that the homomorphism
$\ov \Psi_* : H_{6n-3}(\WW_0(n) \big / S^1)$
$\longrightarrow$ $H_{6n-3} (P(3n) \big / S^1)$ has
a good relation with the degree homomorphisms from
$H_{6n}(\WW(n))$ to $H_{6n}(P(3n))$,
(all the homologies are with suitable local coefficients).
Even better, when we restrict the map
$\ov \Psi$ to any neighborhood of the fixed point
set of $\WW_0(n)$ and $P(3n)$, we can still control the
original degree homomorphism of $\WW(n)$.

$P(3n)$ has only one fixed point $\xi_0$,
we may choose the canonical Euclidean space
$\setR^{6n}=\setC^{3n}$ as the neighborhood of $\xi_0$
and its inverse image in $\WW_0(n)$ is a $2n$-dimensional complex vector
$\DD(n)$ over the $(2n-2)$-dimensional fixed point set
$\HH(n)$ of the $S^1$-action on $\WW_0(n)$.

\vskip 6pt
Now, consider the restriction
$\buildrel =\over\Psi $ of $\ov \Psi$ to
$\DD(n)-\HH(n) \big / S^1$,
we have a map from the $S^1$-orbit space of the associated spherical bundle to
$\setC^{3n} - \{0\} \big / S^1$, which is homotopically
equivalent to $\setC\setP^{3n-1}$.
Its homology homomorphism is just a stable characteristic class of the
complex vector bundle $\DD(n)$ over $\HH(n)$.

Up to now, we have the spaces in the following commutative diagram:

$$\matrix{
\HH(n) & \hookrightarrow & \DD(n) & \maprightu{\Psi} & \setC^{3n} &\cr\cr
 & & \mapdownr{} & & \mapdownr{} & \cr\cr
 & & \WW_0(n) & \maprightu{\Psi} & P(3n) & \cr\cr
 & & \mapdownr{} & \mapnear{\Psi} & & \cr\cr
 & & \WW(n) & & & \cr\cr}
$$

(The vertical maps all are inclusions.)

\vskip 6pt
We also have the orbit spaces which form a commutative diagram as follows:

$$\matrix{
\DD(n)- \HH(n) \big / S^1 & \maprightu{\buildrel =\over\Psi} & \setC^{3n}
-0\big / S^1 & \simeq & \setC\setP^{3n-1} &\cr\cr
\mapdownr{} & & \mapdownr{} & & & \cr\cr
\DD(n) \big / S^1 & \maprightu{} & \setC^{3n} \big / S^1 &  &  & \cr\cr
\mapdownr{} & & \mapdownr{} & & & \cr\cr
\WW_0(n) \big / S^1 & \maprightu{\over\Psi} & P(3n) \big / S^1 & & &
\cr\cr}
$$

Now, we can state the Theorems which are needed for the proof of our
main result.

\medskip
{\bf Theorem A:}\ \  If
$\buildrel =\over\Psi_* $ $:$ $H_{6n-4} (\DD(n) - \HH(n) \big / S^1)$
$\longrightarrow$ $H_{6n-4}(\setC\setP^{3n-1})$ is a zero-homomorphism,
then the degree homomorphism
$\Psi_*:H_{6n}(\WW(n))$ $\longrightarrow$ $H_{6n}(P(3n))$
(with suitable local coefficient) is also a zero homomorphism.

\medskip
{\bf Theorem B:}\ \  The stable characteristic classes of
$\DD(n)$ over $\HH(n)$ are all zero. Thus, $\buildrel =\over\Psi_* $
= 0, in any dimensions.

The proof of Theorem A has two steps:
\begin{namelist}{(ii)}
\item[{\rm (i)}]
$\ov \Psi_*:H_{6n-3}(\WW_0(n) \big / S^1) \longrightarrow H_{6n-3}
(P(3n) \big / S^1)$ is $0$-homomorphism $\Longrightarrow$
degree homo. $(\Psi_*)_{6n}$ is $0$.
\item[{\rm (ii)}]
``$(\buildrel =\over \Psi_* )_{6n-4} = 0$''
$\Longrightarrow$
``$(\ov \Psi_*)_{6n-3}=0$''
\end{namelist}

\vskip 6pt
Step (i) is a result of functorial property of $SO(3)$
$\mathrel {\mathop\times\limits_{S^1}} (\cdot)$.

\vskip 6pt
Step (ii) depends on the following

\medskip
{\bf Proposition:}\ \  $H_{6n-3} (P(k)-\xi_0 \big / S^1)=0$
or
$H_{6n-3} (\prod\limits_{3n} S^2 - F \big / S^1)=0$,
where $F$ is the fixed point set of $S^1$-action on
$\prod\limits_{3n} S^2$.

\vskip 6pt
Both results in the proposition can prove Step (ii), we actually
prove the second. (see section 4.4)

\vskip 6pt
And the proof of Theorem B is to show that the vector
bundle $\DD(n)$ over $\HH(n)$ has finite structure group.

\Vjp
\vskip 10pt
\noindent{\Large \bf 0.5 \ Conclusion}

\vskip 6pt
Suppose $\omega$ is a weight system for the Vassiliev invariant of order $n$. Then
$\sum\limits_{\Gamma} \omega (\Gamma) W(\Gamma) \times P(3n-k)
\times \Sigma_{3n} \cdot \Psi$ forms a $6n$-dimensional cycle in
$(\WW(n), \Psi^* (\OO_{3n}))$
and its image in $H_{6n}(P(3n), \OO_{3n})$ is equal
to a non-zero multiple of
$\sum\limits_{\Gamma} \frac{f_{\Gamma}}{|\Gamma|} \omega(\Gamma)$.
Thus, $\sum\limits_{\Gamma} \frac{f_\Gamma}{|\Gamma|} \omega(\Gamma)=0$,
for any weight system $\omega$,
and hence, $\sum\limits_{\Gamma} \frac{f_\Gamma}{|\Gamma|}[\Gamma]=0$
in the algebra of chord diagram.

\newpage

\noindent{\LARGE \bf \S 1. Degree theory and knot invariant}

\vskip 6pt
In this part, we define the knot graphs and their
configuration spaces, discuss the boundaries and see how to
identify the boundaries of different configuration spaces
to get the knot invariants.

\Vjp
\vskip 10pt
\noindent
{\Large \bf  1.1 \ Knot graph and Configuration space}

\vskip 20pt
\Vjp
\noindent
{\large \bf 1.1.1 \ Knot graph}

\vskip 6pt
A graph $\Gamma$ is an abstract $1$-dimensional finite simplicial complex
with vertices in $V(\Gamma)$ and edges in $E(\Gamma)$.
For $e \in E(\Gamma)$, $e$ is an unordered pair of distinct vertices,
that is, $e = \{v, w\}$, $v \ne w$ in $V(\Gamma)$.

\vskip 6pt
$\kappa:S^1 \longrightarrow \setR^3 \mbox{ is an embedding}, \, \,
K=\kappa(S^1)$.
Fix the linear order for the points on $K$ as follows:
\[
\kappa(1) < \kappa(e^{i\theta}) < \kappa(e^{i\theta'}),
\mbox{ for } 0 < \theta < \theta' < 2 \pi.
\]

\medskip
{\bf Definition (1.1):}\ \ A graph $\Gamma$ is said to be a
knot  graph on $K$, if $V(\Gamma) \subset \setR^3$ and $V(\Gamma)$ is
disjoint union of $V_0(\Gamma)$ and $V_1(\Gamma)$ such that $V_0(\Gamma)$ is
a subset of $K$.
The vertices in $V_0(\Gamma)$ are called base points.
The vertices in $V_1(\Gamma)$ are called inner vertices which
are only assumed to be points of $\setR^3$.
Any two vertices are distinct points of $\setR^3$.

\vskip 6pt
$\setR\setP^2$ denotes the space of straight lines in $\setR^3$
passing through the origin, it is also diffeomorphic to $S^2$ quotiented by the
$\setZ_2$--action: $-1 \cdot x =-x$.
Any non-zero vector $x$ of $\setR^3$ determines an
element $[x]=\{ tx, t \in \setR\}$ of $\setR\setP^2$.

\vskip 6pt
The interest of knot graph $\Gamma$ and $\setR\setP^2$ is the following
map for the edges:
$\phi_\Gamma: E(\Gamma) \longrightarrow \setR\setP^2$,
for the $e=\{v, w\} \in E(\Gamma)$,
$\phi_\Gamma (e)=[v-w]$.
If $X$ is a space of knot graphs, we have a canonical map
associated with $X$, $\Phi : X \longrightarrow \{$
finite subsets of $\setR\setP^2\}$,
$\Phi(\Gamma)=\{\phi_\Gamma(e):e\in E(\Gamma)\}$,
for $\Gamma \in X$.
When all the graphs in $X$ have the same number of edges, say $k$,
then $\Phi$ is a well-defined map from $X$ to the
$k$-fold symmetric product of $\setR\setP^2$,
$\prod\limits_{k} \setR\setP^2 \big / \Sigma_k$,
where $\Sigma_k$ is the symmetry group of
$\{1, 2, \cdots, k\}$ and
$\Sigma_k$ acts on
$\prod\limits_{k} \setR\setP^2$ by permuting the coordinates.

\Vjp
\vskip 20pt
\noindent{\large \bf 1.1.2 \ Equivalence of knot graphs}

\vskip 6pt
Suppose $K, K'$ are knots or simple curves
with a natural linear order for their points.

\medskip
{\bf Definition (1.2):}\ \ Assume $\Gamma$ is knot
graph on $K$ and $\Gamma'$ is a knot graph on $K'$.
A simplicial map $g:\Gamma \longrightarrow \Gamma'$
is said to be an equivalence of knot graphs, if $g$ is an
isomorphism of simplicial complexes and the restriction
of $g$ to $V_0(\Gamma)$ is an
order preserving bijection from
$V_0(\Gamma)$ onto $V_0(\Gamma')$.
(Thus, $g$ sends base points to base point and sends
inner vertices to inner vertices.)

\vskip 6pt
For any knot graph $\Gamma$, the base points of $\Gamma$ have a
linear order which is inherited from the linear order of
points in the knot or simple curve.
The equivalence $g$ must preserve the linear orders of base points.
Furthermore, $g$ is a bijection of $E(\Gamma)$ and $E(\Gamma')$.

\vskip 6pt
Two equivalent knot graphs may be not on the same knot.
Thus, it is convenient to consider a equivalence class of
knot graphs or a single knot graph as an abstract
$1$-dimensional simplicial complex with two type of vertices:
one is the base points with a linear order and one is the inner
vertices.

\Vjp
\vskip 10pt
\noindent
{\large \bf 1.1.3 \ Configuration space of knot graphs}

\medskip
{\bf Definition (1.3):}\ \ Suppose $\Gamma$ is
a knot graph (on some knot or simple curve with linear order) and
$K$ is a knot or a simple curve in $\setR^3$.
$C(\Gamma, K)$ denotes the set of all pairs
$(g, \Gamma')$,
$\Gamma'$ is a knot graph on $K$ and $g$ is equivalence from $\Gamma$
to $\Gamma'$.

\vskip 6pt
An alternating definition of configuration space is the
following:
$\ov C(\Gamma, K)$ $=$ $\{\Gamma'$ $:$ $\Gamma'$ is
a knot graph on $K$ and $\Gamma'$ is equivalence to $\Gamma\}$.
It is easy to see that $C(\Gamma, K)$ is a covering space of
$\ov C(\Gamma, K)$, and the number of elements in the
fibre is the
same as the number of automorphisms of $\Gamma$.
$|\Gamma|$ will be used to denote the order of automorphism
group
$\{g$ $:$ $\Gamma$ $\longrightarrow$ $\Gamma$ equivalence $\}$ of $\Gamma$.

\vskip 6pt
If $\Gamma$ has $m$ base points and $s$ inner vertices, then
the dimension of $C(\Gamma, K)$ is $m+3s$.

\vskip 6pt
If $\Gamma$ has $k$ edges
$e_1, e_2, \cdots, e_k$.
Let $B_k=\prod\limits_{k}\setR\setP^2 \big / \Sigma_k$.
$\Phi$ denotes the canonical map from $C(\Gamma, K)$ to
$B_k$, that is, $\Phi$ $(g$, $\Gamma')$ $=$ $Im$ $(\phi_{\Gamma'})$.
Ordinarily, $\Phi(C(\Gamma, K))$
also has dimension $m+3s$ in $B_k$.
But, sometimes, it is degenerate. For example, when $K$ is
straight line in $\setR^3$,
$\dim \Phi(C(\Gamma, K)) \le m+3s -2$.
The reason is that when a knot graph on a straight line translates
on the line or dilates all give the same image in $B_k$.
In this article, the knot invariant will be a correct
measure of $\Phi(C(\Gamma_i, K))$ or
$\Phi(\ov C(\Gamma_i, K))$,
$i=1, 2, \cdots, l$, for some
$\Gamma_i$, $i=1, 2, \cdots, l$.
Thus, one of important work is to find out some principles or
some rules when or where there are two distinct knot graphs or a
bunch of knot graphs which have the same image in $B_k$.
Then, we can make the identifications at first and use them or disuse
them for the purpose of knot invariant.

\medskip
{\bf Example (1.4):}\ \
Suppose $\Gamma$ has a univaliant inner vertice $v$,
that is, some edge
$e=\{v,w\}$ is only one edge containing $v$ as the end point.
Then, each $(g, g(\Gamma))$ in
$C(\Gamma, K)$ has the same image in $B_k$ as the bunch of elements in
$\{g_t, g_t(\Gamma)), t \ne 0\}$,
where $g_t(v)=g(w)+t(g(v)-g(w))$
and
$g_t(v')=g(v')$, for $v'\ne v$.
Thus, we may identify
$\{g_t, g_t(\Gamma))\}$ to a point for each
$(g, g(\Gamma))$ in $C(\Gamma, K))$,
and the new space has less $1$ dimension.
In the following, we shall find more and more identification rules to make up
spaces for knot invariant.

\vskip 6pt
The most important identification rules will happen in the
boundaries of configuration spaces $C(\Gamma, K)$.

\vskip 20pt
\Vjp
\noindent{\Large \bf 1.2.  \
Boundary of Configuration Spaces}

\Vjp
\vskip 10pt
\noindent {\large \bf 1.2.1.  \
Description of Boundaries}

\vskip 6pt
In the remaining part of this chapter, $K$ will be
fixed and $C(\Gamma, K)$ will be denoted simply by $C(\Gamma)$,
for any knot graph $\Gamma$.

\vskip 6pt
$C(\Gamma)$ is, by no mean, compact, but there is a compactification
(Fulton and MacPherson [6])
for the configuration spaces such that
$\Phi:C(\Gamma)\longrightarrow B_k$ can be extended continuously
to the compactification of $C(\Gamma)$.
For convenience, we assume that $C(\Gamma)$ has been substituted
by its Fulton-MacPherson compactification.

\medskip
{\bf (1.5)}\ \ (Fulton-MacPherson)
The codimension $1$ boundary of $C(\Gamma)$ is the union of
$C(\Gamma; A)$, for any $A$ which is a subset of the vertice set $V(\Gamma)$
and contains at least two vertices.
$C(\Gamma; A)$ is a fibre bundle over $C(\Gamma \big / A)$
with fibres $D(A(\Gamma),x)$, where
$\Gamma \big / A$ and $D(A(\Gamma),x)$ are defined in (1.6) and (1.8).

\medskip
{\bf (1.6)}\ \ $\Gamma \big / A$ denotes the knot graph with vertice
set $V(\Gamma \big / A)  = (V(\Gamma)-A) \cup \{a\}$,
$a \not \in V(\Gamma)$, and edge set $E(\Gamma \big / A)=
\big\{e=\{v,w\}:v, w \not \in A$,
$e \in E(\Gamma)\big\}$ $\cup$
$\big\{e=\{v,a\}:v \not \in A$
and there exists $w \in A$ such that $\{v,w\} \in E(\Gamma)\big\}$.
The set of base points in $\Gamma \big / A$ is the following:
\[
V_0(\Gamma \big / A)=\left\{ \begin{array}{ll}
               V_0(\Gamma), & \mbox{if $A\cap V_0(\Gamma)=\phi$}, \\
               \{a\} \cup V_0(\Gamma)-A, & \mbox{if $A\cap V_0(\Gamma)\ne \phi$}.
            \end{array}
     \right.
\]
Because the new vertice $a$ replaces the vertices in $A$, the linear order
on $V_0(\Gamma \big / A)$ can be defined in the obvious way, under
the Assumption(1.7) below.

\medskip
{\bf Assumption (1.7):}\ \
$A \cap V_0 (\Gamma)$ is an interval of the linear order set
$V_0(\Gamma)$, that is,
$\{v\in V_0(\Gamma):v_1 \le v\le v_2\}$, for some
$v_1,v_2$ in $V_0(\Gamma)$.

\medskip
{\bf (1.8)}\ \ Description of $D(A(\Gamma), x)$

\vskip 6pt
Suppose $(g, g(\Gamma \big / A))$ is an element in
$C(\Gamma \big / A)$ and $x$ is the point $g(a)$.
Then $D(A(\Gamma), x)$ is the fibre of $C(\Gamma; A)$
over $(g, g(\Gamma \big / A))$.

\vskip 6pt
(i) $A(\Gamma)$ is the graph with vertice set
$V(A(\Gamma))=A$ and edge set $E(A(\Gamma))$
$=$ $\{e=\{v,w\} \in E(\Gamma):v$ and $w$
are in $A\}$.
When $A \cap V_0(\Gamma)$ is non-empty,
$A(\Gamma)$ is a knot graph with set of base points
$V_0(A(\Gamma))=A \cap V_0(\Gamma)$.
When $A \cap   V_0(\Gamma)$ is empty, $A(\Gamma)$ is a
graph without base point. Thus, whether $A$ contains base points or not,
the situation is quite different.

\vskip 6pt
(ii) (Case 1): \ \ $A \cap  V_0(\Gamma)$ is non-empty.

\vskip 6pt
In this case, $a$ is a base point. Thus, $x=g(a)$ is a point on $K$
and $D(A(\Gamma), x)$ is the set of ``infinitesimal'' knot graphs
concentrated at $x$.
Precisely, let $l(x)$ denote the tangent line of $K$ at $x$,
$D(A(\Gamma), x)$ is the configuration space
$C(A(\Gamma), l(x))$ quotiented by the translation and dilation
relation:
$(g_1, g_1(A(\Gamma)))$ and $(g_2, g_2(A(\Gamma)))$
are two element in $C(A(\Gamma), l(x))$,
$g_1$ $\sim$ $g_2$ if there exist $\lambda >0$ and
$y_0$ in $\setR^3$ such that
$g_1(v)=y_0+\lambda g_2(v)$ for all
$v \in V(A(\Gamma))=A$.
This definition is good for any kind of graph  with or
without base point.
But, in this case, $y_0$ could only be a vector parallel
to $l(x)$.
Thus,
$\dim D(A(\Gamma), x)=$
$\dim C(A(\Gamma), l(x))-2$,
and
$\dim C(\Gamma; A)=$
$\dim C(\Gamma\big / A) +$
$\dim D(A(\Gamma), x)=$
$\dim C(\Gamma)-1$.

\vskip 6pt
(iii) (Case 2)\ \ $A \cap V_0(\Gamma)$ is empty.

\vskip 6pt
In this case, $D(A(\Gamma), x)$, is nothing to do
with $K$ or the tangent line of $K$.
The original meaning of $D(A(\Gamma), x)$ is also
the set of infinitesimal graphs concentrated at $x$.
Precisely, $D(A(\Gamma), x)$ is the configuration space
$C(A(\Gamma))=\{(g, \Gamma'): g: A(\Gamma) \longrightarrow \Gamma'$
is an equivalence $\}$
quotiented by the translation and dilation relation defined in (case 1).
Note: this relation reduces the dimension by $4$, that is,
$\dim D(A(\Gamma), x)=\dim C(A(\Gamma))-4$.
Formally, $D(A(\Gamma), x)$ is independent of $x$, that is,
$D(A(\Gamma), x)=D(A(\Gamma), x')$, for any
$x$, $x'$ in $\setR^3$.
But, similar to case 1,
$\dim C(\Gamma; A)=\dim C(\Gamma \big / A) + \dim D(A(\Gamma), x)$,
it is also equal to
$\dim C(\Gamma)-1$.

\Vjp
\noindent {\large \bf 1.2.2. \
The canonical map $\Phi:C(\Gamma; A) \longrightarrow B_k$}

\vskip 6pt
Suppose $\Gamma$ has $k$ edge. Then
$\Gamma \big / A$ and $A(\Gamma)$ have the sum of edge numbers equal to $k$.
Thus, for $\alpha=(g, g(\Gamma \big / A))$
in $C(\Gamma \big / A)$ and
$\beta=(g', g'(A(\Gamma)))$ in
$D(A(\Gamma), x)$,
$\Phi(\alpha, \beta)=(\Phi(\alpha), \Phi(\beta))$, it is a
well-defined element in $B_k$.

\vskip 6pt
To show that
$\Phi : C(\Gamma; A) \longrightarrow  B_k$ is a continuous limit
of $\Phi:C(\Gamma) \longrightarrow B_k$,
let $\alpha+t\beta$ denote the following element in
$C(\Gamma)$, as $t$ is a small positive number:

\vskip 6pt
Fix an element $a_0$ in $A$ $(a_0 \in A \cap V_0(\Gamma)$,
if $A \cap V_0(\Gamma) \ne \phi)$.
\[
g_t(v)=\left\{\begin{array}{lll}
               g(v), & \mbox{$v\not \in A$,}\\
               g(a), & \mbox{$v=a_0$,} \\
               g(a)+\frac{t}{|g'|}(g'(v)-g'(a_0)), &\mbox{$v \in A.$}
              \end{array}
     \right.
\]
$$(|g'|=\max \{|g'(v)-g'(w)|, \, \, v, w, \in A\}) $$
$\alpha + t\beta=(g_t, g_t(\Gamma))$, and it is easy to see that
\begin{eqnarray*}
&&\lim\limits_{t \to 0} (\alpha + t\beta)=(\alpha, \beta) \mbox{ and} \\
&&\lim\limits_{t \to 0} \Phi (\alpha+t\beta)=(\Phi(\alpha), \Phi(\beta)).
\end{eqnarray*}
In fact, $((g, g') \longrightarrow g_t)$ defines a collar of
$C(\Gamma; A)$ in $C(\Gamma)$.

\Vjp
\vskip 20pt
\noindent
{\Large \bf  1.3. \  Identification Rules}

\vskip 6pt
We like to find all possible principle that which
knot graphs have the same image in $B_k$, and
construct the identification maps.

\Vjp
\vskip 10pt
\noindent
{\large \bf 1.3.1 \ \  Free vertice}

\vskip 6pt
Suppose $A$ has more than two vertices and $A(\Gamma)$ contains
a free vertice $v$, that is, $v$ is not an endpoint of any
edge in $A(\Gamma)$.

\vskip 6pt
Let $A_1 = A-\{v\}$ and we identify an element
$(g, g(A(\Gamma)))$ in
$D(A(\Gamma), x)$ with its restriction to $A_1$,
$(g |_{A_1} , g(A_1(\Gamma)))$ in
$D(A_1(\Gamma), x)$.
We can also choose an embedding $\phi$ from
$D(A_1(\Gamma), x)$ into
$D(A(\Gamma), x)$, which is a right inverse to the restriction,
as follow:
\subitem
$(g_1, g_1(A_1(\Gamma))) \in D(A_1(\Gamma), x)$,  \hfill\break
$\mbox{let }||g_1||=\max\{
 |g_1(v_1)-g_1(v_2)| :v_1, v_2 \in A_1\}$,

\subitem
fix $v_0 \in A_1$, $g: A \longrightarrow \setR^3$
defined by:
$g(v) = g_1(v_0)+(2||g_1||, 0, 0)$,
$g(v')=g_1(v')$, for $v' \ne v$, and
$\phi (g_1, g_1(A_1(\Gamma))=(g, g((A(\Gamma)))$.
(Note: when $|A|=2$, $||g_1||=0$.)

Let $\tau_0:D(A(\Gamma), x)\longrightarrow D(A(\Gamma), x)$
be defined by
$\tau_0(g, g(A(\Gamma)))=\phi(g |_{A_1}, g(A_1(\Gamma)))$.
Thus, $D(A(\Gamma), x)$ is identified as its codimension $1$
subset $\phi(A(A_1(\Gamma), x))$, and so is
$C(\Gamma; A)$.
This identification map $\tau_0:C(\Gamma; A)\longrightarrow C(\Gamma; A)$
shall be called the identification map of type $0$.

\Vjp
\vskip 10pt
\noindent
{\large \bf 1.3.2 \ Univalent inner vertice}

\vskip 6pt
Suppose $A$ has more than two vertices and $A(\Gamma)$
contains a univalent inner vertice $v$, that is,
there exists a unique edge
$e=\{v, w\}$ containing $v$ as its endpoint.
Let $A_1 = A-\{v\}$, then $A_1(\Gamma)$ has
one less edge than $A(\Gamma)$.
We shall identify an element
$(g, g(A(\Gamma))$ in $D(A(\Gamma), x)$ with
$(g|_{A_1}, g(A_1(\Gamma)), \frac{g(v)-g(w)}{|g(v)-g(w)|})$
in $D(A_1(\Gamma), x) \times S^2$. Similar to the above
case in 1.3.1, we may choose an embedding
$\phi_1:D(A_1(\Gamma), x)\times S^2 \longrightarrow D(A(\Gamma), x)$
such that $\phi_1(\alpha)$ is identified with $\alpha$, for
$\alpha$ in $D(A_1(\Gamma), x) \times S^2$.
And the identification map
$\tau_1:D(A(\Gamma), x) \longrightarrow D(A(\Gamma), x)$
is defined as $\tau_1(g, g(A(\Gamma)))=\phi_1(g|_{A_1}, g(A_1(\Gamma))$,
$\frac{g(v)-g(w)}{|g(v)-g(w)|})$.
An explicit form of $\phi_1$ is given in the following:
\subitem
$\beta=(g_1, g_1(A_1(\Gamma)))\in D(A(\Gamma), x)$,
$y \in S^2$, \hfill\break
let
$||g_1||=\max\{|g_1(v_1)-g_1(v_2)|$,
$v_1, v_2$ in $A_1\}$,\hfill\break
$g:A \longrightarrow \setR^3$,
defined by $g(v)=g_1(w)+2||g_1||y$,
$g(v')=g_1(v')$, for $v'\ne v$. \hfill\break
Then, $\phi_1(\beta, y)=(g, g(A(\Gamma)))$, it is an element
in $D(A(\Gamma), x)$.

\noindent
Now, $\tau_1$ is defined on the fibres $D(A(\Gamma), x)$ of
$C(\Gamma; A)$, we extend $\tau_1$ to the whole $C(\Gamma; A)$ in
the obvious way.
Thus, both $D$ $(A$ $(\Gamma)$, $x)$ and $C$ $(\Gamma$; $A)$ are identified by
$\tau_1$ into their proper subsets of codimension $1$.

\vskip 6pt
This identification map $\tau_1:C(\Gamma; A)\longrightarrow C(\Gamma; A)$
shall be called the identification map of type I.

\Vjp
\noindent
{\large \bf 1.3.3 \  Bi-valent inner vertice}

\vskip 6pt
Suppose $A(\Gamma)$ has a bi-valent inner vertice $v$,
that is, $v$ is an inner vertice and there are exactly
two edges $e_1$, $e_2$ which contain $v$ as the endpoint.

\vskip 6pt
Assume $e_1=\{v, w_1\}$ and $e_2=\{v, w_2\}$.
The identification map on
$D(A(\Gamma), x)$ (or $C(\Gamma; A))$ defined by:

\subitem
$(g, g(A(\Gamma))) \in D(A(\Gamma), x)$, \hfill\break
let $\ov g : A \longrightarrow \setR^3$ be:
$\ov g(v)=g(w_1)+g(w_2)-g(v)$,
and $\ov g(v')=g(v')$, for $v' \ne v$.

And, $\tau_2(g, g(A(\Gamma)))=(\ov g, \ov g(A(\Gamma)))$.

\noindent
$\tau_2$ is a $B_k$-map in the sense:
$\Phi \circ \tau_2=\Phi$,
$\Phi:D(A(\Gamma), x)\longrightarrow B_k$ is the
canonical $B_k$-structure.
This identification map $\tau_2:C(\Gamma; A)\longrightarrow C(\Gamma; A)$
shall be called the identification map of type II.

\Vjp
\vskip 10pt
\noindent
{\large \bf 1.3.4 \ Edge Identification (Type III)}

\vskip 6pt
Suppose $A$ is equal to an edge $e=\{v, w\}$ of
$\Gamma$ and $A$ contains at least one inner vertice.

\vskip 6pt
Then, the identification map
$\tau_3:D(A(\Gamma), x)\longrightarrow \setR\setP^2$ is defined
by :
$(g, g(A(\Gamma)))\in D(A(\Gamma), x)$,
$\tau_3:(g, g(A(\Gamma)))=[g(v)-g(w)]$,
and $C(\Gamma; A)$ is identified with
$C(\Gamma \big / A) \times \setR\setP^2$, by
$\tau_3$ also. This identification map
$\tau_3:C(\Gamma; A) \longrightarrow C(\Gamma \big /A) \times \setR\setP^2$
is called the identification map of type III.

\Vjp
\vskip 10pt
\noindent
{\large\bf 1.3.5 \ Non-edge Identification (Type IV)}

\vskip 6pt
Suppose $A$ has exactly two vertices and $A$ is
not an edge in $\Gamma$.

\vskip 6pt
Then, the identification map $\tau_4$ is exactly the bundle
projection from $C(\Gamma; A)$ onto
$C(\Gamma \big / A)$.
Note: $A(\Gamma)$ has no edge, this identification does not
lose anything.

\vskip 6pt
$\dim C(\Gamma; A)=\dim C(\Gamma \big / A)$, if and only if,
$A$ contains only base points.

\vskip 6pt
When $A$ contains at least one inner vertice, it could be thought as
a "special case" of type 0, in spirit.

\Vjp
\vskip 10pt
\noindent
{\large \bf 1.3.6 \  Summary of the identification maps}

\vskip 6pt
\begin{namelist}{(iii)}
\item[{\rm (i)}]
$|A|\ge 3$, $\tau_0, \tau_1, \tau_2:C(\Gamma; A) \longrightarrow
C(\Gamma; A)$.
\item[{\rm (ii)}]
$|A|=2$, $\tau_3, \tau_4:C(\Gamma; A) \longrightarrow
C(\Gamma \big /  A) \times P(A)$,

where
\[
P(A)=\left\{\begin{array}{ll}
               \setR \setP^2, & \mbox{if $A(\Gamma)$ has one edge,} \\
               \{\mbox{one point}\}, & \mbox{if $A(\Gamma)$ has no edge.}
            \end{array}
     \right.
\]
\item[{\rm (iii)}]
$\tau_0, \tau_1$ send $C(\Gamma; A)$ into codimension $1$ subset;
$\tau_4 (C(\Gamma; A))$ has dimension  $\le$ $\dim C(\Gamma)-2$,
if $A \not \subset V_0(\Gamma)$.
\end{namelist}

\vskip 20pt
\Vjp
\noindent
{\Large \bf 1.4 \ Stable configuration space}

\vskip 6pt
For any positive integer $q$, let $C(\Gamma)^{(q)}$ $=$
$C$ $(\Gamma)$  $\times$ $B_q$,
where $B_q$ $=$ $\prod\limits_{q}\setR \setP^2 \big /\Sigma_q$.
An element in $C(\Gamma)^{(q)}$ is a knot graph
$(g$, $\Gamma')$ together with $q$ points
$\xi_1$, $\xi_2$, $\cdots$, $\xi_q$ in
$\setR\setP^2$.
$C(\Gamma)^{(q)}$ is called the $q$-stable extension of $C(\Gamma)$.
Similarly,
$C(\Gamma; A)^{(q)} = C(\Gamma ;A) \times B_q$.

\vskip 10pt
\Vjp
\noindent
{\large\bf 1.4.1 \ Identifications}

\vskip 6pt
Following section 1.3.6, we have the identification maps directly:
\begin{namelist}{(ii)}
\item[{\rm (i)}]
$|A| \ge 3$, $\tau_0, \tau_1, \tau_2:C(\Gamma; A)^{(q)}
\longrightarrow C(\Gamma; A)^q$
\item[{\rm (ii)}]
$|A| =2$, $\tau_3, \tau_4 :C(\Gamma; A)^{(q)}
\longrightarrow C(\Gamma\big / A)^{(q+e(A))}$,
\end{namelist}
where $e(A)$ is the number of edges in $A(\Gamma)$.

\vskip 10pt
\Vjp
\noindent
{\large\bf 1.4.2 \ Combined configuration spaces}

\medskip
{\bf Definition (1.13):}\ \ Suppose $\Gamma$ has no free vertice.
order $(\Gamma)=|E(\Gamma)| - |V_1(\Gamma)|$, that is,
the number of edges in $\Gamma$ minus the number of inner vertices in
$\Gamma$.

\medskip
{\bf Definition (1.14):}\ \
Suppose $\Gamma$ is a knot graph with order $n$ and
$k$ edges.
$S(\Gamma)$ is the $(3n-k)$-stable extension
of $C(\Gamma)$,
i.e. $C(\Gamma)^{(3n-k)}$. (it is a $B_{3n}$-space).

\medskip
{\bf (1.15):}\ \
Let $\ov \SS(n)$ be disjoint union of
$S(\Gamma)$, for all $\Gamma$ with order $n$, and
$\SS(n)$ be the quotient space of $\ov \SS(n)$ by
taking all possible identifications of type 0, I, II, III and IV.
Thus $\SS(n)$ has a canonical map
\[
\Phi:\SS(n) \longrightarrow B_{3n}.
\]
Is there any cycle $\alpha$ in
$\SS(n)$ such that $\Phi_*(\alpha) \ne 0$ in $H_*(B_{3n})$ ?
\newpage

\vskip 20pt
\Vjp
\noindent{\LARGE \bf \S 2 \ $P(k)$-structure}

\vskip 20pt
\Vjp
\noindent{\Large \bf 2.1  \ The canonical $P(k)$-structures}

\vskip 6pt
$P(k)$ denotes the space $\prod\limits_{k} \setR\setP^2$, the
product of $k$ copies of $\setR\setP^2$.

\vskip 6pt
Assume $\Gamma$ has $k$ edges. Choose an order for the
edges $(e_1, e_2, \cdots, e_k)$,
$e_i=\{v_i, w_i\}$, $i=1, 2, \cdots, k$.
Let $\Psi:C(\Gamma) \longrightarrow P(k)$ be the
continuous map defined by: for
$(g, g(\Gamma)) \in C(\Gamma)$, $g(\Gamma)$ is knot
graph on $K$,
$$\Psi (g, g(\Gamma))=([g(v_1)-g(w_1)], [g(v_2)-g(w_2)], \cdots,
      [g(v_k)- g(w_k)]).$$
\noindent Changing the order of edges, we have $k!$ different maps from
$C(\Gamma)$ to $P(k)$.
Let $\Sigma_k \cdot \Psi$ denote the set of $k!$ maps,
it is the orbit of $\Sigma_k$-action. All the maps in
$\Sigma_k \cdot \Psi$ are called the canonical
$P(k)$-structures of $C(\Gamma)$.

\vskip 6pt
Now, consider the boundary $C(\Gamma; A)$, it is a fibre bundle
over $C(\Gamma \big / A)$ with fibre
$D(A(\Gamma), x)$. Assume $\Gamma \big / A$ has
$k_1$ edges and $A(\Gamma)$ has $k_2$ edges,
$k_1+k_2=k$.

$\Psi_1:C(\Gamma \big / A) \longrightarrow P(k_1)$
denotes a canonical $P(k_1)$-structure of $C(\Gamma \big / A)$, and
$\Psi_2:D(A(\Gamma), x)\longrightarrow P(k_2)\,$
denotes a canonical $P(k_2)$-structure of
$D(A(\Gamma), x)$,
then the restriction of $\Psi$ to $C(\Gamma; A)$,
$\Psi |_{C(\Gamma; A)} =\sigma \cdot (\Psi_1, \Psi_2),\,$
for some permutation $\sigma$ in $\Sigma_k$.

\medskip
{\bf Note:} To make sure the equality, $k_1 + k_2 = k$,
the ordered edge set $(e_1, e_2, \cdots, e_{k_1})$ of
$E(\Gamma \big / A)$ may have the same
edge in different coordinates, that is, one edge may appear more one
times in the ordered edge set. And, when the knot graph has multi-edge,
$\Sigma_k \cdot \Psi$ does not have $k!$ different maps. But, there is
no essential influence on the theories of this article.

\Vjp
\vskip 20pt
\noindent{\Large \bf 2.2  \
Stable $P(r)$-extension of configuration space}

\medskip
{\bf Definition (2.1):}\ \ Suppose $\Gamma$ is a knot graph with $k$
edges and $r$ is a positive integer.
$C(\Gamma) \times P(r)$ is called the $P(r)$-extension
of $C(\Gamma)$.
If $A$ is a subset of $V(\Gamma)$,
$D(A(\Gamma), x) \times P(r)$ is called the $P(r)$-extension
of $D(A(\Gamma), x)$.

\vskip 6pt
We are interested in the $P(r)$-extension of $C(\Gamma)$
and their canonical $P(k+r)$-structure
$\Sigma_{k+r} \cdot \Psi$,
where
$\Psi:C(\Gamma) \times P(r) \longrightarrow P(k+r)$
is defined by:
\begin{eqnarray*}
&&(g, g(\Gamma)) \in C(\Gamma), (\xi_1, \xi_2, \cdots, \xi_r) \in
P(r),\\
&&\Psi((g, g(\Gamma)), \xi_1, \xi_2, \cdots, \xi_r)=
(\Psi(g, g(\Gamma)), \xi_1, \xi_2, \cdots, \xi_r).
\end{eqnarray*}
Similarly, for $D(A(\Gamma), x) \times P(r)$.

\vskip 6pt
Now, consider the $P(r)$-extension of $C(\Gamma)$
prescribed a canonical $P(k+r)$-structure, and the disjoint union
of the all possible spaces is denoted by
$C(\Gamma) \times P(r) \times \Sigma_{k+r} \cdot \Psi$,
where
$\Sigma_{k+r} \cdot \Psi$ is a discrete space.

\medskip
{\bf Definition (2.2):}\ \ Suppose $\Gamma$ is a knot graph on $K$,
the order of $\Gamma$ is $n$ and $\Gamma$ does not have
any free vertice.
\[
P(\Gamma)=C(\Gamma)\times P(3n-k) \times \Sigma_{3n} \cdot \Psi,
\]
where $k=|E(\Gamma)|$ and $\Psi$ is a canonical
$P(3n)$-structure of $C(\Gamma) \times P(3n-k)$.

\vskip 6pt
In the following, we consider the identifications
between the spaces $P(\Gamma)$, or more general,
$C(\Gamma) \times P(r) \times \Sigma_{k+r} \cdot \Psi$.
The identification maps in section 1.3 have two forms:
\begin{namelist}{(ii)}
\item[{\rm (i)}] $|A|\ge 3$,
$\tau:C(\Gamma; A) \longrightarrow C(\Gamma; A)$ \hfill\break
(type 0, type I or type II)
\item[{\rm (ii)}] $|A|=2$,
$\tau: C(\Gamma; A) \longrightarrow C(\Gamma \big / A) \times
P(A)$ \hfill\break
(type III or type IV) \hfill\break
($P(A)=P(|E(A(\Gamma))|)$.)
\end{namelist}
All the identification maps directly give identification maps
for $C(\Gamma) \times P(r)$, that is,
\[
\tau \times id :
\left\{
\begin{array}{ll}
C(\Gamma; A) \times P(r) \longrightarrow C(\Gamma; A) \times P(r), \\
C(\Gamma; A) \times P(r) \longrightarrow C(\Gamma\big / A)\times P(A) \times P(r).
\end{array}
\right.
\]
Let $\ov \tau =\tau \times id_{P(r)}$, for each case.

\medskip
{\bf Theorem (2.3):}\ \ Suppose $\Gamma$ is a
knot graph and a subset $A$ of $V(\Gamma)$ satisfies the condition of
type 0, I, II, III, or IV in section 1.3.
$\ov \tau$ is the corresponding identification map above.
$\Psi$ is a canonical $P(k+r)$-structure for $C(\Gamma; A) \times P(r)$
and $\Psi'$ is a canonical $P(k+r)$-structure for
$\ov \tau(C(\Gamma; A) \times P(r))$.
Then, there is a $1-1$-correspondence
\[
\lambda:\Sigma_{k+r} \cdot \Psi \longrightarrow \Sigma_{k+r} \cdot \Psi'
\]
such that, for any $\Psi_1$ in $\Sigma_{k+r} \cdot \Psi$,
$\lambda(\Psi_1) \circ \ov \tau = \Psi_1$
on $C(\Gamma; A) \times P(r)$.
In fact, $\lambda$ satisfies $\lambda(\sigma \cdot \Psi_1)=
\sigma \cdot \lambda(\Psi_1)$,
for any $\sigma \in \Sigma_{k+r}$ and
$\Psi_1 \in \Sigma_{k+r} \cdot \Psi$, that is,
$\lambda$ is $\Sigma$-equivariant.

\medskip
{\bf (2.4):}\ \ Identification Rule
\begin{eqnarray*}
&&\wt \tau : C(\Gamma; A) \times P(r) \times \Sigma_{k+r} \cdot \Psi
\longrightarrow C(\Gamma; A) \times P(r) \times \Sigma_{k+r} \cdot \Psi' \\
&&\alpha \in C(\Gamma; A) \times P(r),
\Psi_1 \in \Sigma_{k+r} \cdot \Psi, \\
&&\wt \tau(\alpha, \Psi_1)=(\ov \tau(\alpha), \lambda(\Psi_1)).
\end{eqnarray*}
(Similar for the other case)

\medskip
{\bf Definition (2.5):}\ \ A knot graph $\Gamma$ is said to
be normal if
$\Gamma$ does not have any one of the following three
kinds of vertices:
(i) free vertice, (ii) univalent inner vertice, (iii)
bi-valent inner vertice.

\medskip
{\bf (2.6):}\ \ Suppose $n$ is a positive integer.
Let $\ov{\PP(n)}$ be the disjoint union of $P(\Gamma)$,
for all normal knot graphs $\Gamma$ with order $n$,
and $\PP(n)$ be the space $\ov{\PP(n)}$ quotiented by all possible identifications
of type 0, I, II, III and IV (given by the maps $\wt \tau$ on the boundary of
$P(\Gamma)$'s).

\vskip 6pt
Therefore, $\PP(n)$ is a $P(3n)$-space, and has a canonical
$P(3n)$-structure
$\wt {\Psi}:\PP(n) \longrightarrow P(3n)$, defined by:
for $\alpha \in C(\Gamma) \times P(3n-k)$ and
$\Psi_1$, a canonical $P(3n)$-structure for
$C(\Gamma) \times P(3n-k)$,
\[
\wt {\Psi} (\alpha, \Psi_1)=\Psi_1(\alpha).
\]
\newpage

\vskip 20pt
\Vjp
\noindent{\LARGE \bf \S 3 \  Infinitesimal knot graph}

\vskip 6pt
Let $l_0$ denote the line of $z$-axis,
$\{(0,0,t):t$ is a real
number, and consider the configuration space
$C(\Gamma, l_0)$ of knot graphs on $l_0$.
We shall define the translation and dilation on all
$C(\Gamma, l_0)$ and the extended translation and dilation relation on part
of $C(\Gamma, l_0)$, for which $\Gamma$ satisfies a splittable
condition.

\vskip 20pt
\Vjp
\noindent{\Large \bf 3.1 \ The translation and dilation relation}

\medskip
{\bf Definition (3.1):}\ \ Assume $(g, g(\Gamma))$
and $(g', g'(\Gamma))$ are two elements in $C(\Gamma, l_0)$.
$(g, g (\Gamma))$ and $(g', g'(\Gamma))$ are equivalent under
the translation and dilation relation, if there are
real number $\lambda$ and $t$, $\lambda >0$, such that
$g(v)=\lambda g' (v)+(0, 0, t)$,
for all $v \in V(\Gamma)$.
We may write the relation as $(g, g(\Gamma))$
$\buildrel {\rm t.d.} \over \sim$
$(g', g'(\Gamma))$.

\vskip 6pt
{\bf Notation (3.2):}\ \
$W_0(\Gamma)$ denotes the quotient space of
$C(\Gamma, l_0)$ by the translation and dilation relation.
Thus,
\[
\dim W_0(\Gamma)=\dim C(\Gamma, l_0) -2
\]

\vskip 6pt
{\bf Definition (3.3):}\ \ A knot graph $\Gamma$ is said to be
splittable, if there are subgraph $\Gamma_1$ and $\Gamma_2$ in $\Gamma$
such that the following two conditions hold:
\begin{namelist}{(ii)}
\item[{\rm (i)}]
$\Gamma=\Gamma_1 \cup \Gamma_2$, that is, $V(\Gamma)=V(\Gamma_1) \cup V(\Gamma_2)$
and $E(\Gamma)=E(\Gamma_1) \cup E(\Gamma_2)$.
\item[{\rm (ii)}]
$\Gamma_1 \cap \Gamma_2$ is either empty or a knot graph
consisting of exactly one base point (no edge and no inner point).

\end{namelist}

\noindent
And $\Gamma$ is said to be splited into $\Gamma_1$ and $\Gamma_2$

\vskip 6pt
Now, suppose $\Gamma$ is splited into $\Gamma_1$ and $\Gamma_2$.
Let $\rho : W_0(\Gamma) \longrightarrow W_0(\Gamma_1) \times W_0(\Gamma_2)$
be the map defined by $\rho(g, g(\Gamma))=((g |_{\Gamma_1}, g(\Gamma_1)),
(g|_{\Gamma_2}, g(\Gamma_2)))$.

\vskip 6pt
{\bf Definition (3.4):}\ \ Suppose $\Gamma$ is splited into
$\Gamma_1$ and $\Gamma_2$. Two elements $\alpha$ and $\beta$ in
$W_0(\Gamma)$ is said to be
equivalent under the extended translation and dilation relation,
(simply, ETD relation), if
$\rho(\alpha)=\rho(\beta)$ in
$W_0(\Gamma_1) \times W_0(\Gamma_2)$.

{\bf Note:}\ \ If $\Gamma_1 \cap \Gamma_2$ is empty,
$\dim (W_0(\Gamma_1) \times W_0(\Gamma_2))=\dim W_0(\Gamma)-2$.
If $\Gamma_1 \cap \Gamma_2$ consists of a base point,
$\dim (W_0(\Gamma_1) \times W_0(\Gamma_2))=\dim W_0(\Gamma)-1$.
Thus, $W_0(\Gamma) \big / ETD$ has dimension at most
$\dim C(\Gamma)-3$, and is also a $P(k)$-space.

\vskip 6pt
Now, we consider the boundary of $W_0(\Gamma)$.
For any $A$, a proper subset of $V(\Gamma)$, that is,
$A {\subset} V(\Gamma)$ and $A {\neq} V(\Gamma)$, let
$W_0(\Gamma; A)$ denote the quotient space of $C(\Gamma; A)$
(for the knot graphs on $l_0$) by the T.D. relation.
Because $l_0$ is a straight line,
$C(\Gamma; A)=C(\Gamma \big / A) \times D(A(\Gamma), x)$
and $D(A(\Gamma), x)$ is exactly $W_0(A(\Gamma))$.
Thus, $W_0(\Gamma; A) = W_0(\Gamma \big / A) \times W_0(A(\Gamma))$,
which has dimension equal to $\dim W_0(\Gamma)-1$.

\vskip 6pt
When $A(\Gamma)$ has no base point, the T.D. relation
for $C(A(\Gamma))$ is modified as follows:

\vskip 6pt
$(g, g(A(\Gamma)))$
$\buildrel {\rm t.d.} \over \sim$
$(g', g'(A(\Gamma)))$, if there exist real number $\lambda$ and
$y \in \setR^3$ such that $g(v)=\lambda g'(v) +y$,
for all $v \in A$.
\begin{eqnarray*}
W_0(A(\Gamma))
&&=C(A(\Gamma), l_0)  \big / \mbox{T.D. relation} \\
&&=D(A(\Gamma), x), \mbox{for any }x.
\end{eqnarray*}
Thus,
\[
\dim W_0(A(\Gamma))=\dim C(A(\Gamma))-4.
\]
\noindent
The extend T.D. relation could happen in $W_0(\Gamma \big / A)$
and in $W_0(A(\Gamma))$.
When $A(\Gamma)$ has no base point, $A(\Gamma)$ is splittable,
if and only if, $A(\Gamma)$ is a disjoint union of
$A_1(\Gamma)$ and $A_2(\Gamma)$, for some $A_1$ and $A_2$.
For example, $\Gamma$ is the following knot graph:
$V_0(\Gamma)$=\, $\{\,x_i,\,i=1,2,3,4,5,6,7\,\}$,$V_1(\Gamma)=
\{\,y_j,\,j=1,2,3,4,5\,\}$ and $E(\Gamma)=\{\,\{x_1,y_1\},\,\{x_2,y_1\},
\{x_3,y_2\},\,\{x_4,y_3\},\,\{x_5,y_4\},\\  \,\{x_6,y_5\},\,\{x_7,y_5\},\,
\{y_1,y_2\},\,\{y_2,y_3\},\,\{y_3,y_4\},\,\{y_4,y_5\}\,\}$.
$A=\{y_1, y_2, y_4, y_5\}$,
$A_1 = \{y_1, y_2\}$, $A_2=\{y_4, y_5\}$.
Then $A(\Gamma)$ is the splited into $A_1(\Gamma)$ and $A_2(\Gamma)$.

\medskip
{\bf (3.5):}\ \
If $A(\Gamma)$ is the disjoint union of
$A_1(\Gamma)$, $A_2(\Gamma), \cdots, A_r(\Gamma)$, then
the extended translation and dilation relation
associated with the splitting
$(A_1(\Gamma), A_2(\Gamma), \cdots, A_r(\Gamma)$
is defined as following: $\alpha=(g, g(A\Gamma)))$,
$\alpha'=(g', g'(A(\Gamma)))$ are two element in
$W_0(A(\Gamma))$,
$$\alpha \buildrel {\rm E.T.D.} \over \sim \alpha',
\mbox{ if } \, \,
g|_{A_i(\Gamma)} \buildrel {\rm T.D.} \over \sim
g'|_{A_i(\Gamma)},$$

\vskip 20pt
\Vjp
\noindent{\Large \bf  3.2  \ The combined space $\WW_0(n)$}

\vskip 6pt
Similar to the construction of $\PP(n)$, let
$\ov\PP(n, l_0)$ be the disjoint union of
$P(\Gamma)$ of $l_0$, for all normal knot graph
$\Gamma$ with order $n$, and $\PP(n, l_0)$
be the quotient space of $\ov \PP(n, l_0)$ by the identifications of type
0, I, II, III or IV, defined in section 1.3.

\vskip 6pt
And, $\WW_0(n)$ be the quotient space of $\PP(n, l_0)$
by the T.D. relations and the possible E.T.D. relations.

\vskip 6pt
We can understand $\WW_0(n)$ as the quotient space of
disjoint union of $W_0(\Gamma) \times P(3n-k)\times \Sigma_{3n} \cdot \Psi$,
with further quotients by the E.T.D. relations.

\vskip 6pt
{\bf Note:}\ \ The E.T.D. relations considered include all the E.T.D.
relation for configuration spaces and for the boundaries
$(W_0(\Gamma \big / A)$ or $W_0(A(\Gamma))$.

\vskip 20pt
\Vjp
\noindent{\Large \bf 3.3 \ The $S^1$-action on $\WW_0(n)$}

\vskip 6pt
Identify $\setR^3$ with $\setC \times \setR$,
$S^1$ may acts on $\setR^3$ by the multiplication on the
complex plane part, that is, the rotation about the $z$-axis.

\vskip 6pt
Thus, $S^1$ also acts on $C(\Gamma, l_0): \lambda \in S^1$,
$(g, g(\Gamma))\in C(\Gamma, l_0)$,
$\lambda \cdot (g, g(\Gamma))=(\lambda g, (\lambda g)(\Gamma))$.

\vskip 6pt
It is easy to see that all the identification
maps of any type in section 1.3 are $S^1$-equivariant maps
(for the case of $l_0$) and T.D. and E.T.D. relations are also
compatible with the $S^1$-action. Thus $S^1$ also acts on
$\WW_0(n)$.
Furthermore, the $S^1$-action on $\setR^3$ also induces an
$S^1$-action on $\setR\setP^2$ and $P(k)$, and hence,
the $P(k)$-structure map $\Psi$ of $\WW_0(n)$ is $S^1$-equivariant.

\medskip
{\bf Proposition (3.6):}\ \
\begin{namelist}{(ii)}
\item[{\rm (i)}]
$W_0(\Gamma)$, $\WW_0(n)$ are semifree $S^1$-space.
\item[{\rm (ii)}]
$P(k)$ has a unique fixed point of the $S^1$-action.
\end{namelist}

(Proof?)

\medskip
{\bf Notation (3.7):}\ \ $\HH(n)$ denotes the set of all fixed points of
$S^1$-action on $\WW_0(n)$, and $\DD(n)$ denotes the normal vector
bundle of $\HH(n)$ in $\WW_0(n)$.

\vskip 6pt
We find the vector bundle $\DD(n)$ over $\HH(n)$ has finite structure
group. Roughly speaking, the reason is the following:
Let $H(\Gamma)$ be the fixed point set of
$W_0(\Gamma)$. It is not hard to find that $\HH(n)$ is the
union of $H(\Gamma)$ and the restriction of the vector bundle
$\DD(n)$ to $H(\Gamma)$ is trivial.
And the identifications of type 0 to IV make the structure group
isomorphic the semi-direct product of $\bigoplus\limits_{r} \setZ_2$ and  $\Sigma_r$,
for some $r$.

\vskip 6pt
If we consider the $S^1$-equivariant map $\Psi:\WW_0(n) \longrightarrow
P(3n)$. $\xi_0$ denotes the unique fixed point in $P(3n)$.
Then $\HH(n)=\Psi^{-1} (\xi_0)$ and $\DD(n)$ is a complex vector
bundle of $2n$-dim, that is, fibre $=\setC^{2n}$.
$\dim \HH(n) =\dim H(\Gamma)=2n-2$, for non-splittable knot graph $\Gamma$.

\vskip 20pt
\Vjp
\noindent{\Large \bf 3.4  \ The space of infinitesimal knot graphs}

\vskip 6pt
For any $x \in S^2$, $l_x$ denotes the line
$\{ tx, t$ is a real number $\}$ with linear order $:t_1 x < t_2 x$,
for $t_1< t_2$.
Let $W_x(\Gamma)=C(\Gamma, l_x) \big / {\rm T.D. relation}$ and
$W(\Gamma)=\bigcup\limits_{x \in S^2} W_x(\Gamma)$.
Then $W(\Gamma)$ is a fibre
bundle over $S^2$ with fibre $W_0(\Gamma)$.

\vskip 6pt
Consider the following $SO(3)$-action on $W(\Gamma)$.
$\sigma \in SO(3)$, $(g, g(\Gamma))\in W_x(\Gamma)$,
\[
\sigma \cdot (g, g(\Gamma))=
(\sigma \cdot g, (\sigma \cdot g) (\Gamma)) \in W_{\sigma(x)} (\Gamma),
\]
where $(\sigma \cdot g)(v)=\sigma(g(v))$, for $v \in V(\Gamma)$.

\vskip 6pt
{\bf Notation (3.8):}\ \ For any $S^1$-space $X$, let
$SO(3)$
$\mathrel {\mathop\times\limits_{S^1}} X$
denote the quotient space
$SO(3) \times X \big / (\sigma, x) \sim (\sigma \lambda^{-1}, \lambda x)$,
for $\sigma \in SO(3)$, $x \in X$, and $\lambda \in S^1$,
where $S^1$ is identified with the
$SO(2)$-subgroup of $SO(3)$,
which rotates about the $z$-axis.

\vskip 6pt
Thus,
$SO(3)$
$\mathrel {\mathop\times\limits_{S^1}} W_0(\Gamma)$
is naturally diffeomorphic to $W(\Gamma)$ by the map sending
$(\sigma, \alpha) \longrightarrow \sigma \cdot \alpha$.

\vskip 6pt
$SO(3)$ also acts on $P(k)$. The map of $P(k)$-structure of
$W_0(\Gamma)$, $\Psi:W_0(\Gamma) \longrightarrow P(k)$, is a
$S^1$-equivariant map. $\Psi$ induces a $SO(3)$-equivariant map
$\wt\Psi:SO(3)$
$\mathrel {\mathop\times\limits_{S^1}} W_0(\Gamma)$
$\longrightarrow P(k)$,
$\wt \Psi (\sigma, \alpha)=\sigma \cdot \Psi(\alpha)$.

\vskip 6pt
$\wt \Psi$ is exactly the canonical $P(k)$-structure of
$W(\Gamma)$.
Thus, we conclude the following:

\medskip
{\bf (3.9):}\ \ $W_0(\Gamma)$ and its $S^1$-equivariant
$P(k)$-structure map $\Psi: W_0(\Gamma) \longrightarrow P(k)$
completely determines the $P(k)$-space $W(\Gamma)$.

\vskip 6pt
With some abuse of notation, we may use $\Sigma_k \cdot \Psi$ to represent
the set of all canonical $P(k)$-structure of $W(\Gamma)$ and let
$\WW(n)$ denote the quotient space of disjoint union of
$W(\Gamma) \times P(3n-k)\times \Sigma_{3n} \cdot \Psi$.
Thus, $\WW(n)$ is also equal to
$SO(3)$
$\mathrel {\mathop\times\limits_{S^1}} \WW_0(n)$.
And, the $P(3n)$-structure of $\WW(n)$ is also completely
determined by the $S^1$-equivariant $P(3n)$-structure map
$\Psi:\WW_0(n) \longrightarrow P(3n)$, as above.
This suggests that the study of
$\WW_0(n)$ and the associated $S^1$-equivariant $P(3n)$-structure
map will completely understand the ``degree'' of $\WW(n)$,
in the dimension $6n$.

\vskip 6pt
The following Theorem will be proved in the next two chapter.

\medskip
{\bf Theorem (3.10):}\ \ The degree homomorphism
\[
\wt \Psi_* : H_{6n} (\WW(n), \wt \Psi^{\#}(\OO_{3n}))
\longrightarrow H_{6n} (P(3n); \OO_{3n}) \approx \setZ
\]
``become'' a characteristic class of the complex vector
bundle $\DD(n)$ over $\HH(n)$.
(Thus, if $\DD(n)$ has finite structure group,
then the degree homomorphism is zero.)

\vskip 6pt
In the next chapter, we shall develop a theory
for the $P(k)$-space with $S^1$-equivariant $P(k)$-structure.

\medskip
{\bf Remark:}\ \ In $P(3n)$, let us consider the open set
$\{ (y_1, y_2, \cdots, y_{3n}) \in P(3n):y_i$
is not perpendicular to the $z$-axis, for all
$i=1, 2, \cdots,  3n\}$, it is the canonical
$\setR^{6n}$ neighborhood of the fixed point
$\xi_0$.
Then the normal bundle $\DD(n)$ can be identified with
$\Psi^{-1}(\setR^{6n})$ canonically, and
$\Psi:\DD(n) \longrightarrow \setR^{6n}$ is an
injective linear map on each fibre of $\DD(n)$.
\newpage

\noindent{\Large \bf \S 4  \ $S^1$-equivariant $P(k)$-structure}

\vskip 6pt
Suppose $(X, f)$ is a $P(k)$-space, $X$ has an
$S^1$-action and $f:X \longrightarrow P(k)$ is an
$S^1$-equivariant map.
Let $SO(X)=SO(3)$
$\mathrel {\mathop\times\limits_{S^1}} X$,
then $SO(X)$ is also a $P(k)$-space.
In this part, we study the relation between the degree
homomorphism of $P(k)$-space $(SO(X), \wt f)$ and the
homomorphism $\ov f_* : H_* (X \big / S^1)\longrightarrow H_*(P(k) \big / S^1)$,
and furthermore.

\vskip 20pt
\Vjp
\noindent{\Large \bf  4.1 \ The $S^1$-action on $P(k)$}

\vskip 6pt
$P(k)=\prod\limits_{k} \setR\setP^2$.
For any $\lambda = e^{i \theta}$, $0 \le \theta \le 2 \pi$,
and $(x, y, z)\in \setR^3$,
$\lambda \cdot (x,y,z)=(x\cos\theta-y\sin \theta, x\sin\theta+y\cos\theta,z)$.
This gives $S^1$-action on $S^2$ and $\setR\setP^2$.
For $\prod\limits_{k} S^2$ and $\prod\limits_{k}\setR\setP^2$,
$S^1$ acts on them diagonally, that is,
$\lambda \cdot (y_1, y_2, \cdots, y_k)=
(\lambda \cdot y_1, \lambda \cdot y_2, \cdots, \lambda \cdot y_k)$.
Suppose $X$ is a $S^1$-space.
$x_0$ in $X$ is said to the fixed point of the $S^1$-action,
if $\lambda \cdot x_0 = x_0$, for all $\lambda \in S^1$.
The $S^1$-action on $X$ is said to be free, if
$\lambda \cdot x \ne x$, for all $\lambda \ne 1$, and for all $x \in X$.
The $S^1$-action on $X$ is said to be semifree,
if $X$-(the fixed point set of the $S^1$-action) is free.
The following are the main $S^1$-spaces we need:
\begin{namelist}{(ii)}
\item[{\rm(i)}]
The $S^1$-action on $\prod\limits_{k}S^2$ is semifree, and its fixed points
set is $\prod\limits_{k}$ (fixed point set of $S^2$), which has $2^k$ elements.
\item[{\rm(ii)}]
The $S^1$-action on
$P(k)=\prod\limits_{k} \setR\setP^2$ has only one fixed point.
\end{namelist}

\medskip
{\bf Defintion (4.1):}\ \ For any $S^1$-space $X$, let
$SO(X)$ denote the quotient space $SO(3)$
$\mathrel {\mathop\times\limits_{S^1}}X$, that is,
$SO(3) \times X$ quotiented by the equivariant
relation: for $\sigma \in SO(3)$,
$x \in X$, $\lambda \in S^1$,
$(\sigma, \lambda x) \sim (\sigma \lambda^{-1} \cdot x)$.

\vskip 6pt
($S^1$ is the identified as the subgroup of rotations about the $z$-axis
in $SO(3)$.)

\vskip 6pt
Let Fix $(X)$ denote the set of fixed points in the
$S^1$-space $X$. For any $i > \dim \Fix (X)$,
there is a homomorphism $\beta_1$ from
$H_i(X \big / S^1)$ to $H_{i+1}(X)$.
Similarly, $SO(3) \times X$ is a $S^1$-space $:$ $\lambda \cdot
(\sigma , x)=(\sigma \lambda^{-1}, \lambda \cdot x)$,
there is a homomorphism
$\beta_2: H_{i+3} (SO(X)) \longrightarrow H_{i+4}(SO(3) \times X)$
and $\ov \beta : H_{i+1} (X) \longrightarrow H_{i+4} (SO(3) \times X)$.

Consider the following diagram
\[
 \begin{array}{ccc}
  H_{i+1}(X)           & \mapright{\ov \beta} & H_{i+4}(SO(3)\times X) \\  \\
  \mapupr{\beta_1} &                      & \mapupr{\beta_2} \\ \\
  H_i(X\big / S^1) & \mapright{\beta}     & H_{i+3}(SO(X)) \\
 \end{array}
\]
It is obvious that $Im(\ov \beta \circ \beta_1) \subset Im(\beta_2)$.

\medskip
{\bf Question:}\ \ Is there a homomorphism
$\beta:H_i(X \big / S^1) \longrightarrow H_{i+3} (SO(X))$
such that the above diagram is commutative?
Actually, $\beta$ is completely an analogue of $\beta_1$.
$SO(X)$ is an $SO(3)$-space,
$\sigma [\sigma_1, x)]=[(\sigma \sigma_1, x)]$,
and orbit space $SO(X) \big / SO(3)$ is exactly
$X \big / S^1$.
If we assume that the $S^1$-action on $X$ is semifree, then
the fibre map $SO(X) \longrightarrow X \big / S^1$ has fibres homeomorphism
to $SO(3)$ except over the fixed points.
Therefore, the homomorphism $\beta:H_i(X \big / S^1) \longrightarrow H_{i+3}
(SO(X))$ exists and satisfies the functorial property:
for any $S^1$-equivariant map
$f:X \longrightarrow Y$ of semifree $S^1$-spaces $X$ and $Y$,
$i > \dim(\Fix(X))$ (and $\dim (\Fix(Y)$), the
following diagram commutes:

\[
 \begin{array}{ccc}
  H_i(X\big / S^1) & \maprightu{} & H_{i+3}(SO(X)) \\  \\
  \mapdownr{f_*} & & \mapdownr{\ov f_*} \\ \\
  H_i(Y\big / S^1) & \maprightu{} & H_{i+3}(SO(Y)) \\ \\
 \end{array}
\]
Why we need the condition $i >\dim \Fix(X)$?

\vskip 6pt
When $i > \dim \Fix(X)$, we may assume there is a cell
decomposition of $X \big / S^1$, such that
$\Fix(X) \subset$ the $(i-1)$-skeleton of
$X \big / S^1$. In a corresponding cell decomposition in $X$,
the behavior of $\Fix (X)$ will not affect the boundary
operator of $(i+1)$-cells.

\medskip
{\bf Remark:} If dim$X=i+1$, then $\beta_1$,$\beta_2$ and $\ov \beta$
are all isomorphism and the existence of $\beta$ is obvious.
When applying the theory in Chapter 5, it is the situation.

\vskip 20pt
\Vjp
\noindent{\Large \bf  4.2 \ $P(k)$-space with a semifree $S^1$-action}

\vskip 6pt
Suppose $X$ is a semifree $S^1$-space and $f: X \longrightarrow P(k)$
is an $S^1$-equivariant map. $\xi_0$ denotes the unique fixed point
of $S^1$-action on $P(k)$ and $Y$ denotes the inverse image
$f^{-1} (\xi_0)$.
Furthermore, let $U_0$ be an $S^1$-invariant contractible open
neighborhood of $\xi_0$ in $P(k)$ and
$U_1=f^{-1} (U_0)$, it is an $S^1$-invariant open neighborhood of $Y$
in $X$.
Then, $U_0 - \xi_0$ and $U_1 -Y$ are free $S^1$-spaces.

\vskip 6pt
Let $\ov f : X \big / S^1 \longrightarrow P(k) \big / S^1$
be the map induced by $f$,
and $\wt f : SO(X) \longrightarrow P(k)$ be the map
defined by:
\[
\wt f(\sigma , x)=\sigma \cdot f(x).
\]

\vskip 6pt
Let $\OO$ be the orientation sheaf of $P(k)$,
$\ov \OO$ be the orientation sheaf of $P(k) \big / S^1$  such that
$\pi^{\#} (\ov \OO) =\OO$,
where $\pi:P(k) \longrightarrow P(k) \big / S^1$ is the quotient map.

\vskip 6pt
The homomorphism $\beta$ described in the last section 4.1,
together with the local coefficient $\OO$ give the following
diagram:

$$\matrix{
  H_{2k-3}(X\big / S^1, \ov f^{\#} (\ov \OO)) & \maprightu{\beta}
  & H_{2k}(SO(X), \wt f^{\#}(\OO))  &\cr\cr
  \mapdownr{\ov f_*} & & &\mapsear{SO( f)_*}\quad\quad\quad\quad \cr\cr
  H_{2k-3}(P(k)\big / S^1, \ov \OO) & \maprightu{\beta}&\mapvert{\wt f_*}  \quad \quad \maprightu{}
   &H_{2k}(SO(P(k)),\OO')  \cr \cr
  &\mapsear{\ov \beta} &\mapdownr{} \quad\quad \quad \quad \quad \, \, \, \,
  &\mapswr{\gamma} \quad\quad\quad\quad  \cr\cr
  &&H_{2k}(P(k), \OO) \quad \quad& \cr\cr}
$$

\noindent
When $\ov f$ restricts to $U_1 - Y \big / S^1$, we have
$\buildrel =\over f: (U_1-Y \big / S^1) \longrightarrow (U_0-\xi_0\big / S^1)$.

Consider
\begin{eqnarray*}
\buildrel = \over f_*: H_{2k-4}(U_1 - Y \big / S^1)
\longrightarrow && H_{2k-4} (U_0-\xi_0 \big / S^1)\\
\approx &&H_{2k-4} (\setC(\setP^{k-1}) \\
\approx &&\setZ.
\end{eqnarray*}

\medskip
{\bf Theorem (4.2):}\ \ If
$\buildrel = \over f_* :H_{2k-4} (U_1 - Y \big / S^1)
\longrightarrow H_{2k-4} (U_0 - \xi_0 \big / S^1)$
is a zero-homomorphism, then the degree homomorphism
$\wt f_* : H_{2k} (SO(X), \wt f^{\#} (\OO))$  $ \longrightarrow H_{2k} (P(3n),
\OO)$ sends $\beta(x)$ to $0$, for all $x \in H_{2k-3} (X \big / S^1,
\ov f^{\#}(\ov \OO))$.
(That is, $\buildrel = \over f_* =0 \Longrightarrow \wt f \circ \beta =0$.)

\vskip 6pt
{\bf Note:}\ \ In the statement, we do not need the commutative property:
$SO(f)_* \circ \beta =\ov f_* \circ \beta_*$,
probably, it is true, but I never prove it.
(because $P(k)$ is not semifree)

\vskip 20pt
\Vjp
\noindent{\Large \bf  4.3 \ Preliminary to the proof of Theorem (4.2)}

\vskip 6pt
$\prod\limits_{k} S^2$ is the universal covering space of $P(k)$,
let $S(k) = \prod\limits_{k} S^2$ and $q: S(k) \longrightarrow P(k)$
be covering projection.

\vskip 6pt
Let $M$ the pullback of $S(k)$ by
$f: X \longrightarrow P(k)$,
that is, $M=\{(x,y) \in X \times S(k) : f(x)=q(y)\}$.
Thus, we have the following commutative diagram.

\medskip
{\bf (4.3)}\ \ \hskip 2cm
$
\matrix{
  M & \maprightu{g} & S(k) \cr\cr
  \mapdownr{q'} & & \mapdownr{q} \cr\cr
  X & \maprightu{f} & P(k) \cr\cr}
$

\noindent
$S(k)$ is also a $S^1$-space and $q$ is a $S^1$-equivariant map.
$S^1$ acts on $X \times S(k)$ by:
$\lambda \in S^1$, $\lambda \cdot (x,y)=(\lambda x, \lambda y)$.

\vskip 6pt
If $(x, y) \in M$, $\lambda \cdot (x, y)$ is also in $M$.
This defines an $S^1$-action on $M$ such that both
$g:M \longrightarrow S(k)$ and the covering projection
$q':M \longrightarrow X$ are $S^1$-equivariant.
This means that the diagram (4.3) is a commutative diagram of
$S^1$-equivariant maps.
This induces the following commutative diagram of
$SO(3)$-equivariant diagram.

\medskip
{\bf (4.4)}\ \ \hskip 2cm
$\matrix{
  SO(M) & \maprightu{\wt g} & S(k) \cr\cr
  \mapdownr{SO(q')} & & \mapdownr{q} \cr \cr
  SO(X) & \maprightu{\wt f} & P(k)\cr\cr}
$

\vskip 6pt
(Note:\ \ $q$ is originally an $SO(3)$-equivariant map).

\vskip 6pt
Because $q' : M \longrightarrow X$ is an $S^1$-equivariant covering
projection and $X$ is assumed to be a semifree $S^1$-space,
$M$ is also a semifree $S^1$-space.
Thus, both $X$ and $M$ provide $\beta$-homomorphism from
$H_{2k-3}((\cdot) \big / S^1)$ to
$H_{2k} (SO(\cdot))$.
They form the following commutative diagram

\medskip
{\bf (4.5)}\ \
$\matrix{
  H_{2k-3}(M \big / S^1)  & \mapright{\beta_M} & H_{2k}(SO(M))
  &\mapright{\wt g_*} &H_{2k}(S(k)) \cr \cr
  \mapdownr{\ov {q'}_*} & &\mapdownr{SO(q')_*} &&\mapdownr{q_*}\cr \cr
  H_{2k-3}(X \big / S^1, \ov f^{\#}(\ov \OO))  & \mapright{\beta_X} & H_{2k}(SO(X), \wt f^{\#}(\OO))
  &\mapright{\wt f_*} &H_{2k}(P(k), \OO)\cr \cr}
$
(The local coefficients for $M \big / S^1$,
$SO(M)$, and $S(k)$ are trivial.)

\medskip
{\bf Proposition (4.6):}\ \
$\wt f_* \circ \beta_X=0$, if and only if,
$\wt g_* \circ \beta_M = 0$.

\medskip
{\bf Proof:}\ \ Because $H_{2k}(S(k))\approx H_{2k} (P(k), \OO) \approx \setZ$,
it is enough to show this proposition for the coefficient over rationals
$Q$.
But, when over rational, both
$q_* : H_{2k}(S(k)) \longrightarrow H_{2k}(P(k), \OO)$, and
$\ov {q'}_* : H_{2k-3}(M \big / S^1) \longrightarrow H_{2k-3}(X \big / S^1)$
are isomorphisms.
By the commutative property of (4.5), we prove this proposition.
(Note: $M \big / S^1 \longrightarrow X \big / S^1$ is also
a covering space)

\vskip 6pt
In this next section, the problem will transform to the
$S^1$-equivariant $S(k)$-map,
$g:M\longrightarrow S(k)$, and the associated maps.

\vskip 20pt
\Vjp
\noindent{\Large \bf  4.4 \ Proof of Theorem (4.2)}

\vskip 6pt
Let $F$ denote the fixed point set of the $S^1$-action
on $S(k)$. $F$ has $2^k$ element,
$V=q^{-1}(U_0)$ has $2^k$ components $V_y$,
for each $y \in F$.
Each $V_y$ is diffeomorphic to $U_0$ by $q$.
$U$ denotes the inverse image
$g^{-1} (V)$,
$N$ denotes $g^{-1}(F)$. Similarly, $U$ is the
disjoint union of $U_y$, for $y \in F$, and
$N =\bigcup\limits_{y \in F} N_y$, where
$N_y=g^{-1}(y)$,
$U_y=g^{-1} (V_y)$, for $y \in F$.
Note: $V_y$ is a contractible open neighborhood of $y$,
$(y \in F)$, $H_i(V_y)=0$, $H_i(V)=0$, for $i >0$.

\vskip 6pt
Consider the commutative diagrams

\medskip
{\bf (4.7)}\ \
$\matrix{
H_{2k-3} (M \big / S^1) &\maprightu{\beta_M} &H_{2k}(SO(M))
&\maprightu{\wt g_*}    & H_{2k} (S(k)) \cr \cr
\mapdownr{\ov g_*}    & &\mapdownr{SO(g)_*} &\mapnear{\alpha_*} & \cr \cr
H_{2k-3}(S(k)\big / S^1) &\maprightu{\beta_S} & H_{2k} (SO(S(k)))
&&\cr\cr}
$
where $\alpha : SO(S(k)) \longrightarrow S(k)$
is given by $\alpha(\sigma, y)=\sigma \cdot y$
(thus, $\wt g_* \circ \beta_M=\alpha_* \circ \beta_S \circ \ov g_*$),
and

\medskip
{\bf (4.8)} \ \
$\matrix{
H_{2k-3}(M \big / S^1) &\maprightu{\ov g_*} &H_{2k-3}(S(k) \big / S^1) \cr\cr
\mapdownr{} &&\mapdownr{(I):\mbox{injective}} \cr\cr
H_{2k-3}(M\big / S^1, M-N \big / S^1) &\maprightu{}
&H_{2k-3}(S(k)\big / S^1, S(k)-F \big / S^1) \cr\cr
\mapupr{\wr (II')} &&\mapupr{\wr (II)} \cr\cr
H_{2k-3} (U \big / S^1, U-N \big / S^1)
&\maprightu{} &H_{2k-3} (V \big / S^1, V-F \big / S^1) \cr\cr
\mapdownr{} &&\mapdownr{\wr(III)} \cr\cr
H_{2k-4} (U - N \big / S^1) &\maprightu{\buildrel = \over {g_*}}
&H_{2k-4}(V-F \big / S^1) \cr\cr
\bigg | \bigg | &&\bigg | \bigg | \cr\cr
\bigoplus\limits_{y \in F} H_{2k-4} (U_y - N_y \big / S^1)
& &\bigoplus\limits_{y \in F} H_{2k-4} (V_y - y \big / S^1) \cr\cr}
$

\medskip
{\bf Lemma (4.9):}\ \
\begin{namelist}{(iii)}
\item[{\rm (i)}]
The map (I) is an injective homomorphism.
\item[{\rm (ii)}]
The map (II) and (II') are excision isomorphism.
\item[{\rm (iii)}]
The map (III) is an isomorphism.
\end{namelist}

\vskip 6pt
We shall prove Lemma (4.9) later and use the results in (4.9)
to finish the proof of Theorem (4.2).

\vskip 6pt
By the commutativity of (4.7), if
$\ov g_* =0$ in $\dim(2k-3)$, then $\wt g_* \circ \beta_M=0$.

\vskip 12pt
By the commutativity of (4.8), if
$\buildrel = \over g_* =0$ in $\dim (2k-4)$, then
$\ov g_*=0$ in $\dim(2k-3)$.
But, $\buildrel = \over g_*$ is the direct sum of
$(\buildrel = \over g_y)_*$ $:$
$H_{2k-4}(U_y - N_y \big / S^1) \longrightarrow H_{2k-4}(V_y -y \big / S^1)$,
for $y \in F$. And, for each $y \in F$,
$\buildrel = \over g_y$ is just a copy of
$\buildrel = \over f : U_1 -Y \big / S^1 \longrightarrow U_0 -\xi_0 \big /
S^1$. Thus, if $\buildrel = \over f_* =0$ in $\dim (2k-4)$
then $\buildrel = \over g_*=0$ in $\dim (2k-4)$, and hence,
$\wt g_* \circ \beta_M=0$ which is equivalent to
$\wt f_* \circ \beta_X=0$.

\vskip 6pt
Now return to Lemma (4.9), (ii) and (iii) are obvious. To prove (i),
it is enough to show that
$H_{2k-3}(S(k)-F \big / S^1)$ is trivial.
In next section, we study
$S(k)-F$ and prove the result (4.13) we need.

\vskip 20pt
\Vjp
\noindent{\Large \bf  4.5 \ The space $S(k)-F$}

\vskip 6pt
Let $F_0 = \{(0,0,1), (0, 0, -1)\}$,
the fixed point set of $S^2$. Let $T=S^2 - F_0$.
Then
$$
S(k)-F = (T \times \prod\limits_{k-1} S^2)
\cup (S^2 \times T \times \prod\limits_{k-2} S^2) \cup
\cdots \cup (T \times \prod\limits_{k-1} S^2).
$$

\medskip
{\bf Lemma (4.10):}\ \ If $Y$ is a $S^1$-space, then
$Y \times S^1 \big / S^1$ is
homeomorphic to $Y$.

\vskip 6pt
{\bf Proof:}\ \
$f:Y \longrightarrow Y \times S^1 \big / S^1$,
$f(y)=(y,1)$, and $g:Y \times S^1 \big / S^1 \longrightarrow Y$
by $g(y, \lambda)=\lambda^{-1} \cdot y$.
Then both $f \circ g$ and $g \circ f$ are identity map.

\vskip 6pt
Because $T=(0,1) \times S^1$,
$T \big / S^1 \approx (0,1)$ and
$T \times S(k-1)$ is homeomorphic to $(0,1) \times S(k-1)$.

\vskip 6pt
Thus $H_i(T \times S(k-1))=0$,
$i >2k-2$, $H_{2k-2}(T \times S(k-1)) \approx \setZ$ and
$H_{2k-3}(T \times S(k-1))=0$.

\vskip 6pt
Let $X_1=T\times S(k-1)$, $X_2=S^2 \times T \times S(k-2)$, $\cdots$,
$X_i=S(i-1) \times T \times S(k-i)$, $\cdots$, $X_k$ $=$
$S(k-1) \times T$,
$T=S^2 - \{(0,0,1), (0,0, -1)\}$.
$X_i=S(i-1) \times T \times S(k-i)$,
it is an open subset of $S(k)$.
$X_{i_1, i_2, \cdots, i_r}$, $i_1 < i_2 < \cdots < i_r$,
is the intersection of $X_{i_j, j=1, 2, \cdots, r}$,
that is,
$X_{i_1} \cap X_{i_2} \cap \cdots \cap X_{i_r}$.
$\{X_1, X_2, \cdots, X_k\}$ is an open covering of $S(k)-F$.
Let $Y_i = X_i \big / S^1$,
$Y=S(k)-F \big / S^1$.
Then $\{Y_1, Y_2, \cdots, Y_k\}$ is an open covering of $Y$.
Similarly, let $Y_{i_1, i_2, \cdots, i_r}=Y_{i_1} \cap Y_{i_2} \cap
\cdots \cap Y_{i_r}$, it is equal to $X_{i_1, i_2, \cdots, i_r} \big / S^1$.

\Vjp
\noindent
{\bf
Homology of $Y_{i_1, i_2, \cdots, i_r}$ at dimensions
$2k-r-2$ and $2k-r-3$}

\vskip 6pt
Without loss of generality, it is enough to understand
$H_{2k-r-3}(Y_{1, 2, \cdots, r})$ and
$H_{2k-r-2}(Y_{1, 2, \cdots, r})$.

\begin{enumerate}
\item
$r=1$. $H_{2k-3}(Y_1)$:

$X_1=T \times S(k-1)$ which has the equivariant homotopy type as
$S^1 \times S(k-1)$, the map
\begin{eqnarray*}
S^1 \times S(k-1) & \hookrightarrow  & X_1 \\
(\lambda, y) &\longrightarrow &(\lambda, 0, y) ,
\end{eqnarray*}
sending $S^1$ into the equator of $T \subset S^2$, is an explicit
equivariant homotopy equivalence.
Thus $Y_1 = X_1 \big / S^1 \simeq S^1 \times S(k-1) \big / S^1 \approx S(k-1)$.
Thus, $H_{2k-3}(Y_i)=0$, $\forall$ $i=1, 2, \cdots, k$.

\vskip 10pt
$H_{2k-4}(Y_1)=$? :

Consider the embedding
\begin{eqnarray*}
& &\psi^2 : S^1 \times S(k-2) \longrightarrow  X_1, \\
& &\psi^2 : (\lambda, y)=(\lambda a, \lambda b, y),
\end{eqnarray*}
where $a \in T$, $b \in S^2$, $(\lambda a, \lambda b, y)\in T \times S^2
\times S(k-2)$.
$Im \psi^2$ is a $(2k-3)$-dim.
$S^1$-invariant set in $X_1$.
Similarly, we may define
\begin{eqnarray*}
& & \psi^i
: S^1 \times S(k-2) \longrightarrow X_1, 2 \le i \le k,  \\
& &\psi^i(\lambda, (y_1, y_2))=(\lambda a, y_1, \lambda b, y_2),
\end{eqnarray*}
where $y_1 \in S(i-2)$, $y_2 \in S(k-i)$,
$a \in T$, $b \in S^2$.

$\psi^2, \psi^3, \cdots, \psi^k$ represent $k-1$ independent homology
classes in $H_{2k-4} (Y_1)$.
To specify the classes in $H_{2k-4}(Y_1)$,
we denote by
$\psi_1^i$, $i=2, 3, \cdots, k$ the embeddings in $Y_1$.

(Note: $a \in T$ and $b \in S^2$ are two chosen points and
$y, y_1, y_2$ denote the variable in
$S(k-1)$, $S(i-2)$, $S(k-i)$, respectively)

Similarly, we have
$\psi_j^i:S^1 \times S(k-2) \longrightarrow X_j$,
$i=1, 2, \cdots, j-1, j+1, j+2, \cdots, k$,
embeddings in $X_j$, which represents $(2k-4)$-dimensional
homology classes in $Y_j$.

\medskip
\item
$H_{2k-4}(Y_{1, 2})$:

Consider the embedding
\begin{eqnarray*}
& &\varphi_{12}: S^1 \times S(k-2) \longrightarrow X_{1,2}=
T \times T \times S(k-2), \\
& &\varphi_{12}(\lambda, y)=(\lambda a, \lambda b, y),
\end{eqnarray*}
where $a, b$ are constant points in $T$ and $y$ is a
variable in $S(k-2)$.
$\varphi_{12}$ represents the only independent class in
$H_{2k-4}(Y_{1,2})$.
Similarly, embedding
\[
\varphi_{i_1, i_2}(i_1 < i_2):S^1 \times S(k-2)
\longrightarrow X_{i_1, i_2}
\]
represents the only one independent homology
class in $Y_{i_1, i_2}$.

\vskip 10pt
$H_{2k-5}(Y_{1,2})$:

Consider the embedding
\begin{eqnarray*}
& & \psi_{1,2}^j: S^1 \times S^1 \times S(k-3) \longrightarrow
X_{1,2}, j\geq 3 \\
& & \psi_{12}^j (\lambda_1, \lambda_2, y_1, y_2)=
(\lambda_1 a_1, \lambda_2 a_2, y_1, \lambda_1 b, y_2),
\end{eqnarray*}
where $\lambda_1, \lambda_2$ are variables in $S^1$,
$a_1, a_2$ are constant points in $T$,
$b$ is a constant point in $S^2$,
$y_1$ is a variable in $S(j-3)$ and $y_2$ is a variable in
$S(k-j)$.

Then the image of $\psi_{1,2}^j$ is an
$S^1$-invariant $(2k-4)$-dimensional subset of $X_{1,2}$,
($j$ could be an integer from $3$ to $k$).
These embeddings represent $(k-2)$ independent $(2k-5)$-dimension
homology classes in $Y_{1,2}$.

\item
$H_{2k-5}(Y_{1,2,3})$:

Consider the embedding
\begin{eqnarray*}
& & \varphi_{1,2,3}^1 : S^1 \times S^1 \times S(k-3) \longrightarrow
X_{1,2,3}, \\
& & \varphi_{1,2,3}^2 (\lambda_1, \lambda_2,  y)
=(\lambda_1 a_1, \lambda_1 a_2, \lambda_2 a_3, y).
\end{eqnarray*}
Similarly, $\varphi_{1,2,3}^3(\lambda_1, \lambda_2, y)
=(\lambda_1 a_1, \lambda_2 a_2, \lambda_1 a_3, y)$.

Then,
$\varphi_{1,2,3}^2$ and $\varphi_{1,2,3}^3$
represent $2$ independent homology classes in $Y_{1,2,3}$.

\vskip 10pt
$H_{2k-6}(Y_{1,2,3})$:

Consider the embedding
\begin{eqnarray*}
& & \psi_{1,2,3}^i : S^1 \times S^1 \times S^1 \times S(k-4) \longrightarrow
X_{1,2,3}, i \ge 4, \\
& & \psi_{1,2,3}^i (\lambda_1, \lambda_2, \lambda_3, y_1, y_2)
=(\lambda_1 a_1, \lambda_2 a_2, \lambda_3 a_3, y_1, \lambda_1 b, y_2),
\end{eqnarray*}
where $y_1$ is a variable in $S(i-4)$, $y_2$ is a variable
in $S(k-i)$, $a_1, a_2, a_3$ are constant points in $T$, $b$
is a constant point in $S^2$,
$\lambda_j \in S^1$, $j=1, 2, 3$.

Then $\psi_{1,2,3}^i$, $i=4, 5, \cdots, k$. represent $(k-3)$ independent
homology classes in $Y_{1,2, 3}$.

\item
$(H_{2k-r-2}(Y_{i_1, i_2, \cdots, i_r}))$

Simply: $H_{2k-r-2}(Y_{1,2,\cdots, r})$,
\begin{eqnarray*}
& & \varphi_{1,2, \cdots, r}^2 :(\prod\limits_{r=1} S^1) \times S(k-r)
\longrightarrow X_{1, 2, \cdots, r}, \\
& & \varphi_{1, 2, \cdots, r}^2
(\lambda_1, \lambda_2, \cdots, \lambda_{r-1}, y)=
(\lambda_1 a_1, \lambda_1 a_2,
\lambda_2 a_3, \lambda_3 a_4, \cdots,
\lambda_{r-1} a_r, y),
\end{eqnarray*}
$\lambda_i$, $i=1, 2, \cdots, r-1$
are variables in $S^1$, $a_1, a_2, \cdots, a_r$ are constant
points in $T$,
$y$ is a variable in $S(k-r)$.
\begin{eqnarray*}
& & \varphi_{1, 2, \cdots, r}^3
(\lambda_1, \lambda_2, \cdots, \lambda_{r-1}, y) \\
& & = (\lambda_1 a_1, \lambda_2 a_2,
\lambda_1 a_3, \lambda_3 a_4,
\lambda_4 a_5, \cdots,
\lambda_{r-1} a_r, y) \\
& & \vdots \\
& & \varphi_{1, 2, \cdots, r}^r
(\lambda_1, \lambda_2, \cdots, \lambda_{r-1}, y) \\
& & = (\lambda_1 a_1, \lambda_2 a_2,
\lambda_3 a_3, \cdots
\lambda_{r-1} a_{r-1},
\lambda_{1} a_r, y)
\end{eqnarray*}
(Remark: $\varphi_{1, 2, \cdots, r}^i
(\lambda_1, \lambda_2, \cdots, \lambda_{r-1}, y)$
$=$ $(\lambda_1 a_1$, $\cdots$,  $\lambda_1 a_i$, $\cdots$ $y)$)

$\varphi_{1, 2, \cdots, r}^2$,
$\varphi_{1, 2, \cdots, r}^3$, $\cdots$, $\varphi_{1, 2, \cdots, r}^r$
represent $(r-1)$ independent homology classes in $Y_{1, 2, \cdots, r}$.

(Note: $Y_{1, 2, \cdots, r}$ is homotopy equivalence to
$(\prod\limits_{r-1} S^1) \times (\prod\limits_{k-r} S^2)$,
rank $H_{2k-r-2} (Y_{1, 2, \cdot, r})=r-1$.)

\vskip 10pt
$H_{2k-r-3}(Y_{i_1, i_2, \cdots, i_r})$,
$1 \le i_1 < i_2 < \cdots < i_r \le k$:

For convenience, consider
$i_1=1, i_2=2, \cdots, i_r=r$.

Let
\begin{eqnarray*}
& & \psi_{1, 2, \cdots, r}^i : \prod\limits_r S^1 \times S(k-r-1)
\longrightarrow X_{1, 2, \cdots, r} \mbox{ be defined by}:  \\
& &\psi_{1, 2, \cdots, r}^i (\lambda_1, \lambda_2, \cdots, \lambda_r, y_1, y_2)
=(\lambda_1a_1, \lambda_2, a_2, \cdots, \lambda_r a_r, y_1, \lambda_1 b, y_2)
\end{eqnarray*}
where $k \ge i \ge r+1$,
$y_1 \in S(i-r-1)$,
$y_2 \in S(k-i)$,
$a_1, a_2, \cdots, a_r$ are constant points in $T$,
and $b$ is a constant point in $S^2$.

$\psi_{1, 2, \cdots, r}^{r+1}$,
$\psi_{1, 2, \cdots, r}^{r+2}$,
$\cdots$,
$\psi_{1, 2, \cdots, r}^{k}$
represent $(k-r)$ independent homology classes in $Y_{1, 2, \cdots, r}$.

Note: $Y_{i_1, i_2, \cdots, i_r}$ is homology equivalence to
$\prod\limits_{r-1} S^1 \times \prod\limits_{k-r} S^2$,
thus, $H_{2k-r-3}(Y_{i_1, i_2, \cdots, i_r})$ has
rank $(k-r)$.

\end{enumerate}

\vskip 6pt
Let $\ov G_r$ be the direct sum of $H_{2k-r-3}(Y_{i_1, i_2, \cdots, i_r})$,
for all $(i_1, i_2, \cdots, i_r)$ satisfying
$1 \le i_1 < i_2 < \cdots < i_r \le k$.
And $G_{i_1, i_2, \cdots, i_{r+1}}$ be the subgroup of $\ov G_r$,
generated by the $(r+1)$ elements
$\psi_{i_2, i_3, \cdots, i_{r+1}}^{i_1}$,
$\psi_{i_1, i_3, \cdots, i_{r+1}}^{i_2}$,
$\cdots$,
$\psi_{i_1, i_2, \cdots, i_{r}}^{i_{r+1}}$,
for any $(i_1, i_2, \cdots, i_{r+1})$ satisfying
$1 \le i_1 < i_2 < \cdots < i_{r+1} \le k$.

\vskip 6pt
Thus,
\[
\ov G_r=\bigoplus\limits_{1 \le i_1 < i_2 < \cdots < i_{r+1} \le k}
G_{i_1, i_2, \cdots, i_{r+1}}.
\]

\vskip 6pt
Let
\[
\rho_{i_1, i_2, \cdots, i_r}^j : Y_{i_1, i_2, \cdots, i_r}
\longrightarrow Y_{i_1, \cdots, \wh i_j, \cdots, i_r}
\]
be the inclusion, $1 \le j \le r$, and
\[
 \eta_{i_1, i_2, \cdots, i_r} :
H_{2k-r-2}(Y_{i_1, i_2, \cdots, i_r})
 \longrightarrow
\bigoplus\limits_{j=1}^r H_{2k-r-2} (Y_{i_1, i_2, \cdots, \wh i_j, \cdots, i_r}),
\]
be the homomorphism
\[
\eta_{i_1, i_2, \cdots, i_r} (\alpha)
=((\rho_{i_1, i_2, \cdots, i_r}^1)_* (\alpha), \cdots,
(\rho_{i_1, i_2, \cdots, i_r}^r)_* (\alpha)).
\]

\medskip
{\bf Proposition (4.11):}\ \
$\eta_{i_1, i_2, \cdots, i_r}$ sends
$H_{2k-r-2}(Y_{i_1, i_2, \cdots, i_r})$
injectively into $G_{i_1, i_2, \cdots, i_r}$,
which is a subgroup of $\ov G_{r-1}$.

\vskip 6pt
Let $H_r$ be the direct sum of homology groups
$H_{2k-r-2}(Y_{i_1, i_2, \cdots, i_r})$, for all
$(i_1, i_2, \cdots, i_r)$
satisfying $1 \le i_1 < i_2 < \cdots < i_r \le k$.
And
$\zeta_r : H_r \longrightarrow \ov G_{r-1}$ be the direct sum of
$\eta_{i_1, i_2, \cdots, i_r}$, that is,
$\zeta_r((\alpha_{I_1}, \alpha_{I_2}, \cdots, \alpha_{I_l}))$
$=$ $\sum\limits_{i=1}^l \eta_{I_i} (\alpha_{I_i})$,
where $I_i$ is some $(i_1, i_2, \cdots, i_r)$.
Thus, we conclude

\medskip
{\bf Corollary (4.12):}\ \
$\zeta_r : H_r \longrightarrow \ov G_{r-1}$ is injective.

\medskip
{\bf Proof of Proposition (4.11):}\ \
It is enough to check the following facts:
\begin{namelist}{(iii)}
\item[{\rm (i)}]
$\rho_{1,2, \cdots, r}^i \circ \varphi_{1, 2, \cdot, r}^i
=\psi_{1, 2, \cdots, \wh i, \cdots, r}^i$,
for $i \ge 2$.
\item[{\rm (ii)}]
$\rho_{1, 2, \cdots, r}^1 \circ \varphi_{1, 2, \cdots, r}^i
= \psi_{2, 3, \cdots, r}^1$, for $i \ge 2$.
(at least, represent the same homology class)
\item[{\rm (iii)}]
$\rho_{1, 2, \cdots, r}^j \circ \varphi_{1, 2, \cdots, r}^i$
is homotopic to a map which has image with less $1$ dimension in
$Y_{1, 2, \cdots, \wh j, \cdots, r}$, for $i \geq 2$ and
$j \neq 1,i$.
\end{namelist}

\vskip 6pt
Thus, $\eta_{1, 2, \cdots, r}(\varphi_{1, 2, \cdots, r}^i)
=(1, 0, \cdots, 1, 0, \cdots)$
in $G_{1,2,\cdots,r}$, for $i \ge 2$.

\vskip 6pt
$Im(\eta_{1,2,\cdots, r})$ has rank at least $r-1$.
But, rank $(H_{2k-r-2}(Y_{1,2, \cdots, r})) =r-1$.

\medskip
{\bf Theorem (4.13):}\ \
$H_{2k-3}(S(k)-F \big / S^1)=0$.

\medskip
{\bf Proof:}\ \
$S(k)-F \big / S^1 = Y_1 \cup Y_2 \cup \cdots \cup Y_k$
in $S(k) \big / S^1$.
By standard spectral sequence argument, the
following facts will implies
$H_{2k-3} (Y_1 \cup Y_2 \cup \cdots \cup Y_k)=0$.
\begin{namelist}{(iiii)}
\item[{\rm (i)}]
$H_{2k-3}(Y_i)=0$, for $i =1, 2, \cdots, k$.
\item[{\rm (ii)}]
$\bigoplus\limits_{1 \le i < j \le k}
H_{2k-4}(Y_i \cap Y_j)
\buildrel \zeta_2 \over \longrightarrow
\bigoplus\limits_{1 \le i \le k} H_{2k-4} (Y_i)$
is injective.
\item[{\rm (iii)}]
$\bigoplus\limits_{1 \le i \le j < l \le k}
H_{2k-5} (Y_i \cap Y_j \cap Y_l)
\buildrel \zeta_3 \over \longrightarrow
\bigoplus\limits_{1 \le i <j\le k}$ $H_{2k-5} (Y_i \cap Y_j)$
is injective.
\item[{\rm (iv)}] , $\cdots$, etc.
\end{namelist}

\noindent
where $\zeta_r$ is the direct sum of the natural maps,
for example $r=3$,
$\eta_{i,j,l}$ is the composite map
$$\matrix{
H_* (Y_i \cap Y_j \cap Y_l)
&\maprightu{}
&H_*(Y_i\cap Y_j) \bigoplus H_*(Y_i \cap Y_l)
\bigoplus H_* (Y_j \cap Y_l) \cr \cr
&&\mapdownr{} \cr\cr
&&\bigoplus\limits_{1 \le a < b \le k} H_* (Y_a \cap Y_b) \cr\cr}
$$
and $\zeta_3 = \sum\limits_{1 \le i <j < l \le k}
\eta_{i,j,l}$.
\newpage

\noindent{\Large \bf \S 5 \ Normal bundle of the fixed point set in
$\mbox{\hskip 1cm}$  $\WW_0$}

\vskip 6pt
$\WW_0(n)$ is a complicated C. W. complex, it is
not clear that a subspace in $\WW_0(n)$ has a
meaningful normal vector bundle.
But, the fixed point set $\HH(n)$ does have a $2n$-dimensional complex
vector bundle as its normal in $\WW_0(n)$.

\vskip 20pt
\Vjp
\noindent{\Large \bf  5.1 \ Inhomogeneous coordinate of $P(k)$}

\vskip 6pt
$\rho: \setC=\setR^2 \longrightarrow \setR \setP^2$,
$\rho(x,y)=[(x,y,1)]$, maps $\setC$ homeomorphically onto an open
set of $\setR\setP^2$.

\vskip 6pt
$\rho_k:\setC^k \longrightarrow P(k)$,
$P_k(\lambda_1, \lambda_2, \cdots, \lambda_k)=
(\rho(\lambda_1), \rho(\lambda_2), \cdots, \rho(\lambda_k))$,
maps $\setC^k$ homeomorphically onto an open dense set of
$P(k)=\prod\limits_{k}\setR\setP^2$.
$\xi_0$ denotes the fixed point of the
$S^1$-action on $P(k)$ which is defined in section 4.1.
$\rho_k(\setC^k)$ is a contractible open neighborhood of $\xi_0$.

\vskip 6pt
Let $U_k$ denote the open neighborhood
$\rho_k(\setC^k)$ of $\xi_0$.

\vskip 6pt
Suppose $\Gamma$ is a knot graph. As in section 3.1,
$W_0(\Gamma)$ is the configuration space of
infinitesimal knot graph
on the line $l_0$ of $z$-axis,
$\Psi : W_0(\Gamma) \longrightarrow P(k)$ is a canonical
$P(k)$-structure.

\vskip 6pt
Suppose $f : \Gamma \longrightarrow \Gamma'$
is an equivalence of knot graphs and $\Gamma'$ is a knot graph
on $l_0$.
For any $v \in V(\Gamma)$, $f(v)$ is a point of $\setR^3$.
Identify $\setR^3$ as $\setC \times \setR$,
and write the function
$f: V(\Gamma) \longrightarrow \setR^3$ as $(g,h)$,
$g:V(\Gamma) \longrightarrow \setC$ and
$h : V(\Gamma) \longrightarrow \setR$.
Thus, the function $g$, satisfying
$g |_{V_0(\Gamma)} =0$, is called a ground function of $\Gamma$ and the
function $h$, preserving the linear order of base points in
$V_0(\Gamma)$, is called a height function of $\Gamma$.
A ground function $g$ of $\Gamma$ together with a height function
$h$ of $\Gamma$ forms a knot graph on $l_0$.
But a height function could be thought as a knot graph whose vertices
are all on the line $l_0$.

\vskip 6pt
Suppose $\Gamma$ has $s$ inner vertices.
$C(\Gamma, l_0)$ is "almost" an $s$-dimensional complex
vector bundle over the space of all height function, and
$W_0(\Gamma)$ is also "almost" an $s$-dimensional complex vector
bundle over the space $\{$ all height function of $\Gamma\}$
$\big/$ translation and dilation relation.
The relation is the following:

\subitem
$h_1, h_2$ are two height function of $\Gamma$, \hfill\break
$h_1 \buildrel {\rm t.d.} \over \sim h_2$,
if there exist $\lambda > 0$ and real number $t$
such that $h_1(v)=\lambda h_2 (v) + t$.

\noindent
But, we are more interested in the knot graph
$f = (g,h)$ on $l_0$ which is in
$\Psi^{-1} (U_k)$. Why?

\vskip 6pt
Suppose $\Gamma$ has edges $e_1= \{v_1, w_1\}$,
$e_2 = \{v_2, w_2\}$, $\cdots$, $e_k=\{v_k, w_k\}$.
Then
\begin{eqnarray*}
& & \Psi (f) (\mbox{ more formally  } \Psi(f, f(\Gamma)))  \\
& & =([f(v_1)-f(w_1)], [f(v_2)-f(w_2)], \cdots,
[f(v_k) - f(w_k)]) .
\end{eqnarray*}
When $\Psi(f)$ is in $U_k=\rho_k(\setC^k)$,

\[
\rho_k^{-1} (\Psi (f))
=\left(
\frac{g(v_1)-g(w_1)}{h(v_1)-h(w_1)},
\frac{g(v_2)-g(w_2)}{h(v_2)-h(w_2)},
\cdots,
\frac{g(v_k)-g(w_k)}{h(v_k)-h(w_k)} \right).
\]

\vskip 6pt
Thus, in the coordinate system
$\rho_k:\setC^k \longrightarrow P(k)$,
the $P(k)$-structure
$\Psi:W_0(\Gamma) \longrightarrow P(k)$ is an
injective fibre-wise linear map.

\vskip 6pt
To formulate the statement above  precisely,
let us consider the following definitions.

\medskip
{\bf Definition (5.1):}\ \
Suppose $g_i : V(\Gamma) \longrightarrow \setC$
are ground functions and
$h_i : V(\Gamma) \longrightarrow \setR$ are height functions,
$i =1,2$.
$(g_1, h_1)$ is equivalent to $(g_2, h_2)$ under the translation
and dilation relation, if
there exist real number $\lambda$ and $t$,
$\lambda >0$, such that $g_1(v)=\lambda g_2 (v)$ and
$h_1(v) = \lambda h_2(v)+t$
for all $v \in V(\Gamma)$.

\vskip 6pt
Thus, $W_0(\Gamma)$ is contained in and "almost" equal to the space
$\{(g,h): g$ is a ground function of $\Gamma$
and $h$ is a height function of $\Gamma\}$ $\big /$
T.D. relation.
(Recall: $g |_{V_0(\Gamma)}=0$ and
$h$ preserves the linear order of base points, that is,
if $x_1 < x_2$ in $V_0(\Gamma)$,
then $h(x_1) < h(x_2)$.)

\medskip
{\bf Definition (5.2):}\ \
Suppose $\Gamma$ is a knot graph.
$D(\Gamma)$ denotes the subspace of $W_0(\Gamma)$,
$\{(g,h) \in W_0(\Gamma): h(v_i) \ne h(w_i)$,
for all $e_i=\{v_i, w_i\}$
in $E(\Gamma)\}$ $\big /$ T.D. relation,
and $H(\Gamma)=\{h: V(\Gamma) \longrightarrow \setR : h(x) < h(x')$,
for any two base points $x <x'$ in $V_0(\Gamma)$, and
$h(v_i) \ne h(w_i)$,
for all edge $e_i=\{v_i, w_i\}$ in $E(\Gamma)\}$ $\big /$ T.D. relation.
Then $0 \times H(\Gamma)=\{(0,h):h \in H(\Gamma)\}$
is a subspace of $D(\Gamma)$, also, a subspace of $W_0(\Gamma)$.

\medskip
{\bf Proposition (5.3):}\ \
\begin{namelist}{(iii)}
\item[{\rm (i)}]
$D(\Gamma) = \Psi^{-1} (U_k)$
\item[{\rm (ii)}]
$0 \times H(\Gamma)=\Psi^{-1}(\xi_0)$, it is the set of all fixed
points of $S^1$-action on $W_0(\Gamma)$ given in section 3.3.
\item[{\rm (iii)}]
Assume $\Gamma$ has $s$ inner vertices. $D(\Gamma)$ is an
$s$-dimensional complex vector bundle over $H(\Gamma)$.
\hfill\break
(Proof is straightforward and is omitted)
\end{namelist}

\vskip 6pt
For convenience, we identify $H(\Gamma)$ as
$0 \times H(\Gamma)$.
Thus, a height function is also a knot graph on $l_0$.
But, there is no way to think a ground function as a knot graph
on $l_0$.

\vskip 6pt
Let $D(\Gamma)_h = \{(g,h):(g,h)\in D(\Gamma)\}$, it
is the $s$-dimensional complex vector space over $h$.
The vector bundle $D(\Gamma)$ over $H(\Gamma)$ is trivial,
but there is no canonical basis.

\medskip
{\bf Proposition (5.4):} \ \
Suppose $\Gamma$ is a knot graph and $h$ is a height function in $H(\Gamma)$.
Then, the restriction $\Psi_h$ of
$\rho_k^{-1} \circ \Psi$
to $D(\Gamma)_h$ is an injective linear map from
$D(\Gamma)_h$ to $\setC^k$.
(Note: order $(\Gamma) = k-s$ is always positive.)

(Rigorous consideration of injectivity is in section 6.2 and 6.3 . )

\medskip
{\bf Question:} \ \
Could $\Psi_h:D(\Gamma)_h \longrightarrow \setC^k$
continuously extend to the boundary of $H(\Gamma)$?
(What is the boundary of $H(\Gamma)$?)

\vskip  6pt
To answer the question, there are still many works to do.
For convenience, we abuse the notations: to denote
$\rho_k^{-1} \circ \Psi$ also by $\Psi$.
Then
\[
\Psi(g,h)=\left(
\frac{g}{h} (e_1),
\frac{g}{h} (e_2), \cdots,
\frac{g}{h} (e_k) \right),
\]
where
\[
\frac{g}{h} (e_i)=\frac{g(v_i)-g(w_i)}{h(v_i)-h(w_i)},
\mbox{ for } \,
e_i=\{v_i, w_i\}.
\]
And the fibrewise linear map associated with
$W_0(\Gamma) \times P(r)$ is $\Psi \times id_{\setC^r}$
$:$ $D(\Gamma) \times \setC^r \longrightarrow \setC^{k+r}$.
And we also denote the stable fibrewise linear map by $\Psi$.
For any $\sigma$ in $\Sigma_{k+r}$,
$\sigma : \setC^{k+r} \longrightarrow \setC^{k+r}$,
$\sigma(\lambda_1, \lambda_2, \cdots, \lambda_{k+r})
=(\lambda_{\sigma(1)},
\lambda_{\sigma(2)}, \cdots,
\lambda_{\sigma(k+r)})$,
is a linear isomorphism.
Thus, $\Sigma_{k+r} \cdot \Psi$ will consists of all
the fibrewise linear maps associated with all the canonical
$P(k+r)$-structure of $W_0(\Gamma) \times P(r)$.
Similar to the construction of $\WW_0(n)$, let
$\DD(n)$ be the disjoint union of
$D(\Gamma) \times \setC^{3n-k} \times \Sigma_{3n} \cdot \Psi$,
for all normal knot graphs $\Gamma$ with order $n$,
and let $\ov \HH(n)$ be the disjoint union of
$H(\Gamma) \times \Sigma_{3n} \cdot \Psi$,
for all normal knot graph $\Gamma$ with order $n$
(note: the function $\Psi$, multiplied with $H(\Gamma)$,
is the same as the function $\Psi$, multiplied with
$D(\Gamma) \times \setC^{3n-k}$).
Then $\ov \DD(n)$ is a $2n$-dimensional complex vector
bundle over $\ov \HH(n)$.
Now, we need identify the components in $\ov \DD(n)$
and $\ov \HH(n)$ along their boundaries.

\vskip 6pt
Their boundaries are similar to the
boundaries of configuration spaces $C(\Gamma)$.
With the boundaries, $D(\Gamma)$ will be substituted by its
base-compactification and so are $H(\Gamma)$.

\vskip 20pt
\Vjp
\noindent{\Large \bf  5.2 \ Boundary of $W_0(\Gamma)$}

\vskip 6pt
There is a special identification in the boundary of $W_0(\Gamma)$
we should do prior to the other identifications.

\vskip 6pt
Suppose $A(\Gamma)$ is the disjoint union of connected
components
$A_1(\Gamma), A_2(\Gamma), \cdots$,
and $A_m(\Gamma)$, where $A$ is the disjoint union
of $A_1, A_2, \cdots$, and $A_m$. The associated
extended translation and dilation relation is called a
special extended translation and dilation relation,
or simply, a S.E.T.D. relation associated to $(A, \Gamma)$.
(see (3.5) in section 3.1)

\medskip
{\bf Notation (5.5):} \ \
\begin{namelist}{(ii)}
\item[{\rm (i)}]
$Q(A(\Gamma))$ is the space $W_0(A(\Gamma))$
quotiented by the S.E.T.D. relation associated with $(A, \Gamma)$.
\item[{\rm (ii)}]
$Q(\Gamma; A)=W_0(\Gamma \big / A) \times Q(A(\Gamma))$.
\end{namelist}

\medskip
{\bf Remark:} \ \
$Q(A(\Gamma))$ is exactly a subspace of $W_0(A_1(\Gamma))
\times W_0(A_2(\Gamma))$
$\times \cdots \times$ $W_0(A_m(\Gamma))$,
and $Q(\Gamma; A)$ is the more qualified boundary
than $W_0(\Gamma; A)$,
because of the following proposition.

\medskip
{\bf Proposition (5.6):}\ \
Outside some odd part,
$W_0(\Gamma) \bigcup\limits_{A} Q(\Gamma; A)$
with the natural topology is an $s$-dimensional complex vector
bundle over its fixed point set
$\Fix (W_0(\Gamma) \bigcup\limits_{A} Q(\Gamma; A))$.
(In this proposition, $W_0(\Gamma)$ return to its original
space without compactification and $W_0(\Gamma) \bigcup\limits_{A}
Q(\Gamma; A)$ is a locally compact space containing
$W_0(\Gamma)$.)

\medskip
{\bf Proof:}\ \
Let $\Fix (\Gamma)$ denote the space
$\Fix (W_0(\Gamma) \bigcup\limits_{A} Q(\Gamma; A))$.

\vskip 6pt
We first define the bundle projection
\begin{eqnarray*}
& & \pi: (W_0(\Gamma) \bigcup\limits_{A} Q(\Gamma; A)
\longrightarrow \Fix (\Gamma), \\
& & (g,h) \in W_0(\Gamma), \pi(g,h)=(0,h),
\end{eqnarray*}
similarly, for
\[
f_0=(g_0, h_0) \in W_0(\Gamma \big / A)
\]
and
\begin{eqnarray*}
& & f_i=(g_i, h_i) \in W_0(A_i(\Gamma)), \, \, i=1, 2, \cdots, m, \\
& & \alpha = (f_0, f_1, f_2, \cdots, f_m) \in Q(\Gamma; A), \\
& &\pi(\alpha)=((0, h_0), (0, h_1), \cdots, (0, h_m)).
\end{eqnarray*}
Then
$\beta=((0, h_0), (0, h_1), \cdots, (0, h_m))$ is in
$\Fix (\Gamma)$.

\vskip 6pt
But $\pi^{-1} (\beta)$ may contains more elements other than
that in $W_0(\Gamma; A)$.
If $f=(g,h) \in W_0(\Gamma)$ such that
$g |_{A_i}$ is constant, for each $i=1, 2, \cdots, m$
and $h =h_0$, then $\pi(f, (g_1, h_1)$, $\cdots$, $(g_m, h_m))$ is also
equal to $\beta$.
Although $\gamma =(f$, $(g_1, h_1)$, $\cdots$, $(g_m, h_m))$ is
not in the interior of a codimension $1$ boundary
$W_0(\Gamma; A')$, for some $A'$,
$\gamma$ is in the boundary associated to the sets
$(A_1, A_2, \cdots, A_m)$.
Thus, $\pi^{-1}(\beta)$ is the vector space
$G=\{ ((g, h_0)$, $(g_1, h_1)$, $\cdots$, $(g_m, h_m))$ $:$
$g: V(\Gamma)$ $\longrightarrow$ $\setC$, $g |_{A_i}$ is constant,
for each $i=1,2, \cdots, m$,
and $g_i : A_i \longrightarrow \setC$, $i=1,2,\cdots, m\}$
quotiented by the following relations:

\subitem
$\gamma=((g, h_0), (g_1, h_1), \cdots, (g_m, h_m))$
and
$\gamma'=((g', h_0)$, $(g'_1, h_1)$, $\cdots$, $(g'_m, h_m))$
are two elements in $G$, \hfill\break
$\gamma \sim \gamma'$, if there exist complex numbers
$c_1, c_2, \cdots, c_m$ such that
$g_i = g'_i +c_i$,
for each $i=1,2, \cdot, m$.

\noindent
Thus, $\pi^{-1} (\beta)$ has a natural vector space structure
and has $s$ dimension.

\vskip 6pt
And the set $G'(h_0) =\{(g, h_0) \in W_0(\Gamma)$,
for some $i$, $g |_{A_i}$ is not constant $\}$
which contains elements in
$\pi^{-1}$ $((0, h_0))$ but not
$\pi^{-1}$ $((0, h_0)$, $(0, h_1)$, $\cdots$, $(0, h_m))$,
for any $h_1, h_2, \cdots, h_m$,
$h_i : A_i \longrightarrow \setR$,
it is the odd part of $W_0(\Gamma) \bigcup\limits_{A} Q(\Gamma; A)$ in the
fibre over $\beta$.
This proves (5.6).

\vskip 6pt
Remarks on the proof of (5.6):
\begin{namelist}{(iiii)}
\item[{\rm (i)}]
In $W_0(\Gamma) \bigcup\limits_{A} W_0(\Gamma; A)$, if
$g(A_i)=g(A_j)$, for some $i \ne j$,
there is an unnecessary compactification in
$W_0(\Gamma; A)$, which destroies the natural vector space structure.
\item[{\rm (ii)}]
Consider the following commutative diagram
$$\matrix{
Q(\Gamma; A) & \maprightu{q'} & W_0(\Gamma \big / A) \cr\cr
\mapdownr{\pi} &&\mapdownr{\pi'} \cr\cr
\Fix(Q(\Gamma; A)) & \maprightu{q} & \Fix(W_0(\Gamma \big / A)) \cr \cr
\beta &\maprightu{} &h_0 \cr}
$$
where $q(h_0, h_1, \cdots, h_m)=h_0$ and
$q'(f_0, f_1, f_2, \cdots, f_m)=f_0$.
Although $q$ is onto obviously,
$\pi^{-1}(h_0)$ is larger than
$\pi^{-1}(q^{-1}(h_0))=q^{'-1}(\pi^{'-1}(h_0))$.
And the odd part $G'(h_0)$ is equal to $\pi^{-1}(h_0)- \pi^{-1}(q^{-1}(h_0))$,
which is contained in $W_0(\Gamma)$ and not in $Q(\Gamma; A)$.
\item[{\rm (iii)}]
Suppose
$\Psi : W_0(\Gamma) \bigcup\limits_{A} Q(\Gamma; A) \longrightarrow P(k)$
is a canonical $P(k)$-structure.
The complement of $\Psi^{-1} (U_k)$ is exactly the odd part.
\item[{\rm (iv)}]
With some abuse of notation, we may denote
$W_0(\Gamma) \bigcup\limits_{A} Q(\Gamma; A)$ by $\ov W_0(\Gamma)$,
the odd part of $\ov W_0(\Gamma)$ by $Odd(W_0(\Gamma))$, and
$\ov W_0(\Gamma)$ $-$ $Odd(W_0(\Gamma))$ by $D(\Gamma)$.
Then $\Psi : D(\Gamma) \longrightarrow \setC^k$ is the
fibre-wise injective linear map, which has a simple form as follows:

\vskip 6pt
(Case 1)\ \ Interior of $W_0(\Gamma)$.
\begin{eqnarray*}
& & f = (g,h), \\
& & \Psi(f)=\left( \frac{g}{h} (e_1), \frac{g}{h}(e_2), \cdots,
\frac{g}{h}(e_k) \right)
\end{eqnarray*}

\vskip 6pt
(Case 2) \ \ Codimension $1$ boundary ($A(\Gamma)$ has only one connected
component).

$\alpha =(f_0, f_1)$, $f_0=(g_0, h_0)$,
$g_0 |_A$ and $h_0 |_A$ are constant,
$f_1 = (g_1, h_1)$,
$g_1 : A \longrightarrow \setC$,
$h_1 : A \longrightarrow \setR$.
Assume $e_1, e_2, \cdots, e_r$ are in $A(\Gamma)$.
\[
\Psi (\alpha) = \left(
\frac{g_1}{h_1} (e_1),
\frac{g_1}{h_1} (e_2), \cdots,
\frac{g_1}{h_1} (e_r),
\frac{g_0}{h_0} (e_{r+1}), \cdots,
\frac{g_0}{h_0} (e_k) \right)
\]

\vskip 6pt
(Case 3) \ \ $A(\Gamma)$ has more than $1$ connected components.

$\Psi$ has similar form as in (case 2).

(Note: (i) if $e =\{v,w\}$,
$\frac{g}{h} = \frac{g(v) - g(w)}{h(v)-h(w)}$,
(ii) $\Psi$ is, in fact, $\rho_k^{-1} \circ \Psi$.)
\end{namelist}

\vskip 20pt
\Vjp
\noindent{\Large \bf  5.3 \ Identification maps of vector bundles}

\vskip 6pt
A vector bundle map is a (bundle) map between two vector
bundles, which is fibrewise linear isomorphism. Now,
we consider all the "identification" vector bundle maps come from the
identification map in section 1.3 except the identification maps of type 0,
which has been done in the last section (5.2).

\Vjp
\noindent{\large \bf \ 5.3.1 \ Identification map of type II}

\vskip 6pt
Suppose $A(\Gamma)$ has a bivalent inner vertice $v$,
$e_1 = \{v, w, \}$ and $e_2 =\{v, w_2\}$ are the two edge connecting
to $v$. Without loss of generality, we may assume
$A(\Gamma)$ is connected. Thus
$Q(\Gamma; A)= W_0(\Gamma; A)$.

\vskip 6pt
Assume $f =(g, h) \in Q(\Gamma; A)$.
Then, the identification map $\tau_2$ sending
$f$ to $\tau_2(f)=$
$(\tau'_2(g), {\tau''}_2(h))$
\begin{eqnarray*}
& & \tau_2(f) (v) = f(w_1)+f(w_2) - f(v) \\
& & \tau(f) (v')=f(v'), \, \, \mbox{ if } \, \,
v' \ne v.
\end{eqnarray*}
Thus,
\begin{eqnarray*}
& & \tau'_2(g) (v) = g(w_1)+g(w_2) -g(v), \\
& & \tau'_2 (g) (v') = g(v'), \, \, \mbox{ if } \,
v' \ne v ;
\end{eqnarray*}
and
\begin{eqnarray*}
& & \tau^{''}_2 (h)(v) = h(w_1) + h(w_2) - h(v), \\
& & \tau^{''}_2 (h)(v') = h(v'), \, \, \mbox{ if } \,  v' \ne v.
\end{eqnarray*}
Therefore, $\tau'_2: D(\Gamma; A)_h \longrightarrow D(\Gamma; A)_{\tau^{''}_2 (h)}$
is a linear isomorphism of fibres,
where
\begin{eqnarray*}
  D(\Gamma; A) & & = \Psi^{-1} (U_k) \cap Q(\Gamma; A), \\
\mbox{also}& & =D(\Gamma) \cap Q(\Gamma; A).
\end{eqnarray*}

\Vjp
\noindent{\large \bf \ 5.3.2 \ Identification map of type III}

\vskip 6pt
Suppose $A=\{v,w\}$ is an edge of $\Gamma$, and $v$ is an inner vertice.
$\tau_3 : W_0(\Gamma; A) \longrightarrow W_0(\Gamma \big / A) \times
\setR\setP^2$ is the identification map defined in section 1.3.4.

Assume $\alpha =(f_0, f_1) \in W_0(\Gamma \big /  A) \times Q(A(\Gamma))$,
$f_1 = (g_1, h_1)$,
$g_1 : A \longrightarrow \setC$,
$h_1 : A\longrightarrow \setR$,
$f_0 \in W_0(\Gamma \big / A)$.
Then
\[
\tau_3 (\alpha) = \left(f_0, \rho \left(\frac{g_1(v)-g_1(w)}{h_1(v)-h_1(w)}
\right)\right),
\]
$\rho: \setC \longrightarrow \setR\setP^2$ is defined in section 5.1.
If $f_0=(g_0,h_0)$,
$g_0 : V(\Gamma \big / A) \longrightarrow \setC$,
$h_0: V(\Gamma \big / A) \longrightarrow \setR$.
Then
\[
\tau_3 : D(\Gamma; A)_{(h_0, h_1)} \longrightarrow D(\Gamma \big / A)_{h_0}
\times \setC
\]
is the following
\[
\tau_3 ((g_0, h_0), (g_1, h_1)) =\left((g_0, h_0),
\frac{g_1(v) - g_1(w)}{h_1(v) - h_1(w)} \right),
\]
which is an linear isomorphism.

\vskip 6pt
(Note: when $\alpha \in D(\Gamma; a) \subset \Psi^{-1} (U_k)$,
$h_1 (v)-h_1(w)$ cannot be zero.)

\Vjp
\noindent{\large \bf \ 5.3.3 \ Identification map of type IV}

\vskip 6pt
In the case, $\tau_4$ is obviously a linear isomorphism.

\Vjp
\noindent{\large \bf \ 5.3.4 \ Identification of type I}

\vskip 6pt
{\bf (i) Alternating method:}

\vskip 6pt
Suppose $\Gamma$ is a normal knot graph,
$|A| \ge 3$, $A(\Gamma)$ has a univalent inner vertice $v$,
$e=\{v, v_1\}$ is the unique edge connecting to $v$.

\vskip 6pt
(Case 1)  \ \ $v_1$ is also a univalent inner vertice in $A(\Gamma)$.
Then $A(\Gamma)$ is disconnected. The special translation and dilation
relation for $Q(\Gamma; A)$ has reduced $Q(\Gamma; A)$ to lower dimensional
boundary.

\vskip 6pt
(Case 2) \ \
$v_1$ is a bivalent inner vertice in $A(\Gamma)$.
Then we apply the identification map of type II.

\vskip 6pt
(Case 3) \ \ $v_1$ is a trivalent inner vertice in $A(\Gamma)$,
that is, there exist edges
$e_1 =\{v_1, w_1\}$ and
$e_2=\{v_1, w_2\}$ connecting to $v_1$.
Considering the following map, which was first given by Bott and
Taubes [5]:
\begin{eqnarray*}
& & \tau : Q(A(\Gamma)) \longrightarrow Q(A(\Gamma)) \\
& & f:A \longrightarrow \setC \times \setR, \, \, f \in Q(A(\Gamma)). \\
& & \tau(f) (v)=f(w_1)+f(w_2) - f(v) \\
& & \tau(f)(v) = f(w_1) + f(w_2) - f(v) \\
& & \tau(f)(v_1) = f(w_1) + f(w_2) - f(v_1) \\
& & \tau(f)(v') =  f(v') , \, \, \mbox{ if } \, \,
v' \ne v
\mbox{ and } v' \ne v_1 .
\end{eqnarray*}

\vskip 6pt
Thus, we have an identification map similar to the map of type II.
Similarly, we have the linear isomorphism of fibres.
This identification is said to be of the type Y.

\vskip 6pt
(Case 4) \ \ $v_1$ is a base points in $\Gamma$. We may assume
$\Gamma$ has univalent base point
only. And, it is similar to the
(Case 1).

\medskip
{\bf (ii) Original method:}

\vskip 6pt
Following the assumption in (i), $v$ is a univalent
inner vertice of $A(\Gamma)$ and $e =\{v, v_1\}$
is an edge in $A(\Gamma)$.

\vskip 6pt
For convenience, assume $A(\Gamma)$ is connected, then
$Q(A(\Gamma))=W_0(A(\Gamma))$.

\vskip 6pt
In section 1.3.2, we have the identification map
$\tau_1 :W_0(A(\Gamma)) \longrightarrow W_0(A(\Gamma))$:
for $f :A \longrightarrow \setC \times \setR$,
let
$f_1 = f|_{A-\{v\}}$ and $||f_1||=\max\{|f_1 (w)-f(w') |$,
$w, w' \in A-\{v\} \}$,
\begin{eqnarray*}
& & \tau_1(f) (v) = f(v_1)+ 2 || f_1 || \frac{f(v)-f(v_1)}{|f(v)-f(v_1)|}, \\
& & \tau_1(f) (w) = f(w), \, \, \mbox{ for } \,
w \ne v .
\end{eqnarray*}

\vskip 6pt
It is easy to see that $\tau_1$ is $S^1$-equivariant.
Thus, $\tau_1$ sends the fixed points to fixed points.
Write $f$ as $(g,h)$,
$\pi(f)=h$.
But $\pi (\tau_1 (f))$ may not be $\tau_1(\pi(f))$.
Let $\tau'_1 : D(\Gamma; A) \longrightarrow D(\Gamma; A)$
be the map defined as follows:
\begin{eqnarray*}
& & \alpha =((g_0, h_0), (g, h) ) \in D(\Gamma; A), \\
& & g_0: \Gamma \big / A \longrightarrow \setC, \, \,
h_0: \Gamma \big / A \longrightarrow \setR,  \\
& & g:  A \longrightarrow \setC,  \, \,
h_0: A \longrightarrow \setR,
\end{eqnarray*}
At first, let $g_1 = g|_{A-\{v\}}$,
$h_1 =h|_{A-\{v\}}$.
\begin{eqnarray*}
& & \tau'_1 (g) (v) = g(v_1) + 2 || h_1 ||
\frac{g(v)-g(v_1)}{|h(v)-h(v_1)|} ,  \\
& &\tau'_1 (g) (w) = g(w) , \, \mbox{ for } \, w \ne v , \\
& & \tau'_1 (h) (v) = h(v_1) + 2 || h_1 ||
\frac{h(v)-h(v_1)}{|h(v)-h(v_1)|} ,  \\
& &\tau'_1 (h) (w) = h(w) , \, \mbox{ for } \, w \ne v ,
\end{eqnarray*}
and $\tau'_1 (\alpha) =\left( (g_0, h_0), (\tau'_1 (g), \tau'_1(h) \right)$.
It is easy to show that
\begin{namelist}{(iiii)}
\item[{\rm(i)}]
when $\alpha \in \Fix (W_0(\Gamma; A))$,
$\tau'_1(\alpha)=\tau_1(\alpha)$.
\item[{\rm(ii)}]
$\tau'_1$ is a vector bundle map.
\item[{\rm(iii)}]
$\tau'_1$ is the differential of $\tau_1$ at the fixed point set.
\item[{\rm(iv)}]
$\Psi (\tau_1(\alpha))=\Psi(\tau'_1 (\alpha))=\Psi(\alpha)$.
\end{namelist}

\vskip 6pt
And the most important fact is that the two identification maps
$\tau_1 |_{D(\Gamma ; A)}$ and
$\tau'_1$ give the same result.
Thus, we may use $\tau'_1$ instead of $\tau_1$, when we consider only the part
$D(\Gamma; A)$, that is,
$\Psi^{-1} (U_k) \cap Q(\Gamma; A)$.

\medskip
{\bf Remark:}\ \ The alternating method for the ``univalent inner
vertice'' identification is just reducing the type I to the other
types. The merit of this method is that the identification maps for
$\ov{W_0(\Gamma)}$ $(= W_0(\Gamma) \bigcup\limits_{A} Q(\Gamma; A))$
are all local diffeomorphisms for both vector bundles
and base spaces
(the fixed point set of $Q(\Gamma; A)$)

\vskip 6pt
First kind: \ \ $(|A| \ge 3)$

(Type V): \ $\tau_2 :Q(\Gamma; A) \longrightarrow Q(\Gamma; A)$

(Type Y): \ $\tau   :Q(\Gamma; A) \longrightarrow Q(\Gamma; A)$

\vskip 6pt
Second kind: \ \ $(|A| = 2)$

$\tau_3 , \tau_4 :Q(\Gamma; A) \longrightarrow W_0(\Gamma \big / A) \times P(A)$.

\vskip 6pt
Furthermore, the two different type identification maps of first
kind are $P$-orientation reversing involutions.
This is a crucial point in the degree theory.

\vskip 20pt
\Vjp
\noindent{\Large \bf  5.4 \ Conclusion}

\vskip 6pt
We continue the notations $\ov \DD(n)$ and
$\ov \HH(n)$ in section 5.1.
Let $\DD(n)$ be the quotient space of $\ov{\DD(n)}$ by the
identification maps of vector bundles given in lase section (5.3),
and the fixed point set $\HH(n)$ is also a quotient space of
$\ov \HH(n)$ by the related identification maps.
Then $\DD(n)$ is a $2n$-dimensional complex vector bundle over $\HH(n)$.

\vskip 6pt
\noindent {\Large Overall Conclusion}

First, we state a few results:

\begin{eqnarray*}
{\bf (5.7)} \quad
& & H_{6n}(\WW(n), \Psi^{\*} (\OO_{3n}))   \\
& & \approx H_{6n-2} (\WW_0(n), \Psi^{\#} (\OO_{3n})) \\
& & \approx H_{6n-3} (\WW_0(n) \big / S^1, \ov \Psi^{\#} (\ov \OO_{3n})) \\
& & \mbox{(The proof is somewhat straightforward.)}
\end{eqnarray*}

Thus, the homomorphism

$\beta : H_{6n-3} (\WW_0(n) \big / S^1, \ov \Psi^{\#} (\ov \OO_{3n}))
\longrightarrow H_{6n} (\WW_0(n), \Psi^{\#} (\OO_{3n})$, \\
defined in section 4.2, is an isomorphism, and hence,
together with Theorem (4.2), we have

\medskip
{\bf Theorem (5.8):}\ \
If $\buildrel = \over \Psi_* : H_{6n-4} (\DD(n) - \HH(n) \big / S^1)
\longrightarrow H_{6n-4}(\setC \setP^{3n-1})$
is a zero-homomorphism, then the degree homomorphism  \\
$\Psi_* : H_{6n} (\WW(n), \Psi^* (\OO_{3n})) \longrightarrow
H_{6n} (P(3n), \OO_{3n})$ is also
a $0$-map.

\vskip 6pt
In Chapter 6, we shall show that $\DD(n)$ has finite structure group.
Thus,
$\buildrel =\over \Psi_*$ in Theorem (5.8) is a
$0$-homomorphism, and hence, the degree homomorphism $\Psi_*$
is also trivial.

\vskip 6pt
Secondly, we consider the set
$\CC=\WW_0(n) - \DD(n)$.
Then $\WW_0(n) \big / \CC$ is exactly the Thom space
$T(\DD(n))$ of the vector bundle $\DD(n)$.
Thus, we have the following exact sequence

$0 = H_{6n-2}(\CC, \Psi^* (\OO_{3n})) \longrightarrow
H_{6n-2}(\WW_0(n), \Psi^{\#} (\OO_{3n}))$ $ \longrightarrow
H_{6n-2}(T(\DD(n))$  $ \longrightarrow
H_{6n-3}(\CC, \Psi^{\#} (\OO_{3n}))$.

By Thom isomorphism Theorem,
$$
H_{6n-2} (T (\DD(n)) \approx H_{2n-2} (\HH(n)).
$$
Therefore, we have

\medskip
{\bf (5.9):}\ \
$\Delta_* : H_{6n-2} (\WW_0(n), \Psi^{\#} (\OO_{3n})) \longrightarrow H_{2n-2} (\HH(n))$
is an injective map.

\vskip 6pt
With a straightforward, but tedious, computation of $H_*(\HH(n))$,
we can show that the image of $\Delta_*$ contains exactly
the cycles from the weight system of Vassiliev invariant of order $n$.
On the other hand, one can convince himself easily that a weight
system $\omega$ will get cycles
$\sum\limits_{\Gamma} \omega (\Gamma)W_0(\Gamma)$ in
$(\WW_0(n), \Psi^{\#} (\OO_{3n}))$ and
$\sum\limits_{\Gamma} \omega(\Gamma) W(\Gamma)$
in $(\WW(n), \Psi^{\#} (\OO_{3n}))$.
And the degree homomorphism send it to $0$,
which is a non-zero multiple of Feynman integral over
$\sum\limits_{\Gamma} \omega (\Gamma) W(\Gamma)$.
This concludes that the Feynman  integral
$\sum\limits_{\Gamma} \frac{f_\Gamma}{|\Gamma|} \omega (\Gamma)$
is $0$,
for any weight system, and hence,
$\sum\limits_{\Gamma} \frac{f_\Gamma}{|\Gamma|} [\Gamma]$
is $0$ in the algebra of chord diagram.

\medskip
{\bf Remark:}\ \ The computation of $H_*(\HH(n))$ will
appear in additional chapter (seven) or another article.
Although it is lengthy, it has no
problem anyway.
\newpage

\noindent{\Large \bf \S 6  \ Finite structure group}

\vskip 20pt
\Vjp
\noindent{\Large \bf  6.1 \ Finite structure group}

\vskip 6pt
Suppose $\Gamma$ is a knot graph.
$\ov H(\Gamma)$ is the set of all real function $h$ on $V(\Gamma)$
such that $h$ preserves the linear order of base points of $V(\Gamma)$ and $h$
does not degenerates any edge in $V(\Gamma)$,
precisely, for any two base points $x$ and $x'$ with
$x <x'$, $h(x) < h(x')$, and for any edge $e =\{v,w\}$ in $\Gamma$,
$h(v) \ne h(w)$.
$\ov D(\Gamma)=\{(g,h):h \in \ov H(\Gamma)$
and $g$ is a complex-value function on $V(\Gamma)$
such that $g |_{V_0(\Gamma)}=0\}$.
$H(\Gamma)$ and $D(\Gamma)$ are the quotient spaces of
$\ov H(\Gamma)$ and $\ov D(\Gamma)$ by all possible extended
translation and dilation relations.

\vskip 6pt
Suppose $\Gamma$ has $s$ inner vertices.
Then $D(\Gamma)$ is an $s$-dimensional complex vector bundle over
$H(\Gamma)$.

\vskip 6pt
Assume $\Gamma$ has $k$ edges
$e_1,e_2, \cdots, e_k$. Then, we have a fibrewise linear
injective map
$\Psi : D(\Gamma) \longrightarrow \setC^k$,
defined by
\[
\Psi (g,h)=\left( \frac{g}{h} (e_1), \frac{g}{h} (e_2), \cdots,
\frac{g}{h} (e_k) \right ),
\]
which depends on an order of edges
$\bigg ( \frac{g}{h} (e)=\frac{g(v)-g(w)}{h(v)-h(w)}$,
if $e = \{v,w\} \bigg)$.
The E.T.D. relations are exactly the a priori relations having the
same image under $\Psi$.

\vskip 20pt
\Vjp
\noindent{\Large \bf  6.2 \ A trivialization of $D(\Gamma)$}

\vskip 6pt
Fix a linear order for the inner vertices in $\Gamma$, say,
$y_1 < y_2 < \cdots < y_s$, for
$y_1, y_2, \cdots, y_s$ in $V_1(\Gamma)$.

\vskip 6pt
Suppose $h : V(\Gamma) \longrightarrow \setR$ is a height function.
For any edge $e=\{v,w\}$, let $| h |(e)=|h(v)-h(w)|$.
A sequence of vertice $\eta=(v_0, v_1, \cdots, v_r)$ is
said to be an arc connecting $v_0$ to $v_r$, if
$\{v_{i-1}, v_i\}$, $i=1,2,\cdots, r$, are edges in $\Gamma$.
Let $|h |(\eta)=\sum\limits_{i=1}^r |h|(v_{i-1},v_i)$,
it is said to be the $h$-length of the arc $\eta$.
$|h|$ gives a metric on $\Gamma$ as follows:
\begin{eqnarray*}
& & v,w \in V(\Gamma), \\
& & d(h,v,w)=\min\{|h|(\eta): \eta \mbox{ is arc connecting } v
\mbox{ to } w\}.
\end{eqnarray*}
For any $y_i, 1 \le i \le s$, let
\begin{eqnarray*}
& & V(y_i)=V_0(\Gamma) \cup \{y_{i+1}, y_{i+2}, \cdots, y_s\}
\, \, \hbox{and} \\
& & d(h,y_i)=\min\{d(h,v,y_i) : v \in V(y_i)\}.
\end{eqnarray*}

Note: $d(h,y_i)$ is much dependent on the linear order.

\medskip
{\bf Definition (6.1):}\ \
For eeach $i$, $1 \le i \le s$,
\begin{namelist}{(iii)}
\item[{\rm (i)}]
$g(h, y_i) : V(\Gamma) \longrightarrow \setC$ is a
ground function defined by: \\
$g(h, y_i)(v) =\max \{0, d(h,y_i) - d(h,y_i,v)\}$,
for any $v \in V(\Gamma)$.
\item[{\rm (ii)}]
$b(h,y_i)$ (or simply $b_i(h)$) is equal to
$(g(h, y_i),h)$, it is in $D(\Gamma)_h$.
\item[{\rm (iii)}]
$\theta(h):\setC^s \longrightarrow \setC^k$ is a linear map defined by:
$\theta(h)(c_1, c_2, \cdots, c_s)=\Psi(c_1b_1(h)+c_2 b_2(h)+\cdots + c_s b_s(h)
)$.
\end{namelist}

\medskip
{\bf Proposition (6.2):}\ \
Suppose every inner vertices are connecting to base points.
Then
\begin{namelist}{(iii)}
\item[{\rm (i)}]
$(b_1, b_2, \cdots, b_s)$ is a trivialization of the vector bundle
$D(\Gamma)$ over $H(\Gamma)$.
\item[{\rm (ii)}]
$\theta:H(\Gamma) \longrightarrow$ End $(\setC^s, \setC^k)$ sends
$H(\Gamma)$ into a bounded set of \hfill\break End $(\setC^s, \setC^k)$.
\item[{\rm (iii)}]
For each $h \in H(\Gamma), \theta(h)$ is an
injective linear map.
\end{namelist}

\medskip
{\bf Proof:}\ \
\begin{namelist}{(iii)}
\item[{\rm (i)}]
For $j >i$, $y_j \in V(y_i)$.
$d(h, y_i) \le d (h, y_j, y_i)$.
Thus, $g(h, y_i)(y_j)=0$.
For any $v \ne w$ in $V(\Gamma)$, $d(h, v, w)>0$,
no matter that $h(v)$ is equal to
$h(w)$ or not. $d(h, y_i) > 0$, for all
$i=1, 2, \cdots, s$. Thus $\{g(h, y_i), i=1, 2, \cdots, s\}$ are linearly
independent. This proves (i).

\vskip 6pt
Note: If there is an inner vertice which is not connecting to any base point,
the dimension of $D(\Gamma)_h$ is smaller than $s$.
\item[{\rm (ii)}]
$\Psi(b_i(h))=\left(\frac{g(h,y_i)}{h} (e_r)\right)^k_{r=1}$. \hfill\break
For any edge $e=\{v,w\}$, \hfill\break
$|d(h, y_i, v)-d(h, y_i, w)| \le d(h,v,w)=|h(w)-h(v)|$. \hfill\break
Thus,
\begin{eqnarray*}
\left | \frac{g(h, y_i)}{h} (e) \right |
& & =\frac{|g(h,y_i)(v)-g(h,y_i)(w) |} {|h(v)-h(w)|} \\
& & \le \frac{|g(h,y_i,v)-g(h,y_i,w) |} {|h(v)-h(w)|}  \le 1
\end{eqnarray*}
and $|\Psi(b_i(h))| \le \sqrt{k}$ in Euclidean norm.
This prove (ii).
\item[{\rm (iii)}]
Assume $\theta(c_1, c_2, \cdots, c_s)=0$ in $\setC^k$.
Let $g=c_1 gh, y_1)+c_2 g(h, y_2) + \cdots + c_s g(h, y_s)$. \hfill\break
Then $\frac{g}{h}(e)=0$, for all $e \in E(\Gamma)$,
that is $g(v)=g(w)$, for any $e =\{v, w\}$ in $E(\Gamma)$.
By assumption, any inner vertice is connecting to a base point.
Thus $g(v)=0$, for all $v \in V(\Gamma)$.
By (i), $c_1 =c_2=\cdots=c_s =0$.
This proves (iii).
\end{namelist}

\vskip 20pt
\Vjp
\noindent{\Large \bf  6.3 \ Boundary behavior of the trivialization}

\vskip 6pt
Suppose $A$ is a subset of vertices in $\Gamma$.
As in section 5.2, $Q(A(\Gamma))=W_0(A(\Gamma)) \big /$ S.E.T.D. relation,
$Q(\Gamma; A)=W_0(\Gamma \big / A) \times Q(A(\Gamma))$.
To restrict to the case that $Q(\Gamma; A)$ is a codimension $1$
boundary, we assume that (i) $A(\Gamma)$ is connected, or (ii) $A$ consists
of two neighboring base points.
But, the second case is quite trivial, we may consider the first
case only.
Thus, $Q(A(\Gamma))=W_0(A(\Gamma))$ and
$Q(\Gamma; A)=W_0(\Gamma; A) \times W_0(A(\Gamma))$.
Furthermore, $D(\Gamma; A)=D(\Gamma \big / A) \times D(A(\Gamma))$
and $H(\Gamma; A)=H(\Gamma \big / A) \times H(A(\Gamma))$.

\Vjp
\noindent{\large \bf \ 6.3.1 \ Having base points}

\vskip 6pt
Suppose there is at least one base point in $A$.

\vskip 6pt
Then, the set of inner vertices $V_1(\Gamma)$ splits into
disjoint union of $V_1(A(\Gamma))$ and $V_1(\Gamma \big / A)$.
Thus, the linear order on $V_1(\Gamma)$ gives linear
orders on $V_1(A(\Gamma))$ and $V_1(\Gamma \big / A)$.
For convenience, consider the following values for a
height function $h$ on $\Gamma$:
\subitem
$\va(h)=\min \{|h(v)-h(w)|$: for all edge $e =\{v,w\}$ in $\Gamma\}$,
\hfill\break
and $|h|=\max\{|h(v) - h(w)|$: for any two vertice $v,w$ in $\Gamma\}$.

\noindent
$\va (h)$ and $|h|$ are not invariant under T.D. relation, they are
not a function on $H(\Gamma)$.

\vskip 6pt
For any $h_1$ in $H(\Gamma \big / A)$ and $h_2$ in $H(A(\Gamma))$,
let $h_\lambda : V(\Gamma) \longrightarrow \setR$ defined by:
\[
h_\lambda(v)=\left\{\begin{array}{ll}
               h_1(v), & \mbox{if $v$ is not in $A$,}\\ \\
               h_1(a)+\lambda \frac{\va(h_1)}{|h_2|}(h_2(v) - h_2(a_0)),
               , & \mbox{if $v$ is in $A$,} \\
              \end{array}
              \right.
\]
where $a$ denote the new base point corresponding to $A$, and $a_0$ is a fixed
base point in $A$. When $\lambda$ is sufficiently small,
$h_\lambda$ is a height function of $\Gamma$.

\vskip 6pt
Let $b(h_\lambda, y_1)$, $b(h_\lambda, y_2)$, $\cdots$, $b(h_\lambda, y_s)$
denote the basis of $D(\Gamma)_{h_\lambda}$ constructed in section 6.2.
The following results are straightforward.

\medskip
{\bf Proposition (6.3):}\ \
\begin{namelist}{(iii)}
\item[{\rm (i)}] If $y \in V_1 (\Gamma \big / A)$,
\[
\lim\limits_{\lambda \to 0} g(h_\lambda, y)(v)=
\left\{ \begin{array}{ll}
              g(h_1, y)(v), & \mbox{for $v \not \in A$,}\\ \\
              g(h_1,y)(a), &\mbox{for $v \in A$.} \\
        \end{array}
\right.
\]

\item[{\rm (ii)}] If $y \in A \cap V_1(\Gamma)$,
$b(h_\lambda, y) |_{A(\Gamma)}=b(h_2,y)$,
for all $\lambda >0$.
\item[{\rm (iii)}]
If $y \in A \cap V_1 (\Gamma)$,
$\lim\limits_{\lambda \to 0} g(h_\lambda, y)(v)=0$,
for $v \not \in A$.
\end{namelist}

\vskip 6pt
In general, if a parameter of height function $\{h_\lambda, \lambda >0\}$
satisfying:
\begin{namelist}{(ii)}
\item[{\rm (i)}]
$\lim\limits_{\lambda \to 0} h_\lambda = h_1$ and
\item[{\rm (ii)}]
$\lim\limits_{\lambda \to 0} (\varphi(\lambda) h_\lambda |_A + \mu(\lambda))
=h_2$, for some $\varphi(\lambda) > 0$ and $\mu(\lambda)$,
\end{namelist}

\noindent
we also have the result
\begin{namelist}{(ii)}
\item[{\rm (i)}]
in Proposition (6.3) and
\item[{\rm (ii)}']
$\lim\limits_{\lambda \to 0} b(h_\lambda, y) |_{A(\Gamma)} = b(h_2, y)$,
for $y \in A \cap V_1(\Gamma)=V_1(A(\Gamma))$.
\end{namelist}

\Vjp
\noindent{\large \bf \ 6.3.2 \ Inner vertices only}

\vskip 6pt
Suppose $A$ contains inner vertices only, and $y_r$ is the
maximal point in $A$. It is convenience to think the new
inner vertice of $\Gamma \big / A$ corresponding to $A$ is $y_r$.
Thus, $\Gamma \big / A$ has inner vertices $\ov y_r$ and the
points in $V_1(\Gamma)-A$, and $A(\Gamma)$ has inner vertices
the all points in $A$.
($\ov y_r$ is the new inner vertice in $\Gamma \big / A$)
The linear order of $A$ is just the restriction from that
in $V_1(\Gamma)$ and the linear order of $V_1(\Gamma \big / A)$ is
also the restriction to $\{y_r\} \cup (V_1(\Gamma)-A)$
and change the $y_r$ to $\ov y_r$, that is,
for $y_i \in V_1(\Gamma)-A$, $y_i < \ov y_r$, if and only if,
$y_i < y_r$.

\vskip 6pt
Similar to the above case, for $h_1 \in H(\Gamma \big / A)$ and
$h_2 \in H(A(\Gamma))$, let $h_\lambda$ be the height function on
$\Gamma$ defines by:
\[
h_{\lambda}(v)=
\left\{ \begin{array}{ll}
              h_1 (v), & \mbox{if $v \in V(\Gamma)-A$,}\\ \\
              h_1(\ov y_r)+\lambda \frac{\epsilon(h_1)}{|h_2|}
              (h_2(v)-h_2(y_r)),
              &\mbox{if $v \in A$.} \\
        \end{array}
\right.
\]
It is obvious that, when $\lambda$ is sufficiently small
$(\lambda >0)$, $h_\lambda$ is a well-defined height function on $\Gamma$.
Now, by the definition of the basis $(b_i)$, we have
\begin{namelist}{(iii)}
\item[{\rm (i)}]
$\{b(h_1, y)$, for $y \in V_1(\Gamma)-A$, and
$b(h_1, \ov y_r)\}$,
\item[{\rm (ii)}]
$\{b(h_2, y)$, for $y \in A\}$ and
\item[{\rm (iii)}]
$\{b(h_\lambda, y), y \in V_1(\Gamma)\}$.
\end{namelist}

\noindent
The first and the third are bases for $D(\Gamma \big / A)_{h_1}$
and $D(\Gamma)_{h_\lambda}$,
obviously. But, the second is not a basis for
$D(A(\Gamma))_{h_2}$.
Let $s_1 = |A|$, the number of elements in $A$.
By T.D. relation, $D(A(\Gamma))_{h_2}$ has dimension $s_1 -1$.

\vskip 6pt
By assumption $A(\Gamma)$ is connected,
$d(h_2, y, y')$ is a well-defined non-negative number, for
$y, y' \in A$. Thus, $d(h_2, y)= \min \{d(h_2, y, y'): y' > y\}$
is also a well-defined finite value, for $y \ne y_r$.
Therefore, it is reasonable to consider the basis
$\{b(h_2, y): y \in A$ and $y \ne y_r\}$
as the
trivialization of $D(A(\Gamma))$ over $H(A(\Gamma))$.

\medskip
{\bf Proposition (6.4):}\ \
Suppose $A$ is a subset of $V_1(\Gamma)$ such that $A(\Gamma)$
is connected and $y_1 < y_2 < \cdots < y_{s_1}$
is the linear order of elements in $A$ restrict from that
in $V_1(\Gamma)$. Then $\{b(h_2, y_1), b(h_2, y_2), \cdots,
b(h_2, y_{s_1 -1})\}$, for $h_2 \in H(A(\Gamma))$,
is a trivialization for $D(A(\Gamma))$ over $H(A(\Gamma))$.

\medskip
{\bf Proposition (6.5):}\ \
$\{h_\lambda , \lambda >0\}$ defined above.
\begin{namelist}{(iiii)}
\item[{\rm (i)}]
If $y \in V_1(\Gamma \big / A)$ and $y \ne \ov y_r$,
\[
\lim\limits_{\lambda \to 0} g(h_\lambda, y)(v)=
\left\{ \begin{array}{ll}
              g(h_1, y)(v), & \mbox{for $v \not \in A$,}\\ \\
              g(h_1,y)(\ov y_r), &\mbox{for $v \in A$.} \\
        \end{array}
\right.
\]
\item[{\rm (ii)}]
\[
\lim\limits_{\lambda \to 0} g(h_\lambda, y_r)(v)=
\left\{ \begin{array}{ll}
              g(h_1, \ov y_r)(v), & \mbox{for $v \not \in A$,}\\ \\
              g(h_1, \ov y_r)(\ov y_r), &\mbox{for $v \in A$.} \\
        \end{array}
\right.
\]
\item[{\rm (iii)}]
For $y \in A-\{y_r\}$,
\[
b(h_\lambda, y) |_{A(\Gamma)} = b(h_2, y),
\mbox{ for } \lambda  \mbox{ sufficiently small},
\]
 where the equality means that
they are equivalent under T.D. relation.
\item[{\rm (iv)}]
For $y \in A - \{y_r\}$ and $v \not \in A$,
\[
\lim\limits_{\lambda \to 0} h(h_\lambda , y) (v)=0.
\]
\end{namelist}

\medskip
{\bf Proof.}\ \ (i) Assume $y \in V_1(\Gamma \big / A)$ and
$y \ne \ov y_r$.
$d (h_\lambda, y) = \min \{d(h_\lambda, y, v): v \in V_0(\Gamma)$
or $(v \in V_1(\Gamma)$ and $v > y$ in $V_1(\Gamma))\}$,
it is smaller than or equal to
$d' (h_\lambda , y) =\min \{d(h_\lambda, y , v)$ $:$
$v \in V_0(\Gamma)$ or $(v \in \{y_r\} \cup
(V_1(\Gamma) -A)$ and $v > y$ in $V_1(\Gamma))\}$.

\vskip 6pt
Claim: \ $\lim\limits_{\lambda \to 0} d(h_\lambda , y) =
\lim\limits_{\lambda \to 0} d'(h_\lambda, y)$.

\vskip 6pt
(Case 1): \ For all $v \in A$, $v < y$.
$$
d (h_\lambda , y)=d'(h_\lambda, y), \, \, \mbox{ for all } \lambda.
$$

\vskip 6pt
(Case 2): \ There exists $v \in A$, $v > y$.

\vskip 6pt
Then $y_r > y$. When $\lambda$ approaches to $0$,
$d(h_\lambda, y, v)$ and $d(h_\lambda, y, y_r)$
have the same limist $(d(h_1, y, \ov y_r))$.
Therefore, although $d'(h_\lambda, y)$ has less
minimum set, it also has the same limist as
$d(h_\lambda, y)$.

\noindent
And, it is obvious that $d'(h_\lambda, y)=d(h_1, y)$.
This implies the equality in (i) easily.

\vskip 6pt
(ii) \ The set $V(y_r) = V_0(\Gamma) \cup \{v \in V_1(\Gamma),
v > y_r\}$ is equal to $V(\ov y_r)=V_0(\Gamma \big / A) \cup
\{v \in V_1 (\Gamma \big / A)$, $v > y_r\}$.
Thus, $\lim\limits_{\lambda \to 0} d(h_\lambda, y_r)=d(h_1, \ov y_r)$.
This implies (ii).

\vskip 6pt
(iii): \ For $y \in A$, $y \ne y_r$,
\begin{eqnarray*}
& & V(y)=V_0(\Gamma) \cup \{v \in V_1(\Gamma): v>y\}
\mbox{ for } \Gamma,  \\
\mbox{and } \quad & & V'(y)=\{v \in A, v > y \}, \mbox{ for } A(\Gamma).
\end{eqnarray*}

\noindent
Now, compare $d(h_\lambda, y)$ with $d(h_\lambda |_A, y)$ $:$
\begin{eqnarray*}
d(h_\lambda, y) = \min \{ d(h_\lambda, y, v): v \in V(y)\} \\
d(h_\lambda |_A, y) = \min \{ d(h_\lambda, y, v): v \in V'(y)\}
\end{eqnarray*}

\noindent
$y_r$ is in $V(y)$, $y_r$ is also in $V'(y)$.
For any $v \in V(y) - V'(y)$, that is, $v \not\in A$ and
$v \in V(y)$, $d(h_\lambda, y, v)>d(h_\lambda, y, y_r)$,
for $\lambda$ sufficiently small. Thus, when $\lambda$
is sufficiently small,
$d(h_\lambda, y)= d(h_\lambda |_A, y)$.
But, $d(h_\lambda |_A, y) = \lambda d(h_2, y)$.
Thus,
\[
(g(h_\lambda , y)(v), h_\lambda(v))
=\lambda (g(h_2, y)(v), h_2(v))
+ (0, h_1(\ov y_r) -\lambda h_2(y_r)),
\]
for all $v \in A$, that is,
$b(h_\lambda , y) |_{A(\Gamma)}$
is equivalent to $b(h_2, y)$ under T.D. relation.

\medskip
{\bf Remark (6.6):}\ \
For $y \in A-\{y_r\}$,
$g(h_\lambda, y)(v)=0$,
for $v \in V(\Gamma)-A$ and $\lambda$ sufficiently small.

\medskip
{\bf Proof:}\ \ $d(h_\lambda, y, v) > d(h_\lambda, y) = \lambda d(h_2, y)$.
Thus, $g(h_\lambda, y)(v)=\max \{\,0,$ $ d(h_\lambda, y) -
d(h_\lambda, y, v)\}=0$.

\Vjp
\noindent{\large \bf \ 6.3.3 \ (Continued)}

\vskip 6pt
$\{h_\lambda, \lambda >0\}$, as above, is a parameter
of height function on $\Gamma$ such that
$\lim\limits_{\lambda >0} h_\lambda = h_1$ in
$H(\Gamma \big / A)$ and
$\lim\limits_{\lambda \to 0}(h_\lambda |_{A(\Gamma)})=h_2$
(in the sense that they are equivalence under T.D. relation)
in $H(A(\Gamma))$.
Part of limit behavior of $b(h_\lambda, y)$,
$y \in V_1(\Gamma)$, is already studied. But, there are still
something, we forget to study in section 6.3.1 and
6.3.2, that is, when $y \in V_1(\Gamma) - A$, the limit of
the restriction  of $b(h_\lambda, y)$ to $A(\Gamma)$.
Actually, the limit does not exist. Thus, we should
modify $b(h_\lambda, y)$ such that its restriction to $A(\Gamma)$
will approach to zero.

\vskip 6pt
An easy approach is the following:

\vskip 6pt
For $y \in V_1(\Gamma) - A$, let $\ov g (h_\lambda, y)$ be a
ground function such that, for any $e=\{v, w\}$,
\begin{eqnarray*}
& & \ov g(h_\lambda, y)(v) - \ov g(h_\lambda , y)(w) \\ \\
& & = \left\{
  \begin{array}{ll}
     g(h_\lambda, y)(v) - g(h_\lambda, y)(w),
     & \mbox{for  $e \not \in E(A(\Gamma)),$} \\ \\
     \lambda(g(h_\lambda, y)(v) - g(h_\lambda, y)(w)),
     &\mbox{for $e \in  E(A(\Gamma)).$}
  \end{array}
  \right.
\end{eqnarray*}
(note: $|g(h_\lambda, y)(v)- g(h_\lambda, y)(w) | \le | h | (e)$)

\vskip 6pt
When $\lambda$ is small, we switch $g(h_\lambda, y)$ to
$\ov g (h_\lambda, y)$,
continuously, and we get a new basis
$\{\ov b (h_\lambda, y), y \in A \cap V_1 (\Gamma)\}$.
Then, Proposition (6.3) and (6.5) also hold for this new
basis. For example:
\[
\lim\limits_{\lambda \to 0} \ov g(h_\lambda, y)(v)= g(h_1, y)(v),
\]
for $y \in V_1(\Gamma \big / A)$ and $v \not\in A$, in (6.3).
Moreover, we also have

\medskip
{\bf Proposition (6.7):}\ \
\begin{namelist}{(ii)}
\item[{\rm (i)}]
For any $y \in V_1(\Gamma)-A$,
\[
\lim\limits_{\lambda \to 0} \ov b(h_\lambda, y) |_{A(\Gamma)}=0.
\]
\item[{\rm (ii)}]
If $A \subset V_1(\Gamma)$ and $y_r$ is the maximal
element in $A$, then
\[
\lim\limits_{\lambda \to 0} b(h_\lambda, y_r) |_{A(\Gamma)}
=b(h_2, y_r).
\]
\end{namelist}

\medskip
{\bf Remark:}\ \ In the case (ii) of the above proposition, we
may define $\ov b(h_\lambda, y_r)$ similarly such that
$\lim\limits_{\lambda \to 0} \ov b(h_\lambda, y_r)|_{A(\Gamma)}=0$
as in (i).

\vskip 6pt
Thus, by (6.3), (6.5) and (6.7), we find that the limit
of the new basis $\{b(h_\lambda, y)\}$ is
completely determined by the two bases at the limit
$(h_1, h_2)$ of $\{h_\lambda\}$,
in $H(\Gamma \big / A) \times H(A(\Gamma))$.

\vskip 6pt
The above modification of basis is done on a collaring neighborhood
of $H(\Gamma; A)$ in $H(\Gamma)$.
Thus, we need a non-overlapping collaring neighborhood for each
codimension-1 boundary $H(\Gamma ; A)$.
We can believe that the space is good enough to do so.
If one does not believe it, one could try an alternating method:
a more general modification on all $g(h,y)$ in the
following way:
\begin{namelist}{(iii)}
\item[{\rm (i)}]
For any two functions $g_1 , g_2 : V(\Gamma) \longrightarrow \setC$
and $e =\{v, w\}$,
$\frac{g_1}{g_2} (e) =\frac{g_1(v)-g_1(w)}{g_2(v)-g_2(w)}$.
\item[{\rm (ii)}]
$E(y)=\{e =\{y, v\}$, $e$ is an edge in $\Gamma\}$.
\item[{\rm (iii)}]
$\ov g(h, y):V(\Gamma) \longrightarrow \setC$ is a ground function on
$\Gamma$ satisfying
\[
\frac{\ov g(h, y)}{g(h,y)}(e)=
\left\{
\begin{array}{ll}
\frac{5 |h| (e)}{\delta(y)}, &\mbox{if
                              $|h|(e) \le \frac{1}{5} \delta (y)$, }\\ \\
1, &\mbox{if otherwise,}
\end{array}
\right.
\]
where $|h|(e)=|h(v) - h(w)|$ $(e=\{v,w\})$
and $\delta (y) =\min\{|h|(e): e \in E(y)\}$.
\end{namelist}

\noindent
(That is, we shrink the effect of $| h |(e)$ on
$g(h, y)$, when $|h |(e)$ is much smaller than
$\delta(y)$, the $h$-lengths of edges connecting to $y$.)

\vskip 6pt
Note: $\frac{\ov g}{g}(e)=1$ means that
$\ov g(v) -\ov g(w)= g(v)-g(w)$,
whenever $g(v)=g(w)$ or not.

\vskip 6pt
Using the modified basis
$\{b(h,y) = (\ov g(h, y), h)\}$,
we also have the result (6.3), (6.5) and (6.7), except the result
(ii) in (6.5) becoming that
\[
{\rm (ii)}' \, \, \,
\lim\limits_{\lambda \to 0} \ov g(h_\lambda, y_r)(v)
 =\left\{
  \begin{array}{ll}
  g(h_1, \ov y_r)(v), \mbox{for $v \not\in A$} \\   \\
  g(h_1, \ov y_r)(\ov y_r), \mbox{for $v \in A$} \\
  \end{array}
  \right.
  \]
(its limit returns to the original one).

\vskip 20pt
\Vjp
\noindent{\Large \bf  6.4 \ Isotopy of trivializations}

\vskip 6pt
Changing the linear order of inner vertices, we get different trivializations
of $D(\Gamma)$. We shall show that they are all isotopic.

\vskip 6pt
As in section 6.2, $y_1 < y_2 < \cdots <y_s$ are the inner vertices
with linear order (in $\Gamma$).
$V(y_i)=V_0(\Gamma) \cup \{y_{i+1}, y_{i+2}, \cdots, y_s\}$,
$i=1, 2, \cdots, s$, and
$d(h, y_i)=\min\{d(h, v, y_i): v \in V(y_i)\}$.

Now, assume $r$ is an integer, $1 \le r \le s-1$, and
consider the linear order
$y_1 < y_2 < \cdots < y_{r-1}<y_{r+1} < y_r <y_{r+2} < \cdots < y_s$,
which switches the positions of $y_r$ and $y_{r+1}$.
Then, we should consider
\begin{eqnarray*}
& & V'(y_r) = V_0(\Gamma) \cup \{y_{r+2}, y_{r+3}, \cdots, y_s\}, \\
& & V'(y_{r+1}) = V_0(\Gamma) \cup \{y_r, y_{r+2}, y_{r+3}, \cdots, y_s\}, \\
& & d'(h, y_r)=\min\{d(h, v, y_r): v \in V'(y_r)\}, \\
& & d'(h, y_{r+1})=\min\{d(h,v,y_{r+1}): v \in V'(y_{r+1})\}.
\end{eqnarray*}
for convenience, let $V'(y_i)= V(y_i)$ and
$d'(h,y_i)=d(h, y_i)$, for $i \ne r$ and $r+1$.
Furthermore, let $g'(h, y_i)$ be the ground function defined by
$g'(h, y_i)(v)=\max \{0, d'(h, y_i) - d(h, v, y_i)\}$.
And $b'(h, y_i)=(g'(h, y_i), h)$,
$i=1, 2, \cdots, s$. Then,
$\{ b'(h, y_i)$, $i=1,2, \cdots, s\}$ is also a trivialization on
$D(\Gamma)$ over $H(\Gamma)$.
We shall show that $\{b(h, y_i)\}$
and $\{b'(h, y_i)\}$ can be homotopic to a trivialization
$\{\ov b (h, y_i)\}$ through parameters of trivialization
$\{b_t(h, y_i)\}$ and $\{b'_t(h, y_i)\}$, respectively.

\vskip 6pt
Consider the following, as above:
\begin{namelist}{(iiiii)}
\item[{\rm (i)}]
$\ov V(y_r) = \ov V(y_{r+1}) = V'(y_r) \cap
V'(y_{r+1})$, that is, $V'(y_r)$.
\item[{\rm (ii)}]
$\ov V(y_i) = V(y_i)$, for $i \ne r$, $r+1$.
\item[{\rm (iii)}]
$\ov d(h, y_i)=\min\{d(h, v, y_i): v \in \ov V(y_i)\}$.
\item[{\rm (iv)}]
$\ov g(h, y_i)(v) = \max \{0, \ov d(h, y_i) - d(h, v, y_i)\}$
\item[{\rm (v)}]
$\ov b (h, y_i)=(\ov g (h, y_i), h)$,
$i=1, 2, \cdots, s$.
\end{namelist}

\noindent
It is easy to see that
$\ov d(h, y_i)\ge d(h, y_i)$ and
$\ov d(h, y_i) \ge d'(h, y_i)$, $i=1, 2, \cdots, s$. \hfill\break
Let $d_t(h, y_i) = td(h, y_i) + (1 -t) \ov d (h, y_i)$,
for $0 \le t \le 1$. \hfill\break
Also, $g_t(h, y_i)(v)=\max\{0, d_t(h, y_i) - d(h, v, y_i)\}$
and $b_t(h, y_i)=(g_t(h, y_i), h)$, $i=1, 2, \cdots, s$.
Similarly, for $d'_t (h, y_i)$,
$g'_t(h, y_i)$ and $b'_t (h, y_i)$.

\medskip
{\bf Proposition (6.8):}\ \ For any $t$, $0 \le t \le 1$,
and $h \in H(\Gamma)$, both
$\{b_t(h, y_i)\}^s_{i=1}$
and $\{ b'_t (h, y_i)\}_{i=1}^s$ are bases for $D(\Gamma)_h$.
Thus, the two trivialization $\{b(h, y_i)\}$
and $\{b'(h, y_i)\}$ are isotopic through the parameter of
trivializations
$\{b_t(h, y_i), i=1, 2, \cdots, s\}$ and
$\{b'_t(h, y_i), i=1, 2, \cdots, s\}$.

\medskip
{\bf Proof:}\ \ It is easy to show that, for $0 \le t \le 1$,
$g_t(h, y_i)$, $i=1, 2, \cdots, s$, are linearly independent
(Similarly, for $\{g'_t(h, y_i)\}$).
Because $\ov d(h, y_i) \ge d_t(h, y_i)$,
\begin{namelist}{(iiii)}
\item[{\rm (i)}]
$g_t(h, y_i)(y_j) =0$, for $i \le r-1$ and $j > i$,
\item[{\rm (ii)}]
$g_t(h, y_r)(y_j)=0$, for $j \ge r+2$,
\item[{\rm (iii)}]
$g_t(h, y_{r+1})(y_j)=0$, for $j \ge r+2$,
\item[{\rm (iv)}]
$g_t(h, y_i)(y_j)=0$, for $i \ge r+2$ and $j > i$.
\end{namelist}

\noindent
Moreover, $g_t(h, y_i)(y_i) > g_t(h, y_i)(y_j) \ge 0$,
for any $1 \le j \ne i \le s$.

\vskip 6pt
Consider the matrix $G = (G_{i,j})_{1 \le i \le s \atop 1 \le j \le s}$,
$G_{i,j}=g_t(h, y_i)(y_j)$.
The above conditions implies that the determinant of $G$ is
greater than $0$. In fact,
\begin{eqnarray*}
 \det (G) & & = G_{11} G_{22} \cdots G_{ss}
-G_{11} G_{22} \cdots G_{r-1, r-1} G_{r, r+1} G_{r+1, r}
G_{r+2, r+2} \cdots G_{s,s} \\
& & = G_{11} G_{22} \cdots G_{ss} \cdot  (G_{r,r} G_{r+1, r+1})^{-1}
\cdot (G_{r,r} \cdot G_{r+1, r+1} - G_{r, r+1} \cdot G_{r+1,r} ).
\end{eqnarray*}
$G_{r,r} > G_{r, r+1}$ and $G_{r+1, r+1} > G_{r+1, r}$.
Thus, $\det (G) >0$. This proves this proposition.

\vskip 6pt
Therefore, there is no real difference to choose any linear order
for $V_1(\Gamma)$ to construct the trivialization.

\vskip 20pt
\Vjp
\noindent{\Large \bf  6.5 \ Transition maps of vector bundle $\DD(n)$}

\vskip 6pt
Because we have trivializations $(b_1, b_2, \cdots, b_s)$ for
every $D(\Gamma)$, the transition maps are given by all the
identifications stated in section 1.3, or in section 5.3
(vector bundle form).

\Vjp
\noindent{\large \bf  6.5.1 \ Identification of type 0}

\vskip 6pt
This is already absorbed by the special extended translation
and dilation relation
(see section 5.2).

\Vjp
\noindent{\large \bf  6.5.2 \ Identification of type I}

\vskip 6pt
Suppose $A$ is a subset of $V(\Gamma)$ and $A(\Gamma)$
has a univalent inner vertice $v$.

\vskip 6pt
Assume $e =\{v, v_1\}$ is the unique edge in $A(\Gamma)$,
containing $v$ as the endpoint.

\vskip 6pt
By the result of section 6.4, we may assume that $v$ is the minimum
element in $V_1(A(\Gamma))$.
For convenience, let $y_1 = v < y_2 < \cdots <y_r$
denote the inner vertices in $A(\Gamma)$. Then,
for any height function $h$ on $A(\Gamma)$,
$g(h, y_1)(y_1)=|h|(e)$, that is,
$|h(v) - h(v_1)|$, and $g(h, y_1)(y_j)=0$, for $ j\ne 1$.

\medskip
{\bf (6.9):}\ \ $\tau'_1 : D(A(\Gamma)) \longrightarrow D(A(\Gamma))$,
defined in the method (ii) of section 5.3.4,
satisfies the following:
\begin{namelist}{(ii)}
\item[{\rm (i)}]
$\tau'_1 (g(h, y_1), h)$ is exactly equal to
$(g(\tau'_1(h), y_1), \tau'_1(h))$.
\item[{\rm (ii)}]
For $i \ge 2$ and $j \ge 2$,
$\tau'_1 (g(h, y_i))(y_j) = g(h, y_i)(y_j)
=g(\tau'_1(h), y_i)(y_j)$.
\end{namelist}

\medskip
{\bf Proof of (6.9):}\ \
$\tau'_1 (h) (w) = h(w)$, for $w \ne y_1 =v$.
Thus,
\[
g(\tau'_1(h), y_1)(y_1)
=|\tau'_1(h)(v) - \tau'_1(h)(v_1) | = 2 ||h_1||
\]
and
$g(\tau'_1(h), y_1)(w)=0$,  for $w \ne y_1$.
\begin{eqnarray*}
& & \tau'_1 (g(h, y_1))(y_1) \\
& & = g(h, y_1)(v_1) + 2 ||h_1||
\frac{g(h, y_1)(v)-g(h, y_1)(v_1)}{|h(v) - h(v_1)|}  \\
& & = 2 || h_1 ||
\end{eqnarray*}

and $\tau'_1(g(h, y_1))(w)= g(h, y_1)(w)=0$,
for $w \ne y_1$.
This proves (i).

\vskip 6pt
The proof of (ii) is straightforward.

\vskip 6pt
Thus, the restriction of
$\tau'_1$ to $D(A(\Gamma))_h$,
$:$ $D(A(\Gamma))_h \longrightarrow D(A(\Gamma))_{\tau'_1(h)}$,
send $b_1$ to $b_1$ and $b_i$ to $b_i+ \rho_i b_1$,
for $i \ge 2$, where $\rho_i$ is a real number depending on $h$.
By the homotopy property of vector bundle, we may
change $\tau'_1$ to a homotopy one $\tau''_1$,
$\tau''_1(b_i) = b_i$ for all $i$.

\vskip 6pt
Therefore, the transition map for Identification of type I
is the identity map.

\Vjp
\noindent{\large \bf  6.5.3 \ Identification of type II}

\vskip 6pt
Suppose $A$ is a subset of $V(\Gamma)$ and $A(\Gamma)$
has a bivalent inner vertice $v$.

\vskip 6pt
Assume $e_1 =\{v, w_1\}$ and $e_2 = \{v, w_2\}$ are the two
edges connecting to $v$.

\vskip 6pt
As above, we may assume that $v$ is the minimum element in
$V_1(A(\Gamma))$ and $y_1 = v < y_2 < \cdots < y_r$ are
the inner vertices in $A(\Gamma)$.
Then, for any height function of $A(\Gamma)$,
\[
g(h, y_1)(y_1) =\min\{|h| (e_1), |h|(e_2)\}
\]
and $g(h, y_1)(y_i)=0$, $i >1$.

\vskip 6pt
As the notations in section 5.3.1, we have

\medskip
{\bf (6.10):}\ \
\begin{namelist}{(ii)}
\item[{\rm (i)}]
$\tau'_2 (g(h, y_1))=-g(\tau'_2(h), y_1)$.
\item[{\rm (ii)}]
For $i \ge 2$ and $j \ge 2$,
\begin{eqnarray*}
\tau'_2 (g(h, y_i))(y_j) & & = g(h, y_i)(y_j) \\
& & = g(\tau'_2 (h), y_i) (y_j).
\end{eqnarray*}
\end{namelist}

\medskip
{\bf Proof of (6.10):} \ \
\begin{eqnarray*}
& \tau'_2 (g(h, y_1))(y_1) & = g(h, y_1)(w_1) + g(h, y_1)(w_2) -g(h, y_1) (v)\\
& & = -g (h, y_1)(y_1) \\
& & =-\min\{| h |(e_2), | h|(e_1)\}. \\
& g(\tau'_2(h), y_1)(y_1) & =\min\big\{|\tau'_2 (h) | (e_1),
|\tau'_2 (h) | (e_2) \big\} \\
& &= \min \{|h| (e_2), |h|(e_1)\}.
\end{eqnarray*}
This proves (i).

\vskip 6pt
To prove (ii), choose arcs $\eta_1$ and $\eta_2$ from
$y_i$ to $w_1$ and $w_2$, respectively, such that
$\eta_1$ minimizes the $h$-length from $y_i$ to $w_1$ and
$\eta_2$ minimizes the $h$-length from $y_i$ to $w_2$.

\vskip 6pt
(Case 1) \ $y_1$ is on $\eta_1$.

\vskip 6pt
Then $w_2$ is also on $\eta_1$ and $\eta_2$ is a subarc of $\eta_1$.\\
Thus, $g(h, y_i)(w_2) \ge g(h, y_i)(y_1) \ge g(h, y_i) (w_1)$.\\
If $d(h, y_i) \le d(h, w_2, y_i)$,
then $g(h, y_i)(w_2) = g(h_1, y_i)(y_1) =g(h, y_i)(w_1)=0$.\\
If $d(h, w_2, y_i) \le d(h, y_i) \le d(h, w_1, y_i)$, then\\
$g(h, y_i)(w_2) = d(h, y_i)-d(h, w_2, y_i)$
and $g(h, y_i)(w_1)=0$.\\
If $d(h, w_1, y_i) \le d(h, y_i)$, then

$g(h, y_i)(w_1) = d(h, y_i) - d(h, w_1, y_i)$ and
$g(h, y_i)(w_2) = d(h, y_i) - d(h, w_2, y_i)$.\\
By assumption, $i \ge 2$, $y_i > y_1$, we have
$d(h, y_i)= d(\tau'_2(h), y_i)$.\\
Thus,$g(h, y_i)(w_1) = g(\tau'_2 (h), y_i)(w_1)$ and
$g(h, y_i)(w_2) = g(\tau'_2 (h), y_i)(w_2)$.

\vskip 6pt
(Case 2) \ \ $y_1$ is on $\eta_2$.

\vskip 6pt
Then $w_1$ is also $\eta_2$ and $\eta_1$ is a subarc of $\eta_2$.
And we can prove that
\[
g(h, y_i)(w_j) = g(\tau'_2(h), y_i)(w_j), \, \, j =1, 2,
\]
as in (case 1).

\vskip 6pt
(Case 3)\ \ $y_1$ is not on $\eta_1 \cup \eta_2$.

\vskip 6pt
Then
\begin{eqnarray*}
& & d(h, y_i, w_1) = d(\tau'_2 (h), y_i, w_1)\\
 \mbox{and} \, \,  & & d(h, y_i, w_2) = d(\tau'_2 (h), y_i, w_2) .
\end{eqnarray*}
Thus, we also have the equalities
\[
g(h, y_i)(w_j) = g(\tau'_2(h), y_i)(w_j), j=1, 2.
\]
The proof of other results is similar and is omitted.

\vskip 6pt
What is the difference of
$\tau'_2 g(h, y_i))(y_1)$ and
$g(\tau'_2(h), y_i)(y_1)$, for $i \ge 2$?

\medskip
{\bf (6.11):}\ \ For $i \ge 2$,
\[
(\tau'_2 (g((h, y_i))(y_1) - g(\tau'_2 (h), y_i)(y_1) |
 \le 2 \min \{ |h| (e_1), |h|(e_2)\}.
\]

\medskip
{\bf Proof:}\\
$ \left | \tau'_2 (g((h, y_i))(y_1) - \tau'_2(g(h, y_i)(w_1)\right |
 = \left | g(h, y_i)(w_2) -g(h, y_i)(y_1) \right | \le | h | (e_2)$. \\
$\left | g(\tau'_2(h), y_i) (y_1) - g(\tau'_2(h), y_i)(w_1) \right |
 \le | \tau'_2 (h) | (e_1) = | h | (e_2)$. \\
By (6.10), $\tau'_2 (g(h, y_i)))(w_1) =
g(\tau'_2 (h), y_i)(w_1)$.

\vskip 6pt
Thus,
$\left | \tau'_2 (g(h, y_i))(y_1) - g(\tau'_2(h), y_i)(y_1) \right |
\le 2 | h |(e_2).$ \hfill\break
Similarly,
\[
\left | \tau'_2 (g(h, y_i))(y_1) - g(\tau'_2(h), y_i)(y_1) \right |
\le 2 | h | (e_1).
\]
This proves (6.11).

\vskip 6pt
Now, consider the restriction of $\tau'_2$ to
$D(A(\Gamma))_h$,
$:$ $D(A(\Gamma))_h \longrightarrow D(A(\Gamma))_{\tau'_2(h)}$,
it sends $b(h, y_1)$ to $-b(\tau'_2(h), y_1)$
and sends $b(h, y_i)$ to
$b(\tau'_2(h), y_i) + \rho'_i b(\tau'_2 (h), y_1)$,
where $\rho'_i$ is a real number depending on $h$ and
$| \rho'_i | \le 2$.

\vskip 6pt
Thus, we may assume that the identification map
$\tau_2$ has the following form:
$\tau_2 (b_1)=-b_1$ and $\tau_2(b_i) = b_i$, for $i >1$,
without changing the isomorphism class of vector bundle $\DD(n)$.

\Vjp
\noindent{\large \bf  6.5.4 \ Identification of type III}

\vskip 6pt
Suppose $A = \{v, w\}$ is an edge of $\Gamma$ and $v$ is an inner vertice.
Assume $v$ is the smallest inner vertice in $A$
(if $w$ is a base point, then the assumption holds
automatically).

\vskip 6pt
As the notations in section 5.3.2,
\[
\tau_3(b(h, v)) =
\frac{|h(v) - h(w) |}{h(v)-h(w)} = \pm 1 \, \, \mbox{in} \, \, \setC,
\]
for any height function $h$ on $A(\Gamma)$.

\vskip 6pt
Let $\{b(\cdot , y_i), y_i \ne v\}$ denote the
trivialization for $D(\Gamma \big / A)$
and $\delta_1$ the standard basis $\{1\}$ in $\setC$.
Then $\tau_3 : D(\Gamma; A) \longrightarrow D(\Gamma \big / A) \times \setC$
has the following form:
\begin{eqnarray*}
& \tau_3(b(\cdot, y_i)) & = b(\cdot, y_i), \, \, \mbox{for}\, \, y_i \ne v, \\
& \tau_3(b(\cdot, y_i)) & = \pm \delta_1 .
\end{eqnarray*}

\Vjp
\noindent{\large \bf  6.5.5 \ Identification of type IV}

\vskip 6pt
In this case, $A=\{\,v,w\,\}$ is not an edge in $\Gamma$,
the identification map is nothing but the identity map.

\Vjp
\noindent{\large \bf  6.5.6 \ Conclusion}

\vskip 6pt
Without changing the isomorphism class of $\DD(n)$,
the transition maps for the
trivialization $(b_1, b_2, \cdots, b_s, \delta_1, \cdots, \delta_{3n-k})$
in $D(\Gamma) \times \setC^{3n-k} \times
\Sigma_{3n} \cdot \Psi$ has only two kinds:
\begin{namelist}{(ii)}
\item[{\rm (i)}]
one is the permutation of basis.
\item[{\rm (ii)}]
One is the composition of permutation and $\theta_i$,
where $\theta_i(b_i) = -b_i$ and
$\theta_i(b_j)=b_j$, for $j \ne i$.
\end{namelist}
Thus, $\DD(n)$ has finite structure group.

\newpage

\centerline{\Large \bf References}

\begin{namelist}{[9]}
\item{[1]}
D. Altschuler and L. Freidel, {\it Vassiliev knot invariants and
Chern-Simons perturbation theory to all order}, 1995, Preprint.

\item{[2]} S. Axelrod and I. Singer, {\it Chern-Simons perturbation theory},
in the proceedings of the XXth International conference on
differential geometric methods in theoretical physics, June 3-7, 1991, New York
City (World Scientific, Singapore, 1992). \hfill\break
S. Axelrod and I. Singer, {\it Chern-Simons perturbation theory II},
Jour. Diff. Geom. {\bf 39} (1994) 173.

\item{[3]} D. Bar-Natan, {\it Perturbative Aspects of the Chern-Simons
Topological Quantum Field Theory}, Ph. D. thesis,
Princeton University, June 1991. \hfill\break
D. Bar-Natan, {\it Perturbative Chern-Simons theory},
Journal of Knot Theory
and its Ramifications {\bf 4} (4) (1995), 503-547.

\item{[4]}
D. Bar-Natan, {\it On the Vassiliev knot invariants}, Topology {\bf 34}
(1995) 423.

\item{[5]}
R. Bott and C. Taubes,
{\it On the self-linking of knots}, Jour. Math. Phys. {\bf 35} (10)
(1994), 5247-87.

\item{[6]}
W. Fulton and R. MacPherson, {\it Compactification of
Configuration spaces}, Ann. Math. 139 (1994), 183-225.

\item{[7]}
E. Guadagnini, M. Martellini and M. Mintchev, Nucl. Phys.
{\bf B330} (1990) 575.

\item{[8]}
M. Kontsevich, {\it Vassiliev's knot invariants}, Adv. in Sov. Math.,
{\bf 16(2)} (1993), 137.

\item{[9]}
V. A. Vassiliev,
{\it Cohomology of knot spaces},
in ``Theory of singularities and its applications''
(ed. V. I. Arnold), Advances in Soviet Mathematics, AMS, 1990.
\end{namelist}

\end{document}